\pdfoutput=1
\RequirePackage{ifpdf}
\ifpdf 
\documentclass[pdftex]{sigma}
\else
\documentclass{sigma}
\fi

\numberwithin{equation}{section}
\newtheorem{Theorem}{Theorem}[section]
\newtheorem{Corollary}[Theorem]{Corollary}
\newtheorem{Lemma}[Theorem]{Lemma}
\newtheorem{Proposition}[Theorem]{Proposition}
{\theoremstyle{definition}
\newtheorem{Definition}[Theorem]{Definition}
\newtheorem{Example}[Theorem]{Example}
}

\DeclareMathOperator{\Reoperator}{Re}
\DeclareMathOperator{\Lip}{Lip}
\begin{document}

\newcommand{\arXivNumber}{1401.4622}

\allowdisplaybreaks

\renewcommand{\thefootnote}{$\star$}

\renewcommand{\PaperNumber}{064}

\FirstPageHeading

\ShortArticleName{Non-Commutative Resistance Networks}

\ArticleName{Non-Commutative Resistance Networks\footnote{This paper is a~contribution to the Special Issue on
Noncommutative Geometry and Quantum Groups in honor of Marc A.~Rief\/fel.
The full collection is available at
\href{http://www.emis.de/journals/SIGMA/Rieffel.html}{http://www.emis.de/journals/SIGMA/Rieffel.html}}}

\Author{Marc A.~RIEFFEL}

\AuthorNameForHeading{M.A.~Rief\/fel}

\Address{Department of Mathematics, University of California, Berkeley, CA 94720-3840, USA}
\Email{\href{mailto:rieffel@math.berkeley.edu}{rieffel@math.berkeley.edu}}
\URLaddress{\url{http://math.berkeley.edu/~rieffel/}}

\ArticleDates{Received January 22, 2014, in f\/inal form June 10, 2014; Published online June 14, 2014}

\Abstract{In the setting of f\/inite-dimensional $C^*$-algebras~${\mathcal A}$ we def\/ine what we call a~Riemannian metric
for~${\mathcal A}$, which when~${\mathcal A}$ is commutative is very closely related to a~f\/inite resistance network.
We explore the relationship with Dirichlet forms and corresponding seminorms that are Markov and Leibniz, with
corresponding matricial structure and metric on the state space.
We also examine associated Laplace and Dirac operators, quotient energy seminorms, resistance distance, and the
relationship with standard deviation.}


\Classification{46L87; 46L57; 58B34}

\renewcommand{\thefootnote}{\arabic{footnote}}
\setcounter{footnote}{0}

\tableofcontents

\section{Introduction}

This paper has its origins in three questions that arose at dif\/ferent times during my research concerning quantum metric
spaces.
The f\/irst of these questions was my puzzlement about how the ``resistance distance'' that is def\/ined on resistance
networks f\/its in with the metrics that arise on the state spaces of quantum metric spaces
(see Section~12 of~\cite{R5}).
The second question concerns what conditions ensure that quotients of Leibniz seminorms are again Leibniz, a~property
that is important in dealing with quantum metric spaces (see~\cite{R21}).
More recently my research led me to examine the Leibniz property for standard deviation (see~\cite{R28}).
This eventually led me to ask whether there was a~relationship between that topic too and quantum metric spaces and
resistance networks.

In this paper I seek to give a~coherent account of how closely these questions are related, and of the answers to them
that I have found.
I generally carry out the discussion in the setting of non-commutative $C^*$-algebras.
In order not to be distracted by all the technicalities encountered when dealing with unbounded operators and their
dense domains, I deal in this paper only with f\/inite-dimensional $C^*$-algebras, somewhat in the spirit of the seminal
paper of Beurling and Deny~\cite{BgD}.
There is plenty to be said just about the purely algebraic aspects.

As a~thread to tie things together I introduce a~structure that I call a~``non-commutative Riemannian metric''.
This structure lies just below the surface of some of the literature concerning quantum dynamical semigroups and
Dirichlet forms~\cite{CpS, Svg1, Svg2}, but I have not seen this structure explicitly mentioned there.
We will see (Section~\ref{comm}) that when the underlying $C^*$-algebra is commutative, a~Riemanian metric for it leads
naturally to a~resistance network.

In order to provide a~coherent narrative, I include much material that already appears in the literature.
Thus many parts of this paper can be considered to be expository.
But even in these parts many small novelties are included.
And perhaps this paper can serve as a~useful guide for those who are beginning to learn about quantum dynamical
semigroups and Dirichlet forms.

I expect that most of the new results in this paper have suitable extensions to the setting of inf\/inite-dimensional
$C^*$-algebras, for which one will need to work with unbounded operators (mostly derivations).
Many new phenomena will then arise.
There is a~very large literature containing many techniques for dealing with that setting
(see~\cite{Cpr2} and the references it contains). But I do not plan to carry out some of these extensions myself, unless
I happen to f\/ind later that they are important for my study of quantum metric spaces.
My impression at this point is that the setting of this paper is too favorable to be applicable to the main issues that
I am exploring concerning quantum metric spaces.

In f\/inding the path taken in this paper I have been strongly inf\/luenced by the work of Sauvageot, especially Section~3
of~\cite{Svg2}.
It was Sauvageot who discovered the dif\/ferential calculus associated to Dirichlet forms.
Other important sources for me have been~\cite{AHK, DvL}.

\section{Dif\/ferential calculi with quasi-correspondences}

In this section we develop the aspects of non-commutative Riemannian metrics that do not depend on the positivity of the
${\mathcal A}$-valued inner product.
So we assume here only that~${\mathcal A}$ is a~f\/inite-dimensional unital $*$-algebra over ${\mathbb C}$.
In Section~\ref{ncrm} we will assume that~${\mathcal A}$ is a~f\/inite-dimensional $C^*$-algebra, so that positivity has meaning.
Finite-dimensionality is needed in only a~few crucial places, and we will usually point out these places.
We recall~\cite{GVF} that by a~f\/irst-order dif\/ferential calculus over~${\mathcal A}$ one means a~pair $(\Omega,
{\partial})$ consisting of an~${\mathcal A}$-bimodule $\Omega$ (thought of as an analog of a~space of dif\/ferential
one-forms) and a~derivation ${\partial}$ from~${\mathcal A}$ into $\Omega$.
Thus ${\partial}$ is a~linear map that satisf\/ies the Leibniz identity
\begin{gather*}
{\partial}(ab)=({\partial}(a))b+a{\partial}(b)
\end{gather*}
for all $a, b \in {\mathcal A}$.
Note that this implies that ${\partial}(1_{\mathcal A})=0$, where here $1_{\mathcal A}$ denotes the identity element
of~${\mathcal A}$.
It will be important for us to make the usual requirement that the sub-bimodule of $\Omega$ generated by the range of
${\partial}$ is all of $\Omega$, unless the contrary is explicitly stated.

A Riemannian metric on a~dif\/ferentiable manifold is usually specif\/ied by giving an inner product on the tangent space at
each point of the manifold, but one can equally well use the cotangent space instead of the tangent space.
Then the Riemannian metric gives an inner product on the space of dif\/ferential one-forms that has values in the algebra
of smooth functions on the manifold.
Thus, in generalization of Riemannian metrics, we want to consider~${\mathcal A}$-valued sesquilinear forms on $\Omega$
that are compatible with the right~${\mathcal A}$-module structure on $\Omega$.
(Since we work over ${\mathbb C}$, we actually have an analog of the complexif\/ied cotangent bundle.
One might well want to introduce a~``real'' structure, but we will not discuss that possibility.)

We will require that the left action of~${\mathcal A}$ on $\Omega$ be a~$*$-action with respect to the inner product.
We will not assume any positivity for our inner products until the next section.
When positivity is present, it is usual~\cite{Blk2} to refer to such a~bimodule with~${\mathcal A}$-valued inner product
(no derivation involved) as a~``correspondence''.
In the present more general setting we will use the term``quasi-correspondence''.
We will usually denote such an inner product by $\langle \cdot,\cdot\rangle_{\mathcal A}$.
Thus:

\begin{Definition}
\label{corresq}
Let~${\mathcal A}$ be a~unital $*$-algebra.
By a~(right) \emph{pre-quasi-correspondence} over~${\mathcal A}$ we mean an~${\mathcal A}$-bimodule $\Omega$ that is
equipped with an~${\mathcal A}$-valued sesquilinear form $\langle \cdot,\cdot\rangle_{\mathcal A}$ that satisf\/ies
\begin{gather*}
\langle \omega,\omega' a\rangle_{\mathcal A}=\langle \omega,\omega'\rangle_{\mathcal A} a,
\\
(\langle \omega,\omega'\rangle_{\mathcal A})^*=\langle \omega',\omega\rangle_{\mathcal A}
\end{gather*}
and
\begin{gather*}
\langle a \omega,\omega'\rangle_{\mathcal A}=\langle \omega,a^*\omega'\rangle_{\mathcal A}
\end{gather*}
for all $\omega, \omega' \in \Omega$ and $a \in {\mathcal A}$.
We will refer to the sesquilinear form as a~pre-inner-product.
The null-space, ${\mathcal N}$, of the sesquilinear form is def\/ined
to be ${\mathcal N}=\{\omega: \langle \omega,\omega'\rangle_{\mathcal A}=0~\text{for~all}~\omega' \in \Omega\}$.
If ${\mathcal N}=\{0\}$ then we say that the sesquilinear form is non-degenerate, and we call it an inner product.
We then call $(\Omega, \langle \cdot,\cdot\rangle_{\mathcal A})$ a~(right) \emph{quasi-correspondence} over~${\mathcal A}$.
Left (pre-)quasi-correspondences are def\/ined analogously, with the~${\mathcal A}$-valued inner product linear in the
f\/irst variable.
\end{Definition}

As is commonly done, we will usually work with right (pre-)quasi-correspondences, and will usually omit the word
``right''.

It is easily seen that the null-space, ${\mathcal N}$, of the pre-inner-product for a~pre-quasi-correspondence is
a~sub-bimodule, so that $\Omega/{\mathcal N}$ is an~${\mathcal A}$-bimodule, to which the pre-inner-product drops to
give an inner product, for which $\Omega/{\mathcal N}$ is then a~quasi-correspondence.

\begin{Definition}
Let~${\mathcal A}$ be a~unital $*$-algebra.
By a~\emph{calculus with pre-quasi-correspondence} for~${\mathcal A}$ we mean a~triple $(\Omega, {\partial}, \langle
\cdot,\cdot\rangle_{\mathcal A})$ such that $(\Omega, {\partial})$ is a~f\/irst-order dif\/ferential calculus for
${\mathcal A}$ and $(\Omega, \langle \cdot,\cdot\rangle_{\mathcal A})$ is a~pre-quasi-correspondence over~${\mathcal A}$.
If the pre-inner product is non-degenerate we will call this a~\emph{calculus with quasi-correspondence}.
\end{Definition}

When in the next section we impose positivity we will call this structure a~``Riemannian (pre-)metric'' for~${\mathcal A}$.
Notice that we make no assumption about how ${\partial}$ is related to $\langle \cdot,\cdot\rangle_{\mathcal A}$.
Later we will discuss some relations that one might want to require.

For a~closely related def\/inition of a~Riemannian metric for a~$*$-algebra, coming from quite dif\/ferent motivation,
see~\cite{BgM}, and see~\cite{Mjd} for its application to graphs.
For another interesting, and very new, def\/inition of a~Riemannian metric, in the context to non-commutative tori,
see~\cite{Rsn}.

We remark that Def\/inition~\ref{corresq} is very close to Section~3 of~\cite{Svg2}.
Sauvageot~\cite{Svg1} prefers to view $\Omega$ as an analog of the tangent bundle, and ${\partial}$ as the gradient, and
he does not introduce an~${\mathcal A}$-valued inner product.

If $(\Omega, {\partial}, \langle \cdot,\cdot\rangle_{\mathcal A})$ is a~calculus with pre-quasi-correspondence, and
if ${\mathcal N}$ is the null-space of the pre-inner product as above, then it is easily verif\/ied that
$(\Omega/{\mathcal N}, {\partial}, \langle \cdot,\cdot\rangle_{\mathcal A})$ is a~calculus with quasi-correspondence,
where here we do not change the notation for the derivation and the inner product, but they are def\/ined in the evident
way.

We now give four simple but very pertinent examples.

\begin{Example}
\label{excomm}
Let~$X$ be a~f\/inite set, and let ${\mathcal A}=C(X)$, the algebra of ${\mathbb C}$-valued functions on~$X$ with
pointwise multiplication and with complex-conjugation as involution.
We def\/ine a~f\/irst-order dif\/ferential calculus for~${\mathcal A}$ in a~familiar way.
Let $Z=\{(x,y) \in X\times X: x \neq y\}$, and let $\Omega=C(Z)$.
Then $\Omega$ is an~${\mathcal A}$-bimodule for the operations
\begin{gather*}
(f\omega)(x,y)=f(x)\omega(x, y)
\qquad
\text{and}
\qquad
(\omega f)(x, y)=\omega(x, y)f(y)
\end{gather*}
for $f \in {\mathcal A}$ and $\omega\in \Omega$.
We def\/ine a~derivation ${\partial}$ from~${\mathcal A}$ into $\Omega$ by
\begin{gather*}
({\partial} f)(x,y)=f(x)-f(y).
\end{gather*}
We f\/ind it helpful to view this in the following heuristic way.
For a~given point $y \in X$ the directions in which a~function $f \in {\mathcal A}$ can be ``dif\/ferentiated'' are given
by the points of $X \setminus \{y\}$.
These points form a~basis for the ``tangent space'' at~$y$, and the ``tangent space'' at~$y$ can be considered to be
$C(X \setminus \{y\})$.
The dif\/ferential of~$f$ at~$y$ is then given by the function $x \mapsto ({\partial} f)(x, y)$.
It is easily verif\/ied that the sub-bimodule generated by the range of~${\partial}$ is all of~$\Omega$.

To def\/ine an~${\mathcal A}$-valued pre-inner product on $\Omega$ we choose an ${\mathbb R}$-valued function,~$c$,
on~$Z$.
Eventually~$c$ will provide the conductances for a~resistance network, but at this stage we do not assume that $c(x, y)
=c(y, x)$, nor that~$c$ be non-negative.
We write $c_{xy}$ for $c(x,y)$, and for $\omega, \omega' \in \Omega$ we set
\begin{gather*}
\langle\omega, \omega'\rangle_{\mathcal A}(y)=\sum\limits_{x, x \neq y}\bar \omega(x, y)\omega'(x, y)c_{xy}.
\end{gather*}
For f\/ixed~$y$ this can be viewed as giving a~pre-inner-product on the cotangent space at~$y$.
It is easily verif\/ied that with this pre-inner-product $\Omega$ becomes a~pre-quasi-correspondence over~${\mathcal A}$,
and in this way $(\Omega, {\partial}, \langle\cdot, \cdot\rangle_{\mathcal A})$ is a~calculus with
pre-quasi-correspondence for~${\mathcal A}$.
For some closely related correspondences, but without mention of derivations, see~\cite{KPQ}.
\end{Example}

\begin{Example}
\label{exnc}
Let~${\mathcal A}$ be any non-commutative unital f\/inite-dimensional $*$-algebra.
Let $\tilde \Omega$ be~${\mathcal A}$ viewed as a~bimodule over itself.
Choose any element~$v$ of~${\mathcal A}$ that is not in the center of~${\mathcal A}$, and def\/ine a~derivation of
${\mathcal A}$ into $\tilde \Omega$ by
\begin{gather*}
{\partial}(a)=[v,a]=va-av
\end{gather*}
for all $a \in {\mathcal A}$.
Def\/ine an~${\mathcal A}$-valued pre-inner-product on $\tilde \Omega$ by
\begin{gather*}
\langle a, b\rangle_{\mathcal A}=a^*b
\end{gather*}
for all $a, b \in {\mathcal A}$.
Let $\Omega$ be the sub-bimodule of $\tilde \Omega$ generated by the range of ${\partial}$.
It is easily verif\/ied that with the restriction of the above pre-inner-product, $\Omega$ becomes
a~pre-quasi-correspondence over~${\mathcal A}$, and in this way $(\Omega, {\partial}, \langle\cdot,
\cdot\rangle_{\mathcal A})$ is a~calculus with pre-quasi-correspondence for~${\mathcal A}$.
\end{Example}

We remark that if $\Xi$ is any quasi-correspondence over a~unital $*$-algebra~${\mathcal A}$, then any f\/ixed element
$\xi \in \Xi$ determines an inner derivation ${\partial}^\xi$ from~${\mathcal A}$ into $\Xi$ def\/ined by
\begin{gather*}
{\partial}^\xi(a)=a\xi-\xi a.
\end{gather*}
If we let $\Omega^\xi$ be the sub-bimodule of $\Xi$ generated by the range of ${\partial}^\xi$, and restrict to
$\Omega^\xi$ the~${\mathcal A}$-valued inner-product on $\Xi$, then we obtain a~calculus with quasi-correspondence.

\begin{Example}
\label{exgp}
This next example is somewhat a~combination of the two above.
Let~${\mathcal A}$ be any possibly non-commutative unital f\/inite-dimensional $*$-algebra.
Let~$G$ be a~f\/inite group, and let~$\alpha$ be an action of~$G$ on~${\mathcal A}$ by $*$-automorphisms.
Let $\Xi=C(G, {\mathcal A})$, the vector space of~${\mathcal A}$-valued functions on~$G$.
Def\/ine a~right action of~${\mathcal A}$ on $\Xi$ by
\begin{gather*}
(\xi a)(x)=\xi(x)a.
\end{gather*}
Let~$c$ be a~f\/ixed ${\mathbb R}$-valued function on~$G$, and def\/ine an~${\mathcal A}$-valued pre-inner-product on $\Xi$
by
\begin{gather*}
\langle\xi, \eta \rangle_{\mathcal A}=\sum\limits_\Gamma \xi(x)^* \eta(x) c_x.
\end{gather*}
We def\/ine a~left action of~${\mathcal A}$, denoted by $a \cdot \xi$, by
\begin{gather*}
(a \cdot \xi)(x)=\alpha_x(a)\xi(x).
\end{gather*}
With these def\/initions $\Xi$ is an~${\mathcal A}$-pre-quasi-correspondence.

Let $\omega_0 \in \Xi$ be def\/ined by $\omega_0(x)=1_{\mathcal A}$ for all $x \in G$.
The inner derivation, ${\partial}$, determined by $\omega_0$ is then given by
\begin{gather*}
({\partial} a)(x)=\alpha_x(a)-a
\end{gather*}
for all $x \in G$.
We let $\Omega$ be the sub-bimodule of $\Xi$ generated by the range of ${\partial}$, and we restrict to $\Omega$ the
${\mathcal A}$-valued pre-inner-product on $\Xi$.
Then $(\Omega, {\partial}, \langle\cdot, \cdot \rangle_{\mathcal A})$ is a~calculus with pre-quasi-correspondence for
${\mathcal A}$.

Notice that the structure of $\Omega$ depends strongly on the choice of~$\alpha$.
If~$\alpha$ is the trivial action, then $\Omega=\{0\}$.
Notice also that for any choice of~$\alpha$, if~$\omega$ is in $\Omega$ then $\omega(e)=0$, where~$e$ is the identity
element of~$G$.
Thus we do not need a~value for $c_e$.
Then a~natural choice for~$c$ is the inverse of a~length function on~$G$, or its square, left undef\/ined at~$e$.
This is related to the seminorms prominently used in~\cite{R4, R7}.
\end{Example}

\begin{Example}
\label{exexp}
This example is related to the previous two examples.
Let ${\mathcal B}$ be any possibly non-commutative unital f\/inite-dimensional $*$-algebra, and let~${\mathcal A}$ be
a~unital $*$-subalgebra of ${\mathcal B}$.
We can in the evident way view ${\mathcal B}$ as a~bimodule over~${\mathcal A}$, and of course~${\mathcal A}$ can be
viewed as a~bimodule over itself.
Suppose that~$E$ is a~conditional expectation from ${\mathcal B}$ onto~${\mathcal A}$, that is, an~${\mathcal A}$-bimodule projection from ${\mathcal B}$ onto~${\mathcal A}$ that preserves the involution.
On ${\mathcal B}$ we def\/ine a~(right) pre-inner-product with values in~${\mathcal A}$ by
\begin{gather*}
\langle b, c\rangle_{\mathcal A}=E(b^*c)
\end{gather*}
for all $b, c \in {\mathcal B}$.
It is easily verif\/ied that $({\mathcal B}, \langle\cdot, \cdot \rangle_{\mathcal A})$ is a~pre-quasi-correspondence over
${\mathcal A}$.
Then, as commented just before the previous example, any element of ${\mathcal B}$ will def\/ine an inner derivation from
${\mathcal A}$ into ${\mathcal B}$.
If we let $\Omega$ be the sub-${\mathcal A}$-bimodule of ${\mathcal B}$ generated by the range of this derivation, and
if we restrict to $\Omega$ the above pre-inner-product on ${\mathcal B}$, we obtain a~calculus with
pre-quasi-correspondence.
\end{Example}

Let $(\Omega, {\partial}, \langle\cdot,\cdot\rangle_{\mathcal A})$ be a~calculus with pre-quasi-correspondence over some
unital $*$-algebra~${\mathcal A}$.
Notice that the Leibniz identity implies that the right sub-module of $\Omega$ generated by the range of ${\partial}$ is
in fact a~sub-bimodule, and so by our assumptions it is all of $\Omega$.
That is, every element of $\Omega$ can be expressed as a~f\/inite sum of terms of the form $({\partial} a)b$ for $a, b \in
{\mathcal A}$.
But
\begin{gather*}
\langle({\partial} a)b,({\partial} c) d\rangle_{\mathcal A}=b^*\langle{\partial} b,{\partial}
c\rangle_{\mathcal A} d
\end{gather*}
for all $a, b, c, d \in {\mathcal A}$.
Thus the pre-inner-product is entirely determined by the~${\mathcal A}$-valued form~$\Gamma$ def\/ined on~${\mathcal A}$
by
\begin{gather*}
\Gamma(b, c)=\langle{\partial} b,{\partial} c\rangle_{\mathcal A}.
\end{gather*}
Notice that $\Gamma(1_{\mathcal A}, a)=0$ for all $a \in {\mathcal A}$.
The form~$\Gamma$ is ${\mathbb C}$-sesquilinear, and~${\mathcal A}$-symmetric in the sense that
\begin{gather*}
(\Gamma(b,c))^*=\Gamma(c,b)
\end{gather*}
for $b, c \in {\mathcal A}$, but it has no properties with respect to the right~${\mathcal A}$-module structure.
However~$\Gamma$ does have an important property ref\/lecting the $*$-representation condition of the correspondence.
For $a, b, c \in {\mathcal A}$ we have
\begin{gather*}
0=\langle a{\partial} b,{\partial} c\rangle_{\mathcal A}-\langle{\partial} b,a^*{\partial} c\rangle_{\mathcal A}
\\
\phantom{0}
=\langle{\partial}(ab)-({\partial} a)b,{\partial} c\rangle_{\mathcal A}-\langle{\partial} b,{\partial}
(a^*c)-({\partial} a^*)c\rangle_{\mathcal A}
\\
\phantom{0}
=\Gamma(ab, c)-b^*\Gamma(a, c)-\Gamma(b, a^*c)+ \Gamma(b, a^*) c.
\end{gather*}

That is,
\begin{gather*}
\Gamma(ab, c)-\Gamma(b, a^*c)=b^*\Gamma(a, c)- \Gamma(b, a^*) c.
\end{gather*}

In the setting of Dirichlet forms and quantum semigroups~\cite{CpS} the corresponding form~$\Gamma$ is often called
a~``carr\'e-du-champ'' (or sometimes a~``gradient form'').
Once we require positivity, we will use this terminology.
So at this point we set:

\begin{Definition}
Let~${\mathcal A}$ be a~unital $*$-algebra over ${\mathbb C}$.
By a~(right) \emph{quasi-carr\'e-du-champ} (qCdC) for~${\mathcal A}$ we mean an~${\mathcal A}$-symmetric~${\mathcal A}$-valued ${\mathbb C}$-sesquilinear form~$\Gamma$, linear in the second variable, that satisf\/ies both the condition
$\Gamma(1_{\mathcal A}, a)=0$ for all $a \in {\mathcal A}$, and also the $*$-representation condition
\begin{gather}
\Gamma(ab, c)-\Gamma(b, a^*c)=b^*\Gamma(a, c)-\Gamma(b, a^*) c.
\label{eqstar}
\end{gather}
for all $a, b, c \in {\mathcal A}$.
A~left qCcD is def\/ined similarly, but it is linear in the f\/irst variable, and its $*$-representation condition is given
by
\begin{gather*}
\Gamma(ba, c)-\Gamma(b, ca^*)=b\Gamma(a,c)-\Gamma(b, a^*)c^*.
\end{gather*}
If~$\Gamma$ comes from a~f\/irst-order dif\/ferential calculus as above, then we will say that~$\Gamma$ is the qCdC for the
f\/irst-order dif\/ferential calculus with quasi-correspondence.
\end{Definition}

We remark that if $\Gamma(a,a)$ is self-adjoint for all $a \in {\mathcal A}$, then the usual argument shows that
$\Gamma$ is~${\mathcal A}$-symmetric.

\begin{Example}
\label{exmore}
The qCdC for Example~\ref{excomm} is given by
\begin{gather*}
\Gamma(f,g)(y)=\sum\limits_{x, x \neq y} (\bar f(x)-\bar f(y))(g(x)-g(y))c_{xy},
\end{gather*}
while that for Example~\ref{exnc} is given by
\begin{gather*}
\Gamma_v(a, b)=[v,a]^*[v,b].
\end{gather*}
\end{Example}

\begin{Proposition}
\label{progcal}
Let~${\mathcal A}$ be a~unital $*$-algebra over ${\mathbb C}$, and let~$\Gamma$ be a~qCdC for~${\mathcal A}$.
Then there is a~calculus with quasi-correspondence for~${\mathcal A}$ whose qCdC is~$\Gamma$.
\end{Proposition}
\begin{proof}
Let $\tilde \Omega={\mathcal A} \otimes {\mathcal A}$, with its usual~${\mathcal A}$-bimodule structure, given on
elementary tensors by
\begin{gather*}
a(b\otimes c)d=ab \otimes cd.
\end{gather*}
Def\/ine a~derivation ${\partial}^u$ from~${\mathcal A}$ into $\tilde \Omega$ by
\begin{gather*}
{\partial}^u a=a\otimes 1-1\otimes a.
\end{gather*}
(This is the negative of the usual convention, but seems to be more appropriate when using right quasi-correspondences,
and f\/its well with Examples~\ref{excomm}.) Let $\Omega^u$ be the sub-bimodule of $\tilde \Omega$ generated by the range
of ${\partial}^u$.
It is well-known and easily seen to be the kernel of the bimodule map $m: {\mathcal A} \otimes {\mathcal A}\to{\mathcal A}$
determined by $m(a\otimes b)=ab$.
Thus $\Omega^u$ consists of f\/inite sums $\sum a_j \otimes b_j$ such that $\sum a_jb_j=0$.
(We remark that when ${\mathcal A}=C(X)$ for a~f\/inite set~$X$ then $\Omega^u$ is exactly the $C(Z)$ of
Example~\ref{excomm}.) Then $(\Omega^u, {\partial}^u)$ is universal~\cite{GVF} in the sense that if $(\Omega,
{\partial})$ is any other f\/irst-order dif\/ferential calculus for~${\mathcal A}$, then the mapping~$\Phi$ that sends
\begin{gather*}
\sum a_j \otimes b_j=\sum(a_j \otimes 1-1 \otimes a_j)b_j=\sum ({\partial}^u a_j)b_j
\end{gather*}
to $\sum ({\partial} a_j)b_j$ is a~surjective bimodule homomorphism with the property that $\Phi({\partial}^u a)={\partial}(\Phi a)$.

Suppose now that~$\Gamma$ is a~qCdC on~${\mathcal A}$.
Let $B_\Gamma$ be the~${\mathcal A}$-valued 4-linear form def\/ined on~${\mathcal A}$ by
\begin{gather*}
B_\Gamma(a,b,c,d)=a\Gamma(b^*, c)d.
\end{gather*}
It extends to an~${\mathcal A}$-valued linear form on ${\mathcal A}^{\otimes4}$, which we can view as a~bilinear form on
${\mathcal A} \otimes {\mathcal A}$.
From $B_\Gamma$ we can then def\/ine an~${\mathcal A}$-valued pre-inner-product on ${\mathcal A} \otimes {\mathcal A}$,
denoted by $\langle \cdot, \cdot \rangle_{\mathcal A}^\Gamma$.
It is given on elementary tensors by
\begin{gather*}
\langle a\otimes b, c \otimes d\rangle_{\mathcal A}^\Gamma=B_\Gamma(b^*, a^*, c, d)=b^*\Gamma(a, c)d.
\end{gather*}
We can then restrict this pre-inner product to $\Omega^u$.
With this def\/inition it is easily seen that we have the following properties:
\begin{gather*}
\langle\omega a, \omega' b\rangle_{\mathcal A}^\Gamma=a^*\langle\omega, \omega' \rangle_{\mathcal A}^\Gamma b,
\\
(\langle\omega, \omega' \rangle_{\mathcal A}^\Gamma)^*=\langle\omega', \omega \rangle_{\mathcal A}^\Gamma,
\\
\langle{\partial}^u a, {\partial}^u b \rangle_{\mathcal A}^\Gamma=\Gamma(a, b),
\\
\langle a\omega, \omega' \rangle_{\mathcal A}^\Gamma=\langle\omega, a^*\omega' \rangle_{\mathcal A}^\Gamma
\end{gather*}
for all $\omega, \omega' \in \Omega$ and $a, b \in {\mathcal A}$.
To obtain the last relation, notice that for $a, b, c \in {\mathcal A}$ we have
\begin{gather*}
\langle a{\partial}^u b, {\partial}^u c\rangle_{\mathcal A}^\Gamma=\langle{\partial}^u (ab), {\partial}^u c\rangle_{\mathcal
A}^\Gamma-\langle({\partial}^u a)b, {\partial}^u c\rangle_{\mathcal A}^\Gamma=\Gamma(ab, c)-b^*\Gamma(a, c),
\end{gather*}
to which we can apply the $*$-representation condition (equation~\eqref{eqstar}) on~$\Gamma$.
Thus we see that $(\Omega^u, {\partial}^u, \langle\cdot, \cdot \rangle_{\mathcal A}^\Gamma)$ is a~calculus with
pre-quasi-correspondence.
Let ${\mathcal N}$ be the null-space for the pre-inner-product, and let $\Omega^\Gamma=\Omega^u/{\mathcal N}$, to
which the pre-inner-product drops as an inner product.
Let ${\partial}^\Gamma$ be the composition of ${\partial}^u$ with the quotient map to $\Omega^\Gamma$.
Then $(\Omega^\Gamma, {\partial}^\Gamma, \langle\cdot, \cdot\rangle_{\mathcal A}^\Gamma)$ is a~calculus with
quasi-correspondence whose qCdC is~$\Gamma$, as desired.
\end{proof}

There is an evident notion of isomorphism between any two calculi-with-correspondence over~${\mathcal A}$.
It is easy to verify that:

\begin{Theorem}
The above construction gives a~natural bijection between qCdC's over~${\mathcal A}$ and isomorphism classes of
calculi-with-quasi-correspondence for~${\mathcal A}$.
\end{Theorem}

\begin{Example}
\label{exconst}
We now describe an important construction of qCdC's which we will use later.
We phrase this construction in terms of the beginnings of Hochschild cohomology (e.g.\
\cite[p.~187]{Cn2}), but it is not clear to me whether it is useful to do this.
Let~$N$ be any operator on~${\mathcal A}$ (${\mathcal A}$-valued, ${\mathbb C}$-linear) with the property that
$N(1_{\mathcal A})=0$.
Let $\hat N$ be the ${\mathbb C}$-trilinear~${\mathcal A}$-valued form on~${\mathcal A}$ def\/ined by
\begin{gather*}
\hat N(a,b,c)=aN(b)c,
\end{gather*}
extended to give an~${\mathcal A}$-bimodule map from ${\mathcal A}\otimes {\mathcal A} \otimes {\mathcal A}$ to
${\mathcal A}$.
We can view $\hat N$ as a~Hochschild 2-cochain for~${\mathcal A}$ with coef\/f\/icients in the bimodule~${\mathcal A}$
(see equation~(8.47) of~\cite{GVF}).
Then its coboundary, $\delta \hat N$, is the bimodule map from ${\mathcal A}^{\otimes
4}$ to~${\mathcal A}$ def\/ined on elementary tensors by
\begin{gather*}
(\delta \hat N)(a \otimes b \otimes c \otimes d)=a(N(b)c-N(bc)+bN(c))d.
\end{gather*}
We can turn this into an~${\mathcal A}$-valued sesquilinear form on $\tilde \Omega={\mathcal A} \otimes {\mathcal A}$,
denoted by $\langle\cdot, \cdot\rangle_{\mathcal A}^N$, def\/ined on elementary tensors by
\begin{gather*}
\langle a \otimes b, c \otimes d\rangle_{\mathcal A}^N=b^*(N(a^*) c-N(a^*c)+a^*N(c))d.
\end{gather*}
We can then restrict this sesquilinear form to $\Omega^u$.
We can hope that this gives a~pre-quasi-correspondence for $(\Omega^u, {\partial}^u)$.
It is clear that then its qCdC would be given by
\begin{gather*}
\Gamma_N(a, c)=N(a^*) c-N(a^*c)+a^*N(c).
\end{gather*}
Notice that $\Gamma_N$ measures the extent to which~$N$ fails to be a~derivation on~${\mathcal A}$.
In particular, two~$N$'s that dif\/fer by a~derivation will give the same $\Gamma_N$.
(We will see later that it can be convenient to include a~factor of $1/2$ in the def\/inition of $\Gamma_N$.)

We now seek to determine when $\Gamma_N$ is indeed a~qCdC.
It is clear that $\Gamma_N(1_{\mathcal A}, a)=0=\Gamma_N(a, 1_{\mathcal A})$, because $N(1_{\mathcal A})=0$.
We check next that $\Gamma_N$ satisf\/ies the $*$-representation condition.
For $a, b, c \in {\mathcal A}$ we have
\begin{gather*}
\Gamma_N(ab, c)-\Gamma_N(b, a^*c)
=b^*(a^*N(c)-N(a^*c))-(N(b^*)a^*-N(b^*a^*))c
\\
\qquad{}
=b^*(N(a^*)c-N(a^*c)+a^*N(c)) -(N(b^*)a^*-N(b^*a^*) +b^*N(a^*))c
\\
\qquad{}
=b^*\Gamma_N(a,c)-\Gamma_N(b, a^*)c,
\end{gather*}
as desired.

However, $\Gamma_N$ will not in general be symmetric.
Notice that for $a, b \in {\mathcal A}$ we have
\begin{gather*}
(\Gamma_N(b,a))^*=a^*(N(b^*))^*-(N(b^*a))^*+(N(a))^*b
\\
\hphantom{(\Gamma_N(b,a))^*}{}
=a^*N^\sharp(b)-N^\sharp(a^*b)+N^\sharp(a^*)b
=\Gamma_{N^\sharp}(a, b),
\end{gather*}
where we def\/ine $N^\sharp$ by $N^\sharp(c)=(N(c^*))^*$ for $c \in {\mathcal A}$.
Thus $\Gamma_N$ will be symmetric exactly if $\Gamma_N=\Gamma_{N^\sharp}$, and so exactly if $N-N^\sharp$ is
a~derivation of~${\mathcal A}$.
We have thus obtained:

\begin{Theorem}\label{thmsym}
Let~$N$ be a~linear operator on~${\mathcal A}$ that satisfies $N(1_{\mathcal A})=0$, and define $\Gamma_N$ on ${\mathcal
A} \times {\mathcal A}$ by
\begin{gather*}
\Gamma_N(a, b)=N(a^*)b-N(a^*b)+a^*N(b).
\end{gather*}
Then $\Gamma_N$ is a~qCdC exactly if $N-N^\sharp$ is a~derivation of~${\mathcal A}$, where we define $N^\sharp$ by
$N^\sharp(c)=(N(c^*))^*$ for $c \in {\mathcal A}$.
\end{Theorem}

\end{Example}

\begin{Example}
\label{exqc}
In anticipation of what will come in Example~\ref{excont}, let us consider the case, associated to Examples~\ref{exnc}
and~\ref{exmore}, in which we f\/ix $v \in {\mathcal A}$ and def\/ine an operator $N_v$ by
\begin{gather*}
N_v(a)=[v^*, [v, a]]
\end{gather*}
for all $a \in {\mathcal A}$.
We will show that $\Gamma_{N_v}$ is a~qCdC.
Notice that
\begin{gather*}
N_v^\sharp(a)=[v^*,[v,a^*]]^*=N_{v^*}(a),
\end{gather*}
that is, $N_v^\sharp=N_{v^*}$.
Then according to Theorem~\ref{thmsym}, in order for $\Gamma_{N_v}$ to be symmetric we must show that $N_v-N_v^\sharp$
is a~derivation of~${\mathcal A}$.
But by the Jacobi identity
\begin{gather*}
N_v(a)-N_v^\sharp(a)=[v^*, [v,a]]-[v, [v^*, a]]=[[v^*,v],a],
\end{gather*}
so that $N_v-N_v^\sharp$ is indeed a~derivation, and thus $\Gamma_{N_v}$ is a~qCdC.
Notice that $N_v=N_v^\sharp$ exactly when $[v^*, v]$ is in the center of~${\mathcal A}$.

Let us now calculate $\Gamma_{N_v}$.
For $a, b \in {\mathcal A}$ we have
\begin{gather*}
\Gamma_{N_v}(a,b)=[v^*, [v,a^*]]b-[v^*, [v,a^*b]]+a^*[v^*, [v,b]]
\\
\phantom{\Gamma_{N_v}(a,b)}
=[v^*, [v,a^*]]b-[v^*, [v,a^*]b+a^*[v,b]]+a^*[v^*, [v,b]]
\\
\phantom{\Gamma_{N_v}(a,b)}
=[v^*, a]^*[v^*, b]+[v,a]^*[v, b]
\\
\phantom{\Gamma_{N_v}(a,b)}
=\Gamma_{v^*}(a, b)+\Gamma_v(a, b),
\end{gather*}
for which we recall that $\Gamma_v$ was def\/ined by $\Gamma_v(a, b)=[v,a]^*[v,b]$ in Example~\ref{exmore}.
We see that $\Gamma_{N_v}=2\Gamma_v$ exactly if $\Gamma_{v^*}=\Gamma_v$ (which is one example of why a~factor of 1/2
would be convenient in the def\/inition of $\Gamma_N$).
Now if $\Gamma_{v^*}=\Gamma_v$, then for all $a, b \in {\mathcal A}$ we have $[v^*,a]^*[v^*,b]=[v,a]^*[v,b]$.
If we set in this $a=v=b$ we obtain $[v^*,v]^*[v^*,v]=0$.
If~${\mathcal A}$ has the property that $a^*a=0$ implies that $a=0$, as happens for $C^*$-algebras, then we see that
$[v^*, v]=0$, that is,~$v$ is ``normal''.
Since non-normal elements are common in $C^*$-algebras, the property $\Gamma_{N_v}=2\Gamma_v$ can easily fail.
\end{Example}

\section{Non-commutative Riemannian metrics}
\label{ncrm}

We now assume that~${\mathcal A}$ is a~(f\/inite-dimensional) $C^*$-algebra.
(Thus~${\mathcal A}$ can be realized as a~unital $*$-subalgebra of the algebra of all linear operators on a~f\/inite
dimensional Hilbert space, and~${\mathcal A}$ is equipped with the corresponding operator norm.) It is thus meaningful
to consider positive elements of~${\mathcal A}$, that is, self-adjoint elements of~${\mathcal A}$ whose spectrum
(i.e.~set of eigenvalues) is contained in the non-negative real numbers.

Accordingly, we will now require that the~${\mathcal A}$-valued pre-inner-products that we consider on~$\Omega$ are
non-negative, that is, that $\langle\omega,\omega\rangle_{\mathcal A} \geq 0$ for all~$\omega$.
Thus as right~${\mathcal A}$-modules our $\Omega$'s will be right pre-Hilbert~${\mathcal A}$-modules, as def\/ined for
example in Section~II.7.1 of~\cite{Blk2}
(see also Def\/inition~2.8 in~\cite{R26}).
We remark that because positive elements are self-adjoint, this implies the
symmetry of the~${\mathcal A}$-valued pre-inner-products.

\begin{Definition}
By a~\emph{pre-correspondence} we will mean a~pre-quasi-correspondence whose pre-inner-product is non-negative.
A~\emph{correspondence} is then a~pre-correspondence whose pre-inner-product is def\/inite
(see Section~II.7.4.4 of~\cite{Blk2}).
\end{Definition}

\begin{Definition}
Let~${\mathcal A}$ be a~(f\/inite-dimensional) $C^*$-algebra.
By a~(right) \emph{Riemannian pre-metric} for~${\mathcal A}$ we mean a~calculus with pre-correspondence, $(\Omega,
{\partial}, \langle\cdot, \cdot\rangle_{\mathcal A})$, over~${\mathcal A}$.
If the pre-inner-product is def\/inite, then we will call $(\Omega, {\partial}, \langle\cdot, \cdot\rangle_{\mathcal A})$
a~(right) \emph{Riemannian metric} for~${\mathcal A}$.
\end{Definition}

We remark that it would be natural to require also that if ${\partial} a~=0$ then $a \in {\mathbb C} 1_{\mathcal A}$,
but it will be more convenient for us not to require this property, and to view the failure of this property to mean
that $({\mathcal A}, {\partial})$ is not ``metrically connected''.

For positive~${\mathcal A}$-valued pre-inner-products there is a~corresponding Cauchy-Schwartz inequality.
See Proposition~2.9 of~\cite{R26},
or Lemma~2.5 of~\cite{RbW},
or Proposition~II.7.1.4 of~\cite{Blk2}.
It states that for any $\omega, \omega' \in \Omega$ we have
\begin{gather}
\langle\omega, \omega'\rangle_{\mathcal A}^*\langle\omega,\omega'\rangle_{\mathcal A}\leq \|\langle\omega,
\omega\rangle_{\mathcal A}\|\langle\omega', \omega'\rangle_{\mathcal A}
\label{eqCS}
\end{gather}
with respect to the partial order on positive elements of~${\mathcal A}$.
From this inequality one sees by the usual argument that the null-space, ${\mathcal N}$, of the pre-inner-product is
a~right~${\mathcal A}$-submodule of $\Omega$, and in fact is an~${\mathcal A}$-sub-bimodule because the left action of
${\mathcal A}$ is a~$*$-representation.
Then the pre-inner-product drops to an~${\mathcal A}$-valued inner product on $\Omega/{\mathcal N}$.
This inner product determines a~norm, $\|\langle\omega, \omega\rangle_{\mathcal A}\|^{1/2}$, on $\Omega/{\mathcal N}$,
and since in our f\/inite-dimensional situation $\Omega/{\mathcal N}$ is complete for this norm, $\Omega/{\mathcal N}$ is
a~right Hilbert $C^*$-module over~${\mathcal A}$, as def\/ined in Section~II.7.1of~\cite{Blk2}.
The left action then makes $\Omega/{\mathcal N}$ into a~correspondence over~${\mathcal A}$ exactly as def\/ined for
$C^*$-algebras
(see Section~II.7.4.4 of~\cite{Blk2}). We will denote the composition of ${\partial}$ with the quotient map from~$\Omega$ to~$\Omega/{\mathcal N}$ again by~${\partial}$.
Then $(\Omega/{\mathcal N}, {\partial}, \langle\cdot, \cdot\rangle_{\mathcal A})$ will be a~Riemannian metric for~${\mathcal A}$.
In this way we can always pass from a~Riemannian pre-metric to a~Riemanian metric.

\begin{Example}
\label{exmulti}
Examples~\ref{excomm},~\ref{exnc},~\ref{exgp}, and~\ref{exexp}, after evident slight modif\/ications, give examples of
pre-Riemannian metrics.
For instance, in Example~\ref{excomm} we must assume that the function~$c$ takes non-negative values, in
Examples~\ref{exnc} and~\ref{exgp} we must assume that~${\mathcal A}$ is a~unital $C^*$-algebra, and in
Example~\ref{exexp} we must assume that~${\mathcal A}$ and ${\mathcal B}$ are unital $C^*$-algebras and that the
conditional expectation~$E$ is non-negative, so that it is a~conditional expectation in the sense used for
$C^*$-algebras~\cite{Tks}.
It would be interesting to see if Example~\ref{excomm} can be generalized to the setting of~\cite{SnP} and related
papers, and whether the results in this paper can then be applied to the setting of those papers.
\end{Example}

We now give a~further example.

\begin{Example}
\label{exspec}
Let $({\mathcal A}, {\mathcal H}, D)$ be a~f\/inite-dimensional spectral triple~\cite{Cn2, GVF}, that is,~${\mathcal A}$
is a~f\/inite-dimensional $C^*$-algebra, ${\mathcal H}$ is a~f\/inite-dimensional Hilbert space on which~${\mathcal A}$ is
represented, and~$D$ is a~self-adjoint operator on ${\mathcal H}$.
For ease of discussion we assume that the representation of~${\mathcal A}$ is faithful, and so we just take~${\mathcal A}$ to be a~$*$-subalgebra of the $C^*$-algebra ${\mathcal B}({\mathcal H})$ of all operators on ${\mathcal H}$, with
$1_{\mathcal A}$ coinciding with the identity operator on ${\mathcal H}$.
Let~$\tau$ be the unique tracial state on ${\mathcal B}({\mathcal H})$.
Then the dual transformation of the inclusion of ${\mathcal L}^1({\mathcal A}, \tau)$ into ${\mathcal L}^1({\mathcal
B}({\mathcal H}), \tau)$ is a~conditional expectation,~$E$, from ${\mathcal B}({\mathcal H})$ onto~${\mathcal A}$
(see Proposition~V.2.36 of~\cite{Tks}).
Then as in Example~\ref{exexp} we obtain a~Riemannian metric for~${\mathcal A}$
whose bimodule is the~${\mathcal A}$-sub-bimodule of~${\mathcal B}({\mathcal H})$ generated by the range of the
derivation~$a \mapsto [D, a]$.
Thus in our f\/inite-dimensional setting every spectral triple has a~canonically associated Riemannian metric.
Note that dif\/ferent~$D$'s on ${\mathcal H}$ can def\/ine the same derivation, and thus the same Riemannian metric.
More generally, dif\/ferent spectral triples for a~given~${\mathcal A}$ can determine isomorphic Riemannian metrics.
(For a~related inf\/inite-dimensional version see Theorem~2.9 of~\cite{FGR}.
I~thank D.~Goswami for bring this theorem to my attention.)
\end{Example}

Suppose now that $(\Omega, {\partial}, \langle\cdot, \cdot\rangle_{\mathcal A})$ is a~Riemannian pre-metric for
a~f\/inite-dimensional $C^*$-algebra~${\mathcal A}$.
Then in particular, it is a~calculus with pre-correspondence.
Let~$\Gamma$ be its qCdC as discussed in the previous section.
Now for every element~$\omega$ of $\Omega$, expressed as a~f\/inite sum $\omega=\sum\limits^n ({\partial} a_j)b_j$, we have
\begin{gather*}
0 \leq \langle\omega, \omega\rangle_{\mathcal A}=\sum\limits^n b_j^*\Gamma(a_j, a_k) b_k.
\end{gather*}
This implies exactly that the matrix $\{\Gamma(a_j, a_k)\}$ is a~positive element of the $C^*$-algebra $M_n({\mathcal
A})$.
The fact that this holds for all~$n$ and all choices of the $a_j$'s is exactly what is meant by saying that~$\Gamma$ is
``completely positive''.

\begin{Definition}
\label{cdc}
Let~${\mathcal A}$ be a~f\/inite-dimensional $C^*$-algebra.
By a~(right) \emph{carr\'e-du-champ} (CdC) for~${\mathcal A}$ we mean a~qCdC for~${\mathcal A}$ that is completely
positive.
(No def\/initeness is required.)
\end{Definition}

Notice that since positive elements of~${\mathcal A}$ are self-adjoint, a~CdC will automatically be symmetric, as
mentioned before Example~\ref{exmore}.
The sum of two CdC's is again a~CdC, and a~positive scalar multiple of a~CdC is again a~CdC, so the CdC's form a~cone.

From this def\/inition and Theorem~\ref{thmsym} we immediately obtain:

\begin{Proposition}
\label{proncdc}
Let~$N$ be an operator on~${\mathcal A}$ such that $N(1_{\mathcal A})=0$.
As in Theorem~{\rm \ref{thmsym}} define $\Gamma_N$ by
\begin{gather*}
\Gamma_N(a, b)=N(a^*)b-N(a^*b)+a^*N(b).
\end{gather*}
Then $\Gamma_N$ is a~CdC if and only if it is completely positive and $N-N^\sharp$ is a~derivation of~${\mathcal A}$.
\end{Proposition}

We will characterize such~$N$'s in Theorem~\ref{thmccn}.

We remark that in our f\/inite-dimensional situation all derivations of~${\mathcal A}$ are inner.
A~relatively simple proof of this can be extracted from Exercise 8.7.53 of~\cite{KR2}.

\begin{Example}[following Lindblad~\cite{Lbd}] \label{exsg}Let $\{\Phi_t\}$ be a~quantum dynamical semigroup on~${\mathcal A}$, that is, for every
$t \in {\mathbb R}_{\geq 0}$ the operator $\Phi_t$ on~${\mathcal A}$ is completely positive and of norm no greater than~1, and $t \mapsto \Phi_t$ is a~continuous semigroup homomorphism from ${\mathbb R}_{\geq 0}$ into the algebra of bounded
operators on~${\mathcal A}$ (with $\Phi_0$ the identity operator on~${\mathcal A}$).
Assume further that this semigroup is ``conservative'' in the sense that $\Phi_t(1_{\mathcal A})=1_{\mathcal A}$ for all~$t$.
Especially in Section~6 of~\cite{JMP} (and references there) the attitude is taken that such semigroups are a~good
substitute for metrics in the non-commutative setting.
It is not dif\/f\/icult to show that, because we assume that~${\mathcal A}$ is f\/inite-dimensional, the function $t \mapsto
\Phi_t$ is dif\/ferentiable.
Let $-\Delta$ denote its derivative at~0, so that for all $a \in {\mathcal A}$ we have
\begin{gather*}
-\Delta(a)=\lim_{t \rightarrow 0} (\Phi_t(a)-a)/t.
\end{gather*}
Then~$\Delta$ is the generator for $\{\Phi_t\}$ in the sense that $\Phi_t=e^{-t\Delta}$ for all~$t$.
Because $\Phi_t(1_{\mathcal A})=1_{\mathcal A}$ for all~$t$ we have $\Delta(1_{\mathcal A})=0$.
Because $\Phi_t$ is completely positive, we have the basic inequality (Proposition~II.6.9.14 of~\cite{Blk2})
\begin{gather*}
\Phi_t(a^*a)-\Phi_t(a^*)\Phi_t(a) \geq 0
\end{gather*}
for all $a \in {\mathcal A}$.
When we dif\/ferentiate this inequality at $t=0$, noting that the left-hand side has value 0 for $t=0$, we obtain
\begin{gather*}
-\Delta(a^*a)+\Delta(a^*)a+a^*\Delta(a) \geq 0.
\end{gather*}
As in Example~\ref{exconst}, set
\begin{gather*}
\Gamma_\Delta(a, b)=\Delta(a^*)b-\Delta(a^*b)+a^*\Delta(b)
\end{gather*}
for all $a, b \in {\mathcal A}$.
We see that $\Gamma_\Delta$ is a~positive~${\mathcal A}$-valued form.
Because $\Phi_t$ is positive for all~$t$ we have that $\Delta(a^*)=(\Delta(a))^*$ for all $a \in {\mathcal A}$, so
that $\Delta^\sharp=\Delta$.
Because each $\Phi_t$ is completely positive, all of the above observations apply equally well to the semigroup $I_n
\otimes \Phi_t$ acting on $M_n \otimes {\mathcal A}$, whose generator is $I_n \otimes \Delta$.
It follows that $\Gamma_\Delta$ is completely positive.
It then follows from Proposition~\ref{proncdc} that $\Gamma_\Delta$ is a~CdC.
\end{Example}

We remark that Lindblad shows in~\cite{Lbd} that, conversely, under the conditions obtained just above on~$\Delta$ it
will always generate a~quantum dynamical semigroup.
For a~very interesting recent account of some uses of quantum dynamical semigroups in quantum physics see~\cite{WsM},
especially the ``master equations'' in the chapter ``Open systems''.
I thank Eleanor Rief\/fel for bringing this reference to my attention.

\begin{Theorem}
\label{thmcst}
For every CdC~$\Gamma$ for~${\mathcal A}$ there exists a~Riemannian metric for~${\mathcal A}$ whose CdC is~$\Gamma$.
\end{Theorem}

\begin{proof}
Let~$\Gamma$ be a~CdC for~${\mathcal A}$.
As in the proof of Proposition~\ref{progcal} we def\/ine on $\tilde \Omega={\mathcal A} \otimes {\mathcal A}$ an
${\mathcal A}$-valued sesquilinear form determined on elementary tensors by
\begin{gather*}
\langle a\otimes b, c \otimes d\rangle_{\mathcal A}^\Gamma=b^*\Gamma(a, c)d,
\end{gather*}
and then we restrict it to $\Omega^u$.
The positivity of the resulting form follows immediately from the complete positivity of~$\Gamma$.
(This is closely related to the Stinespring construction~\cite{Blk2}.) The other properties of a~correspondence then
follow from the fact that~$\Gamma$ is a~qCdC.
As before, we set
\begin{gather*}
{\partial}^u a=a\otimes 1- 1 \otimes a,
\end{gather*}
and then, as for a~qCdC, we have
\begin{gather*}
\Gamma(a, b)=\langle{\partial}^u a, {\partial}^u b\rangle_{\mathcal A}^\Gamma.
\end{gather*}
We thus f\/ind that $(\Omega^u, {\partial}^u, \langle\cdot, \cdot\rangle_{\mathcal A}^\Gamma)$ is a~Riemannian pre-metric
for~${\mathcal A}$.
From this one can then pass to a~Riemannian metric in the way described before Example~\ref{exmulti}.
\end{proof}

In particular, from Example~\ref{exsg} we see that every quantum semigroup on a~f\/inite-dimensional $C^*$-algebra has
a~canonically associated Riemannian metric.

For $C^*$-algebras, conditional expectations are def\/ined~\cite{Blk2} as in Example~\ref{exexp} but with the additional
condition that they carry positive elements to positive elements.
The following proposition is easily proved by direct calculation:

\begin{Proposition}
Let ${\mathcal B}$ be a~unital $C^*$-subalgebra of~${\mathcal A}$, and let~$E$ be a~conditional expectation from
${\mathcal A}$ onto ${\mathcal B}$.
For every CdC $\Gamma$ for~${\mathcal A}$ we obtain a~CdC $\hat \Gamma$ for ${\mathcal B}$ by setting $\hat
\Gamma(b, c)=E(\Gamma(b, c))$ for all $b, c \in {\mathcal B}$.
\end{Proposition}

We remark that in our f\/inite-dimensional case (in which our $C^*$-algebras are abstract von Neumann algebras) there will
always exist a~conditional expectation from~${\mathcal A}$ onto ${\mathcal B}$.

An interesting relationship between quantum dynamical semigroups and conditionally completely negative operators was
f\/irst presented by Evans~\cite{Evn}.
For a~more recent account see Section~1 of~\cite{Ptr}.
We recall that an operator~$N$ on a~$C^*$-algebra~${\mathcal A}$ is said to be ``conditionally completely negative'' if
whenever we have $\sum\limits^n a_j b_j=0$ for elements of~${\mathcal A}$ then $\sum\limits^n b_j^* N(a_j^* a_k) b_k\leq0$.
Within our setting, the relationship to $\Gamma_N$ as def\/ined in Theorem~\ref{thmsym} is given by:

\begin{Proposition}
\label{proccn}
Let~$N$ be an operator on a~unital $C^*$-algebra~${\mathcal A}$ with the property that $N(1_{\mathcal A})=0$.
Then $\Gamma_N$ is completely positive if and only if~$N$ is conditionally completely negative.
\end{Proposition}

\begin{proof}
Suppose that $\sum\limits^n a_j b_j=0$ for elements of~${\mathcal A}$.
Then
\begin{gather*}
\sum\limits^n b_j^*\Gamma_N(a_j,a_k)b_k
=\sum\limits_j^n b_j^*N(a_j^*)\sum\limits_k^na_k b_k-\sum\limits^n b_j^* N(a_j^*a_k)b_k+\sum\limits_j^n
b_j^*a_j^*\sum\limits_k^n N(a_k)b_k
\\
\hphantom{\sum\limits^n b_j^*\Gamma_N(a_j,a_k)b_k}{}
=- \sum\limits^n b_j^* N(a_j^*a_k)b_k.
\end{gather*}
From this calculation it is clear that if $\Gamma_N$ is completely positive then~$N$ is conditionally completely
negative.
Conversely, suppose that~$N$ is conditionally completely negative, and suppose that we are given elements $a_1,\dots,a_n$,
$b_1,\dots,b_n$ of~${\mathcal A}$.
Set $b_{n+1}=- \sum\limits^n a_j b_j$ and $a_{n+1}=1_{\mathcal A}$.
Then $\sum\limits^{n+1} a_j b_j=0$, and so, on using the above calculation towards the end, we have
\begin{gather*}
\sum\limits^n b_j\Gamma_N(a_j, a_k)b_k
=\sum\limits^{n+1} b_j\Gamma_N(a_j, a_k)b_k-\sum\limits^{n+1}_j b_j\Gamma_N(a_j, 1_{\mathcal A})b_{n+1}
\\
\hphantom{\sum\limits^n b_j\Gamma_N(a_j, a_k)b_k=}
{}-\sum\limits^{n+1}_k b_{n+1}\Gamma_N(1_{\mathcal A}, a_k)b_k+b_{n+1}\Gamma_N(1_{\mathcal A}, 1_{\mathcal A})b_{n+1}
\\
\hphantom{\sum\limits^n b_j\Gamma_N(a_j, a_k)b_k}
=\sum\limits^{n+1} b_j^*\Gamma_N(a_j, a_k)b_k=- \sum\limits^{n+1} b_j^* N(a_j^*a_k)b_k\geq 0.
\end{gather*}
Thus $\Gamma_N$ is completely positive.
\end{proof}

When we combine this with Proposition~\ref{proncdc} we obtain:
\begin{Theorem}
\label{thmccn}
Let~$N$ be an operator on~${\mathcal A}$ such that $N(1_{\mathcal A})=0$.
Then $\Gamma_N$ is a~CdC if and only if~$N$ is conditionally completely negative and $N-N^\sharp$ is a~derivation of
${\mathcal A}$.
\end{Theorem}

\begin{Example}
\label{exgen}
For any $v, h \in {\mathcal A}$ it is easily verif\/ied that the maps $a \mapsto -v^*av$ and $a \mapsto ha+ah^*$ are
conditionally completely negative.
Thus if we def\/ine $N_v$ by
\begin{gather*}
N_v(a)=-v^*av+(1/2)(v^*va+av^*v),
\end{gather*}
we see that $N_v$ is conditionally completely negative.
But it is clear that $N_v(1_{\mathcal A})=0$ and $N^\sharp=N$.
Then from Theorem~\ref{thmccn} it follows that $\Gamma_N$ is a~CdC.
It is easily calculated that
\begin{gather*}
\Gamma_{N_v}(a, b)=[v,a]^*[v,b].
\end{gather*}

More generally, for any given $v_1, \dots, v_m \in {\mathcal A}$, if we set $N=\sum N_{v_j}$, then $\Gamma_N$ will be
a~CdC.
This should be compared with Theorem~2 of~\cite{Lbd}, which implies that for the case in which ${\mathcal A}=
M_n({\mathbb C})$ the above~$N$'s are the most general form of generators of conservative dynamical semigroups on
${\mathcal A}$, up to a~Hamiltonian.
For examples from physics see equations~(3.29) and~(3.81) of~\cite{WsM} and the text around them.
See also equation~(3.154) of~\cite{Prs}.
\end{Example}

\section{Energy forms, and Markov and Leibniz seminorms}
\label{energy}

In this section we assume that $(\Omega, {\partial}, \langle\cdot, \cdot\rangle_{\mathcal A})$ is a~Riemannian
pre-metric for~${\mathcal A}$.
In the present setting this structure does not seem to lead to a~canonical integration procedure, in contrast to the
case of ordinary Riemannian manifolds.
We will just need to choose an integration procedure, perhaps satisfying some compatibility requirement.
We begin by considering a~somewhat more general procedure.

Let ${\mathcal D}$ be a~unital central subalgebra of~${\mathcal A}$, and let~$E$ be a~conditional
expectation~\cite{Blk2} from~${\mathcal A}$ onto ${\mathcal D}$ that satisf\/ies the extra ``tracial'' condition that
$E(ab)=E(ba)$ for all $a, b \in {\mathcal A}$.
This means that~$E$ is a~``${\mathcal D}$-valued trace'' as def\/ined in Def\/inition~V1.24 of~\cite{Tks}.
The two classes of examples that are most immediately evident are, f\/irst, that in which ${\mathcal D}$ is
one-dimensional,~$\tau$ is a~tracial state on~${\mathcal A}$, and $E(a)=\tau(a)1_{\mathcal A}$; and, second, the class
in which~${\mathcal A}$ is commutative, ${\mathcal D}={\mathcal A}$, and~$E$ is the identity map on~${\mathcal A}$.
We will not discuss other classes in this paper.

We then def\/ine a~${\mathcal D}$-valued pre-inner-product on $\Omega$ by
\begin{gather*}
\langle\omega, \omega'\rangle_{\mathcal D}=E(\langle\omega,\omega'\rangle_{\mathcal A}).
\end{gather*}
It is easy to check that this makes $\Omega$ into a~pre-Hilbert ${\mathcal D}$-module for the right action of~${\mathcal
D}$ on~$\Omega$ coming from the right action of~${\mathcal A}$.
Furthermore, the left action of~${\mathcal A}$ on $\Omega$ will again be a~$*$-representation with respect to
$\langle\cdot, \cdot\rangle_{\mathcal D}$.
The additional feature that we gain is that also the right anti-representation of~${\mathcal A}$ on $\Omega$ will be
a~$*$-anti-representation for $\langle\cdot, \cdot\rangle_{\mathcal D}$ because
\begin{gather*}
\langle\omega, \omega'a\rangle_{\mathcal D}=E(\langle\omega,\omega'\rangle_{\mathcal A} a)
=E(a\langle\omega,\omega'\rangle_{\mathcal A})=E(\langle\omega a^*,\omega'\rangle_{\mathcal A})
=\langle\omega a^*,\omega'\rangle_{\mathcal D},
\end{gather*}
where we have used the tracial condition on~$E$.

We are now in position to apply a~line of reasoning from Lemma~3.3.3 of the paper~\cite{Svg2} by Sauvageot.
Let $a \in {\mathcal A}$ with $a=a^*$, and let ${\mathcal B}$ be the unital $C^*$-subalgebra of~${\mathcal A}$
generated by~$a$.
Then ${\mathcal B}$ is commutative, and so we can identify ${\mathcal B}$ with $C(S)$ where~$S$ is the maximal ideal
space of ${\mathcal B}$, which we can identify with the spectrum, $\sigma(a)$, of~$a$, a~subset of ${\mathbb R}$.
Because ${\mathcal B}$ is commutative, its right action on $\Omega$ coming from the right action of~${\mathcal A}$ is
a~$*$-representation for the ${\mathcal D}$-valued inner product.
Let us denote this right $*$-representation of ${\mathcal B}$ by~$\rho$ and denote the left $*$-representation
by~$\lambda$.
Because these two representations commute, they combine to give a~representation, $\rho \otimes \lambda$, of ${\mathcal
B} \otimes {\mathcal B}=C(S \times S)$ on $\Omega$.
Because ${\mathcal D}$ is central, the representation $\rho \otimes \lambda$ commutes with the right action of
${\mathcal D}$ on $\Omega$, that is, it is a~$*$-representation into the $C^*$-algebra of endomorphisms of the right
Hilbert ${\mathcal D}$-module $\Omega$.

Now let~$p$ be a~polynomial of one variable with real coef\/f\/icients.
Let $\tilde p$ be the corresponding polynomial of two variables def\/ined by
\begin{gather*}
\tilde p(s,t)=(p(s)-p(t))/(s-t).
\end{gather*}
If~$p$ is the monomial $p(t)=t^n$ for some~$n$, then
\begin{gather*}
\tilde p(s, t)=s^{n-1}+s^{n-2}t+\dots+t^{n-1}.
\end{gather*}
Thus $\tilde p$ for a~general polynomial~$p$ will be a~linear combination of such expressions.
(The map $p \mapsto \tilde p$ is actually a~nice coproduct with a~Leibniz property.
See Proposition~3.11 of~\cite{HbT}.)

Since~$S$ is a~subset of ${\mathbb R}$ we can view $\tilde p$ as an element of $C(S \times S)$, and thus we can form the
operator $(\lambda \otimes \rho)(\tilde p)$.
Notice that
\begin{gather*}
{\partial}(a^n)=a^{n-1}({\partial} a)+a^{n-2}({\partial} a)a+\dots+({\partial} a)a^{n-1},
\end{gather*}
which we then recognize as being $((\lambda \otimes \rho)(\tilde p))({\partial} a)$ when~$p$ is the monomial $p(t)=t^n$.
It follows that for any polynomial~$p$ we have
\begin{gather*}
{\partial}(p(a))=((\lambda \otimes \rho)(\tilde p))({\partial} a).
\end{gather*}
Consequently
\begin{gather*}
\langle{\partial}(p(a)), {\partial}(p(a))\rangle_{\mathcal D} \leq (\|\tilde p\|_S)^2\langle{\partial} a, {\partial}
a\rangle_{\mathcal D},
\end{gather*}
where $\|\tilde p\|_S$ denotes the supremum norm of $\tilde p$ as an element of $C(S \times S)$.
But this supremum norm is exactly the Lipschitz constant, $\Lip(p)$, of~$p$ with respect to the restriction
to~$S$ of the metric from ${\mathbb R}$.
When we set $\|\omega\|^{\mathcal D}=\|\langle\omega, \omega\rangle_{\mathcal D}\|^{1/2}$, we see that we have
\begin{gather*}
\|{\partial}(p(a))\|^{\mathcal D} \leq \Lip(p) \|{\partial} a\|^{\mathcal D}.
\end{gather*}
Now let~$F$ be any ${\mathbb R}$-valued Lipschitz function on~$S$.
It has an extension~\cite{Wvr2} to a~Lipschitz function, $\hat F$, on any interval containing~$S$, such that $\|\hat
F\|_\infty=\|F\|_\infty$ and $\Lip(\hat F)=\Lip(F)$.
By the usual smoothing argument (e.g.\ as in the proof of Proposition~2.2 of~\cite{R4}),
$\hat F$ can be uniformly approximated on a~neighborhood of the
interval by functions with continuous f\/irst derivative, with no increase in the Lipschitz constant, and these functions
can in turn be approximated by polynomials uniformly and uniformly in the f\/irst derivative.
Thus~$F$ can be uniformly approximated by such polynomials, whose Lipschitz constants are no bigger than
$\Lip(F)$.
From what we found above for polynomials we thus obtain:

\begin{Proposition}
\label{prolip}
With notation as above, for any $a \in {\mathcal A}$ such that $a^*=a$ and for any ${\mathbb R}$-valued Lipschitz
function~$F$ on $S=\sigma(a)$ we have
\begin{gather*}
\langle{\partial}(F(a)), {\partial}(F(a))\rangle_{\mathcal D} \leq (\Lip(F))^2\langle{\partial} a, {\partial}
a\rangle_{\mathcal D}.
\end{gather*}
\end{Proposition}

Notice that for the case in which ${\mathcal D}={\mathbb C} 1_{\mathcal A}$ the above result holds for any trace, not
just for tracial states, as seen by scaling.

We now concentrate on the tracial case.

\begin{Definition}
Let~$\Gamma$ be a~CdC for~${\mathcal A}$, and let~$\tau$ be a~faithful trace on~${\mathcal A}$.
By the corresponding \emph{energy form} for $(\Gamma, \tau)$ we mean the ${\mathbb C}$-valued pre-inner product
${\mathcal E}_\Gamma$ on~${\mathcal A}$ def\/ined by
\begin{gather*}
{\mathcal E}_\Gamma(a, b)=\tau(\Gamma(a,b))
\end{gather*}
for all $a,b \in {\mathcal A}$.
If~$\Gamma$ is the CdC for a~Riemannian pre-metric $(\Omega, {\partial}, \langle\cdot, \cdot\rangle_{\mathcal A})$, then
we will also say that ${\mathcal E}_\Gamma$ is the energy form for $\Omega$ and~$\tau$.
When no confusion is likely we will often just write ${\mathcal E}$.
\end{Definition}

By normalizing the trace and then applying Proposition~\ref{prolip} we immediately obtain:

\begin{Corollary}
\label{corlip}
Let notation be as just above.
Then for any $a \in {\mathcal A}$ such that $a^*=a$, and for any ${\mathbb R}$-valued Lipschitz function~$F$ on
$\sigma(a)$ we have
\begin{gather*}
{\mathcal E}(F(a), F(a)) \leq (\Lip(F))^2 {\mathcal E}(a, a).
\end{gather*}
\end{Corollary}

We remark that this means that ${\mathcal E}$ is a~Dirichlet form according to Def\/inition~2.3 of~\cite{AHK}.
The def\/inition of a~Dirichlet form given before Theorem~3.3 of~\cite{DvL} has a~further requirement, but we will see
shortly that this further requirement is also satisf\/ied.
(Note that the Dirichlet forms we deal with here are all ``conservative'', i.e.~take value 0 if one of the entries is
$1_{\mathcal A}$.) In view of the terminology often used for Dirichlet forms, we set:

\begin{Definition}
We will say that a~form that satisf\/ies the property obtained in Corollary~\ref{corlip} is ``Markov'', or satisf\/ies the
``Markov property''.
\end{Definition}

\begin{Definition}
\label{defnorm}
Let notation be as above.
We def\/ine a~seminorm, $L_{\mathcal E}$, on~${\mathcal A}$ by
\begin{gather*}
L_{\mathcal E}(a)=({\mathcal E}(a,a))^{1/2}.
\end{gather*}
We will call $L_{\mathcal E}$ the \emph{energy norm} on~${\mathcal A}$ (even though it is a~seminorm).
\end{Definition}

\begin{Theorem}
\label{thmsemi}
Let ${\mathcal E}$ be the energy form for a~given CdC and faithful trace, and define $L_{\mathcal E}$ as above.
Then $L_{\mathcal E}$ is indeed a~seminorm, and it satisfies the Markov property that for any $a \in {\mathcal A}$ such
that $a^*=a$ and for any ${\mathbb R}$-valued Lipschitz function~$F$ on the spectrum $\sigma(a)$ we have
\begin{gather*}
L_{\mathcal E}(F(a)) \leq (\Lip(F)) L_{\mathcal E}(a).
\end{gather*}
Furthermore, $L_{\mathcal E}$ satisfies the Leibniz property that for any $a, b \in {\mathcal A}$ we have
\begin{gather*}
L_{\mathcal E}(ab) \leq L_{\mathcal E}(a) \|b\|+ \|a\| L_{\mathcal E}(b).
\end{gather*}
\end{Theorem}

\begin{proof}
$L_{\mathcal E}$ is a~seminorm because it is the seminorm for the ordinary pre-inner-product ${\mathcal E}$.
The Markov property follows immediately from Corollary~\ref{corlip}.

For the Leibniz property, if we have started with a~CdC~$\Gamma$, we apply Theorem~\ref{thmcst} to obtain the
corresponding Riemannian metric $(\Omega, {\partial}, \langle\cdot, \cdot \rangle_{\mathcal A}^\Gamma)$ for~${\mathcal A}$.
We use~$\tau$ to def\/ine an ordinary inner product on $\Omega$
by $\langle\omega, \omega'\rangle_\tau=\tau(\langle\omega, \omega'\rangle_{\mathcal A})$, with corresponding norm $\|\cdot\|_\tau$.
This norm is an~${\mathcal A}$-bimodule norm.
This means that
\begin{gather*}
\|\omega a\|_\tau\leq\|\omega\|_\tau \|a\|
\qquad
\text{and}
\qquad
\|a\omega\|_\tau\leq\|a\| \|\omega\|_\tau
\end{gather*}
for all $\omega \in \Omega$ and $a \in {\mathcal A}$.
The f\/irst of these inequalities uses the tracial property of~$\tau$ to calculate that
\begin{gather*}
\|\omega a\|_\tau^2=\tau(\langle\omega a, \omega a\rangle_{\mathcal A})=\tau(a^*\langle\omega,\omega\rangle_{\mathcal A} a)=\tau(a^*a\langle\omega,\omega\rangle_{\mathcal A}) \leq \|a^*a\|
\|\omega\|_\tau^2.
\end{gather*}
For the second of these inequalities, note that $\|a\|^2 1_{\mathcal A} -a^*a \geq 0$ in~${\mathcal A}$, and so has
a~positive square-root, say~$c$, in~${\mathcal A}$.
Then
\begin{gather*}
0 \leq \langle c\omega, c\omega\rangle_{\mathcal A}=\langle\omega, c^2\omega\rangle_{\mathcal A}=\|a\|^2\langle\omega,
\omega\rangle_{\mathcal A}-\langle a\omega, a\omega\rangle_{\mathcal A}.
\end{gather*}
We can now apply~$\tau$ to this to obtain the desired inequality.
Now ${\mathcal E}(a,b)=\langle{\partial} a, {\partial} b\rangle_\tau$, so that
\begin{gather*}
L_{\mathcal E}(a)=\|{\partial} a\|_\tau.
\end{gather*}
Because ${\partial}$ is a~derivation, the Leibniz inequality for $L_{\mathcal E}$ follows immediately.
(Notice the importance of the \emph{complete} positivity of~$\Gamma$ for this proof, because it leads to the inner
product on~$\Omega$.)
\end{proof}

We will see at the end of Section~\ref{secd2c} that basically the completely Markov property implies the completely Leibniz property.
The Markov property of standard deviation given in Theorem~3.9 of~\cite{R28} is a~special case of Theorem~\ref{thmsemi}
above, as can be seen from the discussion in Section~\ref{secsdev}.

The Leibniz property is important for the considerations in~\cite{R21, R17}, which is one reason that I began to study
the topic of this paper.
But there are many other seminorms that satisfy both the Markov and Leibniz conditions.
As examples when~${\mathcal A}$ is commutative, let $(X, \rho)$ be a~f\/inite metric space, let~$Z$ be def\/ined as in
Example~\ref{excomm}, and def\/ine~$c$ on~$Z$ by $c_{xy}=1/\rho(x,y)$.
Def\/ine~$L$ on ${\mathcal A}=C(X)$ by
\begin{gather*}
L(f)=\sup\{|f(x)-f(y)\,|\, c_{xy}: (x,y) \in Z\}.
\end{gather*}
This is the usual Lipschitz constant for~$f$ for the metric~$\rho$.
The metric from~$L$ on the state space $S({\mathcal A})$ of~${\mathcal A}$, when restricted to~$X$ identif\/ied with the
extreme points of $S({\mathcal A})$, is the original metric~$\rho$.
More generally, for any $p \geq 1$ def\/ine~$L$ by
\begin{gather*}
L(f)=\bigg(\sum\{|f(x)-f(y)|^p c_{xy}: (x,y) \in Z\}\bigg)^{1/p}.
\end{gather*}
It is easily seen that these seminorms satisfy both the Markov and Leibniz conditions.
A~further example is given in Example~\ref{exchd}.

Note that $L_{\mathcal E}$, as def\/ined in Def\/inition~\ref{defnorm}, need not be a~$*$-seminorm.
For example in the setting of Example~\ref{exqc} where $\Gamma_v(a, b)=[v,a]^*[v,b]$ we see that if $[v, v^*] \neq 0$
then~$v$ is in the null space of the form $\Gamma_v$ while $v^*$ is not, and so for a~faithful trace on~${\mathcal A}$
we will have $L_{\mathcal E}(v)=0 \neq L_{\mathcal E}(v^*)$.
But we can always obtain a~Markov and Leibniz $*$-seminorm by taking the max $L_{\mathcal E}(a)\vee L_{\mathcal
E}(a^*)$.

So far we have put no conditions on ${\partial}$.
But there is an important condition that is used in various papers concerning Dirichlet forms, which does put
a~condition on ${\partial}$, involving the chosen trace as well as the CdC.

\begin{Definition}
\label{deftreal}
For a~CdC $\Gamma$ on~${\mathcal A}$ and a~trace~$\tau$ on~${\mathcal A}$, we say that~$\Gamma$ is \emph{$\tau$-real}
if for all $a, b \in {\mathcal A}$ we have
\begin{gather*}
\tau(\Gamma(a^*, b^*))=\tau(\Gamma(b, a)).
\end{gather*}
If~$\Gamma$ comes from a~f\/irst-order dif\/ferential calculus with correspondence, $(\Omega, {\partial}, \langle\cdot,
\cdot\rangle_{\mathcal A})$, then we will say that ${\partial}$ is \emph{$\tau$-real} if~$\Gamma$ is~$\tau$-real.
\end{Definition}

It is clear that~$\Gamma$ is~$\tau$-real exactly if
\begin{gather*}
{\mathcal E}(a^*, b^*)={\mathcal E}(b, a),
\end{gather*}
and that $L_{\mathcal E}$ is then a~$*$-seminorm.

Following the terminology given in~\cite{DvL}, we accordingly set:

\begin{Definition}
\label{defereal}
Let ${\mathcal E}$ be a~sesquilinear ${\mathbb C}$-valued form on~${\mathcal A}$.
We say that ${\mathcal E}$ is \emph{real} if
\begin{gather*}
{\mathcal E}(a^*, b^*)={\mathcal E}(b, a),
\end{gather*}
for all $a, b \in {\mathcal A}$.
\end{Definition}

In view of the considerations above it would be reasonable to def\/ine a~Riemannian metric to be a~pair $(\Gamma, \tau)$
consisting of a~CdC and a~trace, and even require~$\tau$-reality.
But we do not adopt this def\/inition.

We now relate the above def\/initions to the setting of Example~\ref{exconst}.
When we want to emphasize viewing~${\mathcal A}$ equipped with the inner product $\langle a,b\rangle_\tau=\tau(a^*b)$, we
will often write~${\mathcal A}$ as $L^2({\mathcal A}, \tau)$.

\begin{Proposition}
\label{proreel}
Let~$N$ be a~self-adjoint operator on the Hilbert space $L^2({\mathcal A}, \tau)$ with the property that $N(1_{\mathcal
A})=0$, and define $\Gamma_N$ as in Example~{\rm \ref{exconst}}.
Then $\Gamma_N$ is~$\tau$-real in the sense that
\begin{gather*}
\tau(\Gamma_N(a^*, b^*))=\tau(\Gamma(b, a)    )
\end{gather*}
for all $a,b\in {\mathcal A}$.
\end{Proposition}

\begin{proof}
For $a, b \in {\mathcal A}$ we have
\begin{gather*}
\tau(\Gamma_N(a^*, b^*)-\Gamma_N(b, a))
=\tau(N(a)b^*-N(ab^*)+aN(b^*) -N(b^*)a+N(b^*a)-b^*N(a))
\\
\hphantom{}{}
=\tau(N(ab^*-b^*a))=0,
\end{gather*}
where for the last equality we have used the fact that for any $c \in {\mathcal A}$ we have
\begin{gather*}
\tau(N(c))=\langle1_{\mathcal A}, N(c)\rangle_\tau=\langle N(1_{\mathcal A}), c\rangle_\tau=0
\end{gather*}
because $N(1_{\mathcal A})=0$.
\end{proof}

\section[$\{L_n\}$ is ${\mathcal L}^2$-matricial]{$\boldsymbol{\{L_n\}}$ is $\boldsymbol{{\mathcal L}^2}$-matricial}

Let $(\Omega, \langle\cdot,\cdot \rangle_{\mathcal A})$ be a~(pre-)correspondence over~${\mathcal A}$.
For any $n \in {\mathbb N}$ we def\/ine a~correspondence, $(\Omega_n, \langle\cdot, \cdot\rangle_n)$, over $M_n({\mathcal
A})=M_n \otimes {\mathcal A}$ by
\begin{gather*}
\Omega_n=M_n \otimes \Omega=M_n(\Omega),
\end{gather*}
with the evident left and right actions of $M_n({\mathcal A})$, and with $M_n({\mathcal A})$-valued (pre-)inner product
determined by
\begin{gather*}
\langle\alpha \otimes \omega, \beta \otimes \omega'\rangle_n=\alpha^*\beta \otimes \langle\omega,
\omega'\rangle_{\mathcal A}
\end{gather*}
for $\alpha, \beta \in M_n$ and $\omega, \omega' \in \Omega$.
It is easily verif\/ied that the left and right actions of $M_n({\mathcal A})$ relate to this sequilinear form in the way
needed for a~pre-correspondence.

\begin{Lemma}
\label{lempos}
The above $M_n({\mathcal A})$-valued sesquilinear form $\langle\cdot, \cdot\rangle_n$ is positive.
Thus $(\Omega_n, \langle\cdot, \cdot\rangle_m)$ is a~pre-correspondence for $M_n({\mathcal A})$.
If $\langle\cdot, \cdot\rangle_{\mathcal A}$ is definite, then so is $\langle\cdot, \cdot\rangle_n$.
\end{Lemma}
\begin{proof}
Given $t=\sum\limits^m \alpha_j \otimes \omega_j$, we have
\begin{gather*}
\langle t, t\rangle_n=\sum\limits_{j,k}^m \alpha_j^*\alpha_k \otimes \langle\omega_j,\omega_k\rangle_{\mathcal A}.
\end{gather*}
But the matrix $C=\{\langle\omega_j,\omega_k\rangle_{\mathcal A}\}$ is a~non-negative element of $M_n({\mathcal A})$
because $M_n({\mathcal A})$ is faithfully represented on the Hilbert~${\mathcal A}$-module ${\mathcal A}^n$ (with inner
product $\langle\{a_j\}, \{b_k\}\rangle_{\mathcal A}=\sum a_j^*b_j$), and for any $\{a_j\} \in {\mathcal A}^n$ we have
\begin{gather*}
\langle\{a_j\}, C\{a_j\}\rangle_{\mathcal A}=\sum\limits_{j,k} a^*_j\langle\omega_j,\omega_k\rangle_{\mathcal A}
a_k=\Big\langle\sum\limits_j\omega_j a_j, \sum\limits_k\omega_k a_k\Big\rangle_{\mathcal A}\geq0.
\end{gather*}
Thus~$C$ can be expressed as $C=D^*D$ for some $D \in M_n({\mathcal A})$, and then, for $D=\{d_{jk}\}$, we have
\begin{gather*}
\langle t, t\rangle_n=\sum\limits_{j,k} \alpha_j^*\alpha_k \otimes \sum\limits_p d^*_{pj}d_{pk}=\sum\limits_p
\left(\sum\limits_j \alpha_j\otimes d_{pj}\right)^*\left(\sum\limits_k \alpha_k \otimes d_{pk}\right) \geq 0.
\end{gather*}

If $\langle t,t\rangle_n=0$, then by the generalized Cauchy Schwartz inequality of equation~\eqref{eqCS} we have
$\langle t,s\rangle_n=0$ for all $s \in M_n(\Omega)$.
Let $\{e_{jk}\}$ be the usual matrix units for $M_n$.
We can express~$t$ as $t=\sum\limits_{j,k} e_{jk} \otimes \omega_{jk}$.
For any f\/ixed $p$, $q$ and any $\omega' \in \Omega$ set $s=e_{pq} \otimes \omega'$.
Then
\begin{gather*}
0=\langle t, s\rangle_n=\sum\limits_{j, k} e_{jk}^*e_{pq} \otimes \langle\omega_{jk}, \omega'\rangle_{\mathcal A}
=\sum\limits_k e_{kq} \otimes \langle\omega_{pk}, \omega'\rangle_{\mathcal A}.
\end{gather*}
By the linear independence of the $e_{jk}$'s it follows that $\omega_{jk}=0$ for all $j$, $k$.
\end{proof}

Suppose now that $(\Omega, {\partial}, \langle\cdot, \cdot\rangle_{\mathcal A})$ is a~Riemannian pre-metric for
${\mathcal A}$.
Then we can def\/ine ${\partial}_n$ on $M_n({\mathcal A})$ with values in $\Omega_n$ by setting it on elementary tensors
to be
\begin{gather*}
{\partial}_n(\alpha \otimes a)=\alpha \otimes {\partial} a.
\end{gather*}
It is easily seen that ${\partial}_n$ is a~derivation.
Note further that the $({\partial} A)B$'s for $A, B \in M_n({\mathcal A})$ span~$\Omega_n$ because
\begin{gather*}
({\partial}_n(\alpha \otimes a))(\beta \otimes b)=\alpha\beta \otimes ({\partial} a)b
\end{gather*}
for $a,b \in {\mathcal A}$.
Thus:

\begin{Proposition}
For notation as above, $(\Omega_n, {\partial}_n, \langle\cdot, \cdot\rangle_n)$ is a~Riemannian pre-metric for
$M_n({\mathcal A})$.
\end{Proposition}

\begin{Proposition}
\label{procdn}
Let notation be as above, let~$\Gamma$ be the CdC for $(\Omega, {\partial}, \langle\cdot, \cdot\rangle_{\mathcal A})$,
and let $\Gamma_n$ be the CdC for $(\Omega_n, {\partial}_n, \langle\cdot, \cdot\rangle_n)$.
Then $\Gamma_n$ is given in terms of~$\Gamma$ by
\begin{gather*}
(\Gamma_n(A,B))_{jk}=\sum\limits_p \Gamma(a_{pj}, b_{pk}).
\end{gather*}
for $A, B \in M_n({\mathcal A})$ and $A=\{a_{jk}\}$ and similarly for~$B$.
\end{Proposition}
\begin{proof}
In terms of the usual matrix units $\{e_{jk}\}$ for $M_n$ we have
\begin{gather*}
\Gamma_n(A, B)=\langle{\partial}_n A, {\partial}_n B\rangle_n=\Big\langle\sum e_{ij} \otimes {\partial} a_{ij}, \sum
e_{k \ell} \otimes {\partial} b_{k \ell}\Big\rangle_n
\\
\phantom{\Gamma_n(A, B)}
=\sum e_{ji} e_{k\ell} \otimes \langle{\partial} a_{ij}, {\partial} b_{k \ell} \rangle_{\mathcal A}
=\sum\limits_{j\ell} e_{j \ell} \otimes \sum\limits_i \Gamma(a_{ij}, b_{i\ell}).\tag*{\qed}
\end{gather*}
\renewcommand{\qed}{}
\end{proof}

\begin{Corollary}
\label{corcdn}
Let ${\mathcal E}$ be the energy form for~$\Gamma$ and a~choice of faithful trace~$\tau$ on~${\mathcal A}$.
Let $\text{tr}_n$ be the usual un-normalized trace on $M_n$, and let $\tau_n=\text{tr}_n \otimes \tau$, the
corresponding trace on $M_n({\mathcal A})$.
Let ${\mathcal E}_n$ denote the corresponding energy form for the above Riemannian metric for $M_n({\mathcal A})$ and
the trace $\tau_n$.
Then ${\mathcal E}_n$ can be expressed in terms of ${\mathcal E}$ by
\begin{gather*}
{\mathcal E}_n(A, B)=\sum\limits_{jk} {\mathcal E}(a_{jk}, b_{jk}).
\end{gather*}
\end{Corollary}

We remark that the right-hand side above is just the usual pre-inner product that one puts on the tensor product of
pre-Hilbert spaces, applied to $M_n \otimes {\mathcal A}$ with the inner product on $M_n$ from its un-normalized trace
and with the pre-inner-product on~${\mathcal A}$ being ${\mathcal E}$, as can be seen by calculations similar to those
in the proof of Proposition~\ref{procdn}.
We can make the same def\/inition for general sesquilinear forms.

We can now apply the results of the previous section to obtain:

\begin{Corollary}
\label{matlei}
With notation as above, ${\mathcal E}_n$ is a~Markov form for all~$n$.
In other words, ${\mathcal E}$ is completely Markov.
\end{Corollary}

The def\/inition given in~\cite{DvL} for a~sequilinear form ${\mathcal E}$ to be a~Dirichlet form is slightly stronger
than that used by many authors, for in addition to the Markov condition (which in~\cite{DvL} is called the Lipschitz
condition), which concerns only self-adjoint elements, it also requires that (in our context) for any $a \in {\mathcal
A}$ one have
\begin{gather*}
{\mathcal E}(|a|, |a|) \leq {\mathcal E}(a, a).
\end{gather*}
In Proposition~3.4
of~\cite{DvL} it is shown that if ${\mathcal E}_2$ is Markov, then ${\mathcal E}$ is Dirichlet in
their sense.
From this and what we have shown above it is not dif\/f\/icult to obtain:

\begin{Corollary}
\label{matdir}
For ${\mathcal E}$ coming from a~Riemannian metric and a~trace on~${\mathcal A}$ as above, each ${\mathcal E}_n$ is
Dirichlet.
In other words, ${\mathcal E}$ is completely Dirichlet.
\end{Corollary}

As in Def\/inition~\ref{defnorm}, for each~$n$ we def\/ine a~seminorm, $L_{{\mathcal E}_n}$, on $M_n({\mathcal A})$ by
\begin{gather*}
L_{{\mathcal E}_n}(A)=({\mathcal E}_n(A,A))^{1/2}.
\end{gather*}

Let ${\mathcal V}$ be a~vector space over ${\mathbb C}$.
We let $M_n({\mathcal V})$ denote the vector space of $n \times n$ matrices with entries in ${\mathcal V}$.
Then $M_n({\mathcal V})$ is in an evident way a~bimodule over $M_n$.
We adapt to seminorms in the obvious way the def\/inition of Ruan~\cite{Rua} of an ${\mathcal L}^2$-matricial norm on
${\mathcal V}$.

\begin{Definition}
Let notation be as above.
A~sequence $\{\sigma_n\}$ in which $\sigma_n$ is a~seminorm on $M_n({\mathcal V})$ for each n, is said to be an
\emph{${\mathcal L}^2$-matricial seminorm} on ${\mathcal V}$ if it satisf\/ies the following two properties:
\begin{itemize}\itemsep=0pt
\item the normed-bimodule condition
\begin{gather*}
\sigma_n(\alpha V \beta) \leq \|\alpha\| \sigma_n(V) \|\beta\|
\end{gather*}
for all $\alpha, \beta \in M_n$ and all $V \in M_n({\mathcal V})$.
\item the ${\mathcal L}^2$-condition that if $V \in M_m({\mathcal V})$ and $W \in M_n{\mathcal V}$,
so that we can
form $V \oplus W=(
\begin{smallmatrix}
V & 0
\\
0 & W
\end{smallmatrix}
)$ in $M_{m+n}({\mathcal V})$, then
\begin{gather*}
\sigma_{m+n}(V \oplus W)=\big((\sigma_m(V))^2+(\sigma_n(W))^2\big)^{1/2}.
\end{gather*}
\end{itemize}
\end{Definition}

\begin{Theorem}
The sequence $\{L_{{\mathcal E}_n}\}$ defined immediately after Corollary~{\rm \ref{matdir}} is an ${\mathcal L}^2$-matricial
seminorm on~${\mathcal A}$ consisting of seminorms that satisfy the Markov and Leibniz conditions.
\end{Theorem}

\begin{proof}
The normed-bimodule condition is easily obtained from arguments similar to those used in the proofs of
Theorem~\ref{thmsemi} and Lemma~\ref{lempos}.
For the ${\mathcal L}^2$-condition, let $A, B \in M_m({\mathcal A})$ and $C, D \in M_n({\mathcal A})$.
From Corollary~\ref{corcdn} it follows immediately that
\begin{gather*}
{\mathcal E}_{m+n}(A \oplus C, B \oplus D)={\mathcal E}_m(A, B)+{\mathcal E}_n(C, D).
\end{gather*}
From this the ${\mathcal L}^2$-condition follows immediately.
\end{proof}

\section{The metric on the state space}
\label{metric}

Let $S({\mathcal A})$ denote the state space of the $C^*$-algebra~${\mathcal A}$.
Thus $S({\mathcal A})$ consists of the positive linear functionals,~$\mu$, on~${\mathcal A}$ having the property that
$\mu(1_{\mathcal A})=1$.
These are the natural generalization of probability measures to the non-commutative setting.

With notation as in the previous sections, let $L_{\mathcal E}$ be the seminorm on~${\mathcal A}$
from $(\Omega, {\partial}, \langle\cdot, \cdot\rangle_{\mathcal A}, \tau)$.
Of course, because ${\partial} 1_{\mathcal A}=0$ we have $L_{\mathcal E}(1_{\mathcal A})=0$.
Such seminorms are exactly the kind used in def\/ining quantum metric spaces~\cite{R4, R5, R21}, and they determine an
ordinary metric, $\rho_{\mathcal E}$, on $S({\mathcal A})$, def\/ined by
\begin{gather*}
\rho_{\mathcal E}(\mu, \nu)=\sup \{|\mu(a)-\nu(a)  |: L_{\mathcal E}(a) \leq 1\}.
\end{gather*}
This metric is a~generalization of the Monge-Kantorovich metric on the set of ordinary probability measures on a~compact
metric space.

\begin{Definition}
\label{defem}
The metric $\rho_{\mathcal E}$ def\/ined above on $S({\mathcal A})$ is called the \emph{energy metric} associated with
$(\Omega, {\partial}, \langle\cdot, \cdot\rangle_{\mathcal A}, \tau)$
\end{Definition}

This metric will take value $+\infty$ if there is an $a \in {\mathcal A}$ with $a \notin {\mathbb C} 1_{\mathcal A}$
such that $L_{\mathcal E}(a)=0$.
In this case we interpret this as meaning that our ``quantum space'' is not metrically connected.
For instance, for the class of examples described in Example~\ref{exnc} in which ${\partial}_v(a)=[v, a]$, so that
$\Gamma_v(a, b)=[v,a]^*[v,b]$ and ${\mathcal E}(a, a)=\tau([v,a]^*[v,a])$ and $L_{\mathcal E}(a)=(\tau([v,a]^*[v,a]))^{1/2}$, we obviously have $L_{\mathcal E}(v)=0$, so that $\rho_{\mathcal E}$ does take the value
$+\infty$.
On the other hand, if ${\partial}$ is given as a~sum of terms of the form $[v,a]$ for dif\/ferent~$v$'s (as in
Examples~\ref{exnc} and~\ref{excont}), it can easily happen that $\rho_{\mathcal E}$ takes only f\/inite values.

\begin{Definition}
\label{defcon}
With notation as above we say that~${\mathcal A}$ is \emph{metrically connected} for $(\Omega, {\partial}, \langle\cdot,
\cdot\rangle_{\mathcal A})$ if (for~$\Gamma$ the corresponding CdC) we have $\Gamma(a, a)=0$ only when $a \in {\mathbb
C} 1_{\mathcal A}$.
\end{Definition}

For the rest of this section we will assume that~${\mathcal A}$ is metrically connected unless the contrary is stated.
Then because~$\tau$ is faithful, $\rho_{\mathcal E}$ will take only f\/inite values (in our f\/inite-dimensional situation).

In order to make clear at what point we need the various properties satisf\/ied by ${\mathcal E}$, let us now assume for
a~while that ${\mathcal E}$ is an arbitrary pre-inner-product on~${\mathcal A}$ that satisf\/ies the property that
${\mathcal E}(a, a)=0$ exactly when $a \in {\mathbb C} 1_{\mathcal A}$.
We def\/ine $L_{\mathcal E}$ as in Def\/inition~\ref{defnorm}, and we def\/ine the metric $\rho_{\mathcal E}$ on $S({\mathcal
A})$ as in Def\/inition~\ref{defem}.
Let ${\tilde {\mathcal A}}={\mathcal A}/{\mathbb C} 1_{\mathcal A}$.
Then ${\mathcal E}$ drops to a~def\/inite inner product on ${\tilde {\mathcal A}}$, which we will again denote by
${\mathcal E}$.
Since~${\mathcal A}$ is f\/inite-dimensional, ${\tilde {\mathcal A}}$ equipped with ${\mathcal E}$ is a~Hilbert space.
Each element of ${\tilde {\mathcal A}}$ has a~unique representative in the null-space of~$\tau$ (i.e.~orthogonal to
$1_{\mathcal A}$ in $L^2({\mathcal A}, \tau)$), and so we can identify ${\tilde {\mathcal A}}$ with the null-space
of~$\tau$ when convenient.

Denote the dual vector space of~${\mathcal A}$ by ${\mathcal A}'$.
Then the dual vector space of ${\tilde {\mathcal A}}$ can be identif\/ied with the subspace ${\mathcal A}'^o$ consisting
of elements of ${\mathcal A}'$ that take value 0 on $1_{\mathcal A}$.
Note that if $\mu, \nu \in S({\mathcal A})$ then $\mu-\nu \in {\mathcal A}'^o$.
Any $\lambda \in {\mathcal A}'^o$ determines a~linear functional on the f\/inite-dimensional Hilbert space $\tilde
{\mathcal A}$, and thus is represented by an element, $h_\lambda$, of $\tilde {\mathcal A}$, so that
\begin{gather*}
\langle a, \lambda\rangle={\mathcal E}(h_\lambda, a)
\end{gather*}
for all $a \in {\mathcal A}$, where here $\langle a, \lambda\rangle$ denotes the usual pairing between ${\tilde {\mathcal A}}$
and its dual space.
(We let~$a$ also denote its image in ${\tilde {\mathcal A}}$.) Thus $\lambda \mapsto h_\lambda$ is a~conjugate linear
map from ${\mathcal A}'^o$ into ${\tilde {\mathcal A}}$.
It is clearly injective.
But ${\mathcal A}'^o$ and ${\tilde {\mathcal A}}$ have the same dimension, and so this map is also surjective.
When convenient we can view $h_\lambda$ as an element (unique) of~${\mathcal A}$ such that $\tau(h_\lambda)=0$.

Notice that $L_{\mathcal E}$ is just the Hilbert space norm on ${\tilde {\mathcal A}}$, and it determines a~dual norm,
$L_{\mathcal E}'$, on~${\mathcal A}'^o$, def\/ined by
\begin{gather*}
L_{\mathcal E}'(\lambda)=\sup \{|\langle a, \lambda\rangle|: L_{\mathcal E}(a) \leq 1\}.
\end{gather*}
For $\mu, \nu \in S({\mathcal A})$ we then see that
\begin{gather*}
\rho_{\mathcal E}(\mu, \nu)=L_{\mathcal E}'(\mu-\nu).
\end{gather*}
But in terms of $h_\lambda$ we have
\begin{gather*}
L_{\mathcal E}'(\lambda)=\sup \{|{\mathcal E}(h_\lambda, a)|: {\mathcal E}(a, a) \leq 1\}.
\end{gather*}
Since $({\tilde {\mathcal A}}, {\mathcal E})$ is a~Hilbert space, we know that the supremum on the right side is
attained at the unit vector pointing in the direction of $h_\lambda$, and consequently
\begin{gather*}
L_{\mathcal E}'(\lambda)=L_{\mathcal E}(h_\lambda).
\end{gather*}
That is, the surjective map $\lambda \mapsto h_\lambda$ is a~(conjugate linear) isometry from ${\mathcal A}'^o$ onto
${\tilde {\mathcal A}}$.
From this we obtain:
\begin{Proposition}
\label{promet}
Let notation be as above.
For $\mu, \nu \in S({\mathcal A})$ we have
\begin{gather*}
\rho_{\mathcal E}(\mu, \nu)=L_{\mathcal E}(h_{\mu-\nu}).
\end{gather*}
\end{Proposition}
Fix now some $\mu_0 \in S({\mathcal A})$, and for any $\mu \in S({\mathcal A})$ set
\begin{gather*}
\sigma(\mu)=h_{\mu-\mu_0}.
\end{gather*}
Note that for $\mu, \nu \in S({\mathcal A})$ we have
\begin{gather*}
\sigma(\mu)-\sigma(\nu)=h_{\mu-\nu}.
\end{gather*}
Let $\|\cdot \|_{\mathcal E}$ denote the norm on ${\tilde {\mathcal A}}$ from the inner product ${\mathcal E}$ on
${\tilde {\mathcal A}}$.
Then:

\begin{Theorem}
\label{thmisom}
With notation as above
\begin{gather*}
\rho_{\mathcal E}(\mu, \nu)=\|\sigma(\mu)-\sigma(\nu)\|_{\mathcal E}
\end{gather*}
for all $\mu, \nu \in S({\mathcal A})$.
Thus~$\sigma$ is an affine isometry from the convex metric space $(S({\mathcal A}), \rho_{\mathcal E})$ into the Hilbert
space $({\tilde {\mathcal A}}, {\mathcal E})$.
\end{Theorem}
Note that $\sigma(\mu_0)=0$, so that the choice of $\mu_0$ determines which element of $S({\mathcal A})$ is sent to
$0$ by~$\sigma$.

Jorgensen and Pearse~\cite{JrP} were the f\/irst to discover that, for resistance networks, i.e.~for~${\mathcal A}$
commutative, at least the set of extreme points of the state space equipped with the metric from the energy form, embeds
isometrically into a~Hilbert space.
See also Section~5.1 of~\cite{JrP2}, where the relationship with negative semidef\/inite forms is discussed.

In order to see the further consequences of requiring that ${\mathcal E}$ actually comes from a~non-commutative
Riemannian metric, we need to introduce the Laplace operator.

\section{The Laplace operator}

Because~${\mathcal A}$ is f\/inite-dimensional and~$\tau$ is a~faithful trace on~${\mathcal A}$, for any pre-inner-product
${\mathcal E}$ on~${\mathcal A}$ there will be a~unique positive linear operator,~$N$, on $L^2({\mathcal A}, \tau)$ such
that
\begin{gather*}
{\mathcal E}(a, b)=\langle a, N b\rangle_\tau
\end{gather*}
for all $a, b \in {\mathcal A}$, where $\langle a, b\rangle_\tau=\tau(a^*b)$.
(The word ``positive'' here refers to~$N$ as an operator on the Hilbert space $L^2({\mathcal A}, \tau)$, and \emph{not}
to how it relates to the order structure on the $C^*$-algebra~${\mathcal A}$.) If ${\mathcal E}$ comes from a~Riemannian
metric $(\Omega, {\partial}, \langle\cdot, \cdot\rangle_{\mathcal A})$, so that ${\mathcal E}(a, b)=\tau(\langle{\partial} a, {\partial} b\rangle_{\mathcal A})$, then it is appropriate to view~$N$ as
${\partial}^*{\partial}$, and make:

\begin{Definition}
If ${\mathcal E}$ comes from a~Riemannian metric and faithful trace, then we denote the operator~$N$ as above
by~$\Delta$ and we call it the \empty{Laplace operator} corresponding to the Riemannian metric and faithful trace.
\end{Definition}

We remark that this is contrary to the conventions frequently made that lead to the Laplace operator being a~negative
operator.

We now investigate the resulting special properties that~$\Delta$ will have.
We assume from now on that~$\Delta$ comes from a~Riemannian metric and a~faithful trace as above.
Notice f\/irst that $\Delta(1)=0$ because ${\mathcal E}(a,1)=0$ for all $a \in {\mathcal A}$.
Def\/ine $\Gamma_\Delta$ as in Theorem~\ref{thmsym} with~$\Delta$ playing the role of~$N$ there (but here with an
additional factor of 1/2 as is commonly done in the literature), so that
\begin{gather*}
\Gamma_\Delta(a, b)=(1/2)(\Delta(a^*)b-\Delta(a^*b)+a^*\Delta(b))
\end{gather*}
for all $a, b \in {\mathcal A}$.
It has the properties described in Theorem~\ref{thmsym}.
Furthermore, $\Gamma_\Delta$ is~$\tau$-real according to Proposition~\ref{proreel}.

\begin{Lemma}
\label{lemlap}
For any $a, b, c \in {\mathcal A}$ we have
\begin{gather*}
\langle c, \Gamma_\Delta(a, b)\rangle_\tau
=(1/2)\tau(\langle{\partial} a, ({\partial} b)c^*\rangle_{\mathcal A}+\langle c{\partial} b^*, {\partial} a^*\rangle_{\mathcal A}).
\end{gather*}
\end{Lemma}

\begin{proof}
For any $a, b, c \in {\mathcal A}$ we have
\begin{gather*}
2\langle c, \Gamma_\Delta(a, b)\rangle_\tau=\tau(c^*(\Delta(a^*)b-\Delta(a^*b)+a^*\Delta(b)))
\\
\phantom{2\langle c, \Gamma_\Delta(a, b)\rangle_\tau}
=\tau(bc^*\Delta(a^*))-\tau(c^*\Delta(a^*b))+\tau(c^*a^*\Delta(b))
\\
\phantom{2\langle c, \Gamma_\Delta(a, b)\rangle_\tau}
=\tau(\langle{\partial}(cb^*),{\partial}(a^*)\rangle_{\mathcal A})-\tau(\langle{\partial} c,
{\partial}(a^*b)\rangle_{\mathcal A})+\tau(\langle{\partial}(ac), {\partial} b\rangle_{\mathcal A})
\\
\phantom{2\langle c, \Gamma_\Delta(a, b)\rangle_\tau}
=\tau(\langle c{\partial} b^*, {\partial} a^*\rangle_{\mathcal A})+\tau(\langle{\partial} a,({\partial}b)c^*\rangle_{\mathcal A}),
\end{gather*}
where for the last equality we have used both the fact that ${\partial}$ is a~derivation and that several terms cancel,
and the tracial property of~$\tau$.
\end{proof}

\begin{Proposition}
\label{propos}
With notation as above, $\Gamma_\Delta(a, a) \geq 0$ for all $a \in {\mathcal A}$.
Consequently $\Gamma_\Delta$ is also symmetric.
\end{Proposition}

\begin{proof}
If we let $c=dd^*$ in the above Lemma, and then rearrange, we obtain
\begin{gather*}
2\langle d, \Gamma_\Delta(a, a)d\rangle_\tau
=\tau(\langle({\partial} a)d, ({\partial} a)d\rangle_{\mathcal A})+\tau(\langle d^*{\partial} a^*, d^* {\partial} a^*\rangle_{\mathcal A})\geq0.
\end{gather*}
Since the representation of~${\mathcal A}$ on $L^2({\mathcal A}, \tau)$ is faithful, it follows that $\Gamma_\Delta(a,
a) \geq 0$ as desired.
The symmetry of $\Gamma_\Delta$ follows by the usual arguments.
\end{proof}

Because
\begin{gather*}
\Delta(a^*)a-\Delta(a^*a)+a^*\Delta(a) \geq 0,
\end{gather*}
as seen in Proposition~\ref{propos}, it is equal to its adjoint, which is
\begin{gather*}
a^*(\Delta(a^*))^*-(\Delta(a^*a))^*+(\Delta(a))^* a.
\end{gather*}
Def\/ine $\Delta^\sharp$ by $\Delta^\sharp(a)=(\Delta(a^*))^*$ for $a \in {\mathcal A}$, as done just before Theorem~\ref{thmsym}.
Then the above expression is equal to $2\Gamma_{\Delta^\sharp}(a, a)$, and so by polarization we have
\begin{gather*}
\Gamma_{\Delta^\sharp}=\Gamma_\Delta.
\end{gather*}
This is consistent with the relation found in Theorem~\ref{thmsym} that ensures that $\Gamma_\Delta$ is symmetric.

Let $\tau_n$ be def\/ined on $M_n({\mathcal A})$ as in Corollary~\ref{corcdn}, and let ${\mathcal E}_n$ be the
corresponding energy form.
Then for any $A, B \in M_n({\mathcal A})$ we see from Corollary~\ref{corcdn} that
\begin{gather*}
{\mathcal E}_n(A, B)=\sum\limits_{jk} {\mathcal E}(a_{jk}, b_{jk})=\sum\limits_{jk}\langle a_{jk}, \Delta
b_{jl}\rangle_\tau=\langle A, (I_n \otimes \Delta)B\rangle_{\tau_n}.
\end{gather*}

Thus we obtain:

\begin{Proposition}
\label{prolpn}
The Laplacian for ${\mathcal E}_n$ is the operator $I_n \otimes \Delta$ on ${\mathcal L}^2(M_n({\mathcal A}, \tau_n)   )$.
In particular, $I_n \otimes \Delta$ comes from a~Riemannian metric.
\end{Proposition}

We denote $I_n \otimes \Delta$ by $\Delta_n$.
It will have the properties described above.
In particular, $\Gamma_{\Delta_n}$ is positive on $M_n({\mathcal A})$, and $\tau_n$-real by Proposition~\ref{proreel}.
But straight-forward calculations using Proposition~\ref{procdn} show that
\begin{gather*}
\Gamma_{\Delta_n}=(\Gamma_\Delta)_n.
\end{gather*}
Thus $\Gamma_\Delta$ is completely positive in the sense def\/ined just before Def\/inition~\ref{cdc}.
Then from Proposition~\ref{proncdc} we obtain:

\begin{Theorem}
Let~$\Gamma$ be the CdC for a~Riemannian metric on~${\mathcal A}$, let~$\tau$ be a~faithful trace on~${\mathcal A}$, and
let~$\Delta$ be the Laplace operator for the corresponding energy form.
Then $\Gamma_\Delta$ is a~CdC.
\end{Theorem}

In view of Proposition~\ref{proccn} we have:

\begin{Corollary}
With notation and assumptions as just above,~$\Delta$ is a~conditionally completely negative operator on~${\mathcal A}$.
\end{Corollary}

\begin{Proposition}
\label{proeng}
With notation as above we have
\begin{gather*}
\tau(\Gamma_\Delta(a, b))=(1/2)(\tau(\Gamma(a, b))+\tau(\Gamma(b^*, a^*))       )
=(1/2)({\mathcal E}(a,b)+{\mathcal E}(b^*, a^*))
\end{gather*}
for all $a, b \in {\mathcal A}$.
In particular, ${\mathcal E}(a,b)=\tau(\Gamma_\Delta(a, b))$ for all $a, b \in {\mathcal A}$ if and only if~$\Gamma$
is~$\tau$-real $($where ``$\tau$-real'' is defined in Definition~{\rm \ref{deftreal})}.
\end{Proposition}

\begin{proof}
We obtain the f\/irst assertion when we set $c=1_{\mathcal A}$ in the equation of Lemma~\ref{lemlap}.
The second assertion then follows from the def\/inition of ${\partial}$ being~$\tau$-real.
\end{proof}

We remark that the formula of Proposition~\ref{proeng} shows the virtue of including the factor of 1/2 that we
introduced in this section.

Returning to our case of CdC's, it follows that if~$\Gamma$ is not~$\tau$-real then it can not coincide with~$\Gamma_\Delta$.

Let us now set $\Delta^\natural=(1/2)(\Delta+\Delta^\sharp)$.
Then we see that again $\Gamma_{\Delta^\natural}=\Gamma_\Delta$.
But $\Delta^\natural$ has the further property that $(\Delta^\natural(a))^*=\Delta^\natural(a^*)$.
This means that it satisf\/ies the conditions given in~\cite{Lbd} for $-\Delta^\natural$ to be the generator of a~quantum
semigroup.
We can consider this semigroup to be the heat semigroup for our Riemannian metric, especially when ${\partial}$
is~$\tau$-real, in which case we also have ${\mathcal E}={\mathcal E}_{\Delta^\natural}$.

\begin{Proposition}
\label{proreal}
For notation as above,~$\Gamma$ is~$\tau$-real if and only if $(\Delta(a))^*=\Delta(a^*)$ for all $a \in {\mathcal
A}$, that is, $\Delta=\Delta^\sharp=\Delta^\natural$.
\end{Proposition}

\begin{proof}
Essentially by def\/inition, $\Gamma$ is~$\tau$-real exactly when ${\mathcal E}$ is real (Def\/inition~\ref{defereal}).
Then for all $a, b \in {\mathcal A}$ we have
\begin{gather*}
{\mathcal E}(a, b)=\langle a, \Delta(b)\rangle_\tau=\langle\Delta(a), b\rangle_\tau=\tau((\Delta(a))^*b),
\end{gather*}
while
\begin{gather*}
{\mathcal E}(b^*, a^*)=\langle b^*, \Delta(a^*)\rangle_\tau=\tau(b\Delta(a^*)  )
=\tau(\Delta(a^*)b).
\end{gather*}
Thus ${\mathcal E}(a,b)={\mathcal E}(b^*, a^*)$ for all $a, b \in {\mathcal A}$ if and only if $(\Delta(a))^*=\Delta(a^*)$ for all $a \in {\mathcal A}$.
\end{proof}

\begin{Example}
\label{excont}
We consider a~continuation of Examples~\ref{exqc} and~\ref{exgen}.
Let $v_1, \dots, v_m$ be elements of~${\mathcal A}$, and set
\begin{gather*}
\Gamma(a, b)=\sum [v_j, a]^*[v_j, b].
\end{gather*}
Each term is a~CdC since it comes from a~Riemannian metric as in Examples~\ref{exmulti} and~\ref{exmore}.
Since sums of CdC's are again CdC's, it follows that~$\Gamma$ is a~CdC.
For any faithful trace~$\tau$ we have
\begin{gather*}
{\mathcal E}(a, b)=\sum \tau((a^*v_j^*-v_j^*a^*)[v_j, b])
=\tau\left(a^*\sum[v_j^*, [v_j, b]]\right).
\end{gather*}
It follows that the corresponding Laplace operator~$\Delta$ is def\/ined by
\begin{gather*}
\Delta(b)=\sum[v_j^*, [v_j, b]].
\end{gather*}
(See equation~(2.1) of~\cite{zhb}  for a~somewhat special case of this.) Let us see when~$\Gamma$ is~$\tau$-real.
According to Proposition~\ref{proreal} it suf\/f\/ices to determine when $\Delta^\sharp=\Delta$.
It is easily seen that
\begin{gather*}
\Delta^\sharp(b)=\sum[v_j, [v_j^*, b]].
\end{gather*}
Then
\begin{gather*}
\Delta(b)-\Delta^\sharp(b)=\sum\big([v_j^*, [v_j, b]]-[v_j, [v_j^*, b]]\big)
\\
\phantom{\Delta(b)-\Delta^\sharp(b)}
=\sum\big([v_j^*, [v_j, b]]+[[v_j^*, b], v_j]\big)=\bigg[\sum [v_j^*, v_j], b\bigg],
\end{gather*}
where the last equality comes from the Jacobi identity as in Example~\ref{exqc}.
It follows that~$\Gamma$ is~$\tau$-real exactly if $\sum [v_j^*, v_j]$ is in the center of~${\mathcal A}$.
But it is clear that for every trace~$\tau$ on~${\mathcal A}$ we have $\tau(\sum [v_j^*, v_j])=0 $ and that $\sum
[v_j^*, v_n]$ is self-adjoint.
Since it is in the center, and every state on the center extends to a~trace on~${\mathcal A}$ in our f\/inite-dimensional
situation, this is equivalent to being 0.
Thus~$\Gamma$ is~$\tau$-real exactly if
\begin{gather*}
\sum [v_j^*, v_j]=0.
\end{gather*}
This last condition is exactly the ``detailed balance'' condition of Proposition~6.9 of~\cite{DvL}.
\end{Example}

We remark that when this example is compared to Example~\ref{exgen} we see that there are generators~$N$ of quantum
dynamical semigroups on some $C^*$-algebras~${\mathcal A}$ for which there is no faithful trace~$\tau$ on~${\mathcal A}$
such that the CdC $\Gamma_N$ is~$\tau$-real.
This suggests that one should consider faithful states that are not tracial, as considered in~\cite{Cpr1, Cpr2}.
But I have not investigated that direction.

The formula in Lemma~\ref{lemlap} suggests a~further condition that can be required of a~Riemannian metric and trace.
This condition essentially appears already in Section~1.2 of~\cite{Ptr}, where Peterson calls it ``real''.
Since we already are using ``$\tau$-real'', I prefer to use the term ``$\tau$-balanced'':

\begin{Definition}
\label{defbal}
Let $(\Omega, {\partial}, \langle\cdot, \cdot \rangle_{\mathcal A})$ be a~Riemannian metric on~${\mathcal A}$.
For~$\tau$ a~faithful trace on~${\mathcal A}$, we say that this Riemannian metric is \emph{$\tau$-balanced} if
\begin{gather*}
\tau(\langle{\partial} a, ({\partial} b)c\rangle_{\mathcal A})=\tau(\langle c^*{\partial} b^*, {\partial}
a^*\rangle_{\mathcal A})
\end{gather*}
for all $a, b, c \in {\mathcal A}$.
If~$\Gamma$ is the CdC for the Riemannian metric, then the~$\tau$-balanced condition can be stated in terms of~$\Gamma$
as
\begin{gather*}
\tau(\Gamma(ab, c))=\tau(\Gamma(c^*, b^*)a^*)+\tau(b^*\Gamma(a, c))
\end{gather*}
for all $a, b, c \in {\mathcal A}$.
\end{Definition}

Because $\tau(\langle{\partial} a, ({\partial} b) c^*\rangle_{\mathcal A})=\langle c, \Gamma(a, b)\rangle_\tau$,
it follows immediately from Lemma~\ref{lemlap} that:

\begin{Theorem}
\label{thmbal}
Let $(\Omega, {\partial}, \langle\cdot, \cdot \rangle_{\mathcal A})$ be a~Riemannian metric on~${\mathcal A}$, and
let~$\Gamma$ be its CdC.
Let~$\tau$ be a~faithful trace on~${\mathcal A}$, and let~$\Delta$ be the corresponding Laplace operator.
Then $\Gamma=\Gamma_\Delta$ if and only if~$\Gamma$ is~$\tau$-balanced.
\end{Theorem}

It is clear that being~$\tau$-balanced is a~stronger condition than being~$\tau$-real.
The following example shows that it is in fact a~strictly stronger condition.

\begin{Example}
We consider CdC's of the form discussed in Example~\ref{exqc}, that is, of the form $\Gamma(a,b)=[v,a]^*[v,b]$ for
some $v \in {\mathcal A}$.
We want~$\Gamma$ to be~$\tau$-real, and so from Example~\ref{excont} we see that~$v$ must commute with $v^*$, that is,
be normal.

By Def\/inition~\ref{defbal}, in order for~$\Gamma$ to be~$\tau$-balanced we must have
\begin{gather*}
0=\tau(\langle c[v,b^*], [v,a^*]\rangle_{\mathcal A}-\langle[v,a],[v,b]c^*\rangle_{\mathcal A})=\tau([b,v^*]c^*[v,a^*]-[a^*, v^*][v,b]c^*).
\end{gather*}
Because this is true for all $c \in {\mathcal A}$, and~$\tau$ is faithful and tracial, this is equivalent to the
requirement that
\begin{gather}
0=[v,a^*][b,v^*]-[a^*, v^*][v,b].
\label{bal}
\end{gather}
for all $a, b \in {\mathcal A}$.
This is satisf\/ied if~$v$ is self-adjoint, so for it to fail we must choose~$v$ to be normal but not self-adjoint.

To continue, we now take~${\mathcal A}$ to be $M_n({\mathbb C})$ for some~$n$, with its usual trace.
Since~$v$ is to be normal, we can assume that it is a~diagonal matrix.
Since the requirement of equation~\eqref{bal} involves only commutators, we can always change~$v$ by adding a~scalar
multiple of $1_{\mathcal A}$.
Thus if~$v$ has only two eigenvalues, we can assume that one of those eigenvalues is~0.
Then~$v$ is a~scalar multiple of a~self-adjoint matrix, and again the requirement of equation~\eqref{bal} is satisf\/ied.
Thus we must assume that~$v$ has at least 3 eigenvalues, but we can assume that one of those eigenvalues is~0.
We can also multiply~$v$ by a~scalar, and so assume that another of the eigenvalues is~1.
By conjugating~$v$ by a~permutation matrix we can then assume that the f\/irst 3 diagonal entries of~$v$ are~$1$,~$\alpha$,~$0$,
where~$\alpha$ is some non-real complex number.
View~${\mathcal A}$ as acting on ${\mathbb C}^n$ in the usual way.
Then let~$b$ be the element of~${\mathcal A}$ that takes the f\/irst standard basis vector to the second, takes the second
to the third, and sends all other standard basis vectors to~0.
A~simple calculation shows that with $a=b^*$ we have
\begin{gather*}
[v,a^*][b,v^*]-[a^*, v^*][v,b] \neq 0.
\end{gather*}
Thus~$\Gamma$ is not~$\tau$-balanced.
\end{Example}

We now relate the Laplace operator~$\Delta$ to the metric $\rho_{\mathcal E}$ on the state space.
For this we assume that~${\mathcal A}$ is metrically connected (Def\/inition~\ref{defcon}), so that the kernel of~$\Delta$
is exactly ${\mathbb C} 1_{\mathcal A}$.
Let ${\mathcal A}_0=\{a \in {\mathcal A}: \tau(a)=0\}$, so that when viewed as a~subspace of $L^2({\mathcal A},
\tau)$ it is exactly the orthogonal complement of ${\mathbb C} 1_{\mathcal A}$.
Thus~$\Delta$ carries ${\mathcal A}_0$ into itself and is invertible there.
Later, when we write $\Delta^{-1}$, it is to be interpreted as an operator on ${\mathcal A}_0$.
In the evident way ${\mathcal A}_0$ can be identif\/ied with $\tilde {\mathcal A}={\mathcal A}/1_{\mathcal A}$.
Because ${\mathcal A}_0$ is a~Hilbert space when equipped with the inner product from $L^2({\mathcal A}, \tau)$, we can
identify ${\mathcal A}'^0$ conjugate linearly with ${\mathcal A}_0$ itself.
Accordingly we change our earlier conventions in the usual way, and for $\lambda \in {\mathcal A}_0$ we def\/ine
$h_\lambda \in {\mathcal A}_0$ such that
\begin{gather*}
\langle\lambda, a\rangle_\tau={\mathcal E}(h_\lambda, a)
\end{gather*}
for all $a \in {\mathcal A}$, so that $\lambda \mapsto h_\lambda$ is linear, in contrast to our convention in Section~\ref{metric}.
But ${\mathcal E}(h_\lambda, a)=\langle\Delta h_\lambda, a\rangle_\tau$ for all~$a$, and so $\lambda=\Delta h_\lambda$, or
\begin{gather*}
h_\lambda=\Delta^{-1} \lambda.
\end{gather*}
Thus
\begin{gather*}
{\mathcal E}(h_\lambda, h_\lambda)=\langle\lambda, h_\lambda\rangle_\tau=\langle\lambda, \Delta^{-1}
\lambda\rangle_\tau.
\end{gather*}
In view of Proposition~\ref{promet}, we can express $\rho_{\mathcal E}$ in terms of~$\Delta$ by:

\begin{Proposition}
\label{metlap}
With notation as above, for any $\mu, \nu \in S({\mathcal A})$, with $\mu-\nu$ viewed as an element of $L^2({\mathcal
A}, \tau)$, we have
\begin{gather*}
\rho_{\mathcal E}(\mu, \nu)=\langle\mu-\nu, \Delta^{-1}(\mu-\nu)\rangle_\tau^{1/2}.
\end{gather*}
\end{Proposition}

We will make use of this in Section~\ref{strn} for the commutative case.

\section{The commutative case}
\label{comm}

We now examine the special case in which~${\mathcal A}$ is a~commutative (and f\/inite-dimensional) $C^*$-algebra.
Let~$X$ be its maximal ideal space.
Then~$X$ is a~f\/inite set, and we can and will identify~${\mathcal A}$ with $C(X)$.
We begin by determining all possible CdC's on~${\mathcal A}$.

\begin{Theorem}
\label{thmch}
As in Example~{\rm \ref{exmore}}, for $Z=\{(x,y) \in X\times X: x\neq y\}$, any non-negative func\-tion~$c$ on~$Z$ defines
a~CdC for~${\mathcal A}$, by
\begin{gather*}
\Gamma(f,g)(y)=\sum\limits_{x,x \neq y} (\bar f(x)-\bar f(y))(g(x)-g(y))c_{xy},
\end{gather*}
for $x \in X$.
Conversely, every CdC for~${\mathcal A}$ is of this form, and thus there is a~bijection between the set of CdC's and the
set of~$c$'s.
Furthermore, every CdC $\Gamma$ for~${\mathcal A}$ satisfies the extra condition that
\begin{gather*}
\Gamma(f^*, g^*)=\Gamma(g, f)
\end{gather*}
for all $f, g \in {\mathcal A}$, so that~$\Gamma$ is~$\tau$-real for any trace~$\tau$ on~${\mathcal A}$.
Furthermore, if $c_{yx}=c_{xy}$ for all $(x,y) \in Z$ then~$\Gamma$ is~$\tau$-balanced for any trace.
\end{Theorem}

\begin{proof}
Note that in the above sum def\/ining~$\Gamma$ we do not need values for $c_{yy}$ for any~$y$.
As suggested in Example~\ref{exmore}, it is easily seen that, given~$c$, the above formula gives a~CdC (which clearly
satisf\/ies the extra condition), and a~direct calculation verif\/ies the statement about being $\tau$-balanced.

Thus we must prove the converse.
So let~$\Gamma$ be some given CdC for~${\mathcal A}$.
For each $x \in X$ let $\delta_x$ be the usual ``delta-function'' at~$x$.
Since the $\delta_x$'s form a~basis for~${\mathcal A}$, the constants
\begin{gather*}
\gamma_{pq}^y=\Gamma(\delta_p, \delta_q)(y)
\end{gather*}
for $p, q, y \in X$ completely determine~$\Gamma$.
Because~$\Gamma$ is symmetric, we see that
\begin{gather*}
\gamma_{qp}^y=\bar\gamma_{pq}^y
\end{gather*}
for all $p, q, y \in X$.
Because~$\Gamma$ is positive we see that
\begin{gather*}
\gamma_{pp}^y \geq 0
\end{gather*}
for all $p, y \in X$.
Because $\Gamma(1, f)=0$ for all $f \in {\mathcal A}$ we see that
\begin{gather*}
\sum\limits_p \gamma_{pq}^y=0
\end{gather*}
for all $q, y \in X$.
Finally, we must examine the consequences of the $*$-representation condition, equation~\eqref{eqstar}.
Let $y, w, p, q \in X$.
Then from equation~\eqref{eqstar} we obtain
\begin{gather*}
\Gamma(\delta_w\delta_p,\delta_q)(y)-\Gamma(\delta_p, \delta_w\delta_q)(y)=\delta_p(y)\Gamma(\delta_w, \delta_q)(y)
- \Gamma(\delta_p,\delta_w)(y)\delta_q(y).
\end{gather*}
Suppose that $p \neq q$.
On setting $w=p$ we see that:
\begin{gather*}
\text{if}
\qquad
p\neq q,
\qquad
y\neq p
\qquad
\text{and}
\qquad
y\neq q
\qquad
\text{then}
\qquad
\gamma_{pq}^y=0	,
\end{gather*}
whereas:
\begin{gather*}
\text{if}
\qquad
p\neq q
\qquad
\text{and}
\qquad
y=q
\qquad
\text{then}
\qquad
\gamma_{pq}^q=-\gamma_{pp}^q.
\end{gather*}
In particular, since $\gamma_{pp}^q \geq 0$ by the positivity of~$\Gamma$, we see that $\gamma_{pq}^y \in {\mathbb R}$
for all $p, q, y \in X$.

Then for $f, g \in {\mathcal A}$ and $y \in X$ we have
\begin{gather*}
\Gamma(f, g)(y)=\sum\limits_{p, q} \bar f(p)g(q)\gamma_{pq}^y
\\
\phantom{\Gamma(f, g)(y)}
=\sum\limits_{p,p\neq y} \bar f(p)\left(\sum\limits_q g(q)\gamma_{pq}^y\right)+\bar f(y)\sum\limits_q g(q)\gamma_{yq}^y
\\
\phantom{\Gamma(f, g)(y)}
=\sum\limits_{p,p\neq y} \bar f(p)(g(q)\gamma_{pp}^y +g(y)\gamma_{py}^y)-\bar f(y)\sum\limits_q
g(q)\gamma_{qq}^y
\\
\phantom{\Gamma(f, g)(y)}
=\sum\limits_{p,p\neq y} \bar f(p)(g(p)-g(y))\gamma_{pp}^y-\bar f(y)\sum\limits_p g(p)\gamma_{pp}^y.
\end{gather*}
We can add to this $0=\bar f(y)g(y)\sum\limits_p\gamma_{pp}^y$ (because $\gamma_{pp}^y=- \gamma_{py}^y$) to obtain:
\begin{gather*}
\Gamma(f,g)(y)=\sum\limits_{p,p\neq y} (\bar f(p)-\bar f(y))(g(p)-g(y))\gamma_{pp}^y.
\end{gather*}
Accordingly, if we set $c_{py}=\gamma_{pp}^y$ for all $p$, $y$ with $p \neq y$ we obtain the desired formula.
Note that for each $p, y \in X$ with $p \neq y$ we have
\begin{gather*}
c_{py}=\gamma_{pp}^y=\Gamma(\delta_p, \delta_p)(y),
\end{gather*}
which is non-negative by assumption, as needed.
\end{proof}

We remark that, as suggested by the remarks about commutative algebras near the beginning of Section~\ref{energy}
leading to the proof of Proposition~\ref{prolip}, the above CdC will satisfy the Markov condition.
This is easily seen directly.
If $f \in {\mathcal A}$, and if~$F$ is even a~${\mathbb C}$-valued function def\/ined on~$\sigma(f)$ (which is the range
of~$f$), then
\begin{gather*}
\Gamma(F\circ f, F\circ f)(y)=\sum\limits_{x}|F(f(x))-F(f(y))|^2 c_{xy} \leq(\Lip(F))^2 \Gamma(f,
f)(y).
\end{gather*}

Given $c_{xy}$'s as above, we can view~$X$ as consisting of the nodes of a~directed graph such that there is an edge
from~$x$ to~$y$ exactly if $c_{xy} \neq 0$; and we can consider the $c_{xy}$'s to be weights on the edges, recognizing
that the weights on the two edges joining in opposite directions two given nodes need not be equal.

For traditional reasons that will be discussed below, we now introduce a~factor of 1/2 into the formula for~$\Gamma$,
much as we did just before Lemma~\ref{lemlap}.
Thus from now on we assume that
\begin{gather*}
\Gamma(f,g)(y)=(1/2)\sum\limits_{x \neq y} (\bar f(x)-\bar f(y))(g(x)-g(y))c_{xy},
\end{gather*}
We have not yet chosen a~trace on~${\mathcal A}$.
But counting measure is implicit in the formula above for~$\Gamma$, and it is anyway almost a~canonical choice.
Thus we will choose (integration against) counting measure as our trace~$\tau$.
The corresponding energy form ${\mathcal E}$ is then given by
\begin{gather*}
{\mathcal E}(f, g)=(1/2)\sum\limits_{x,y} (\bar f(x)-\bar f(y))(g(x)-g(y))c_{xy}.
\end{gather*}
(So we are in the ``jump part'' of Example 4.20 of~\cite{CpS}.) Notice that the part of the summand involving~$f$
and~$g$ is even in~$x$ and~$y$ (i.e.~unchanged under exchanging~$x$ and~$y$).
Consequently its sum with the odd part of~$c$ will be~0.
Thus we can replace~$c$ by $\tilde c$ def\/ined by $\tilde c_{xy}=(c_{xy}+c_{yx})/2$.
Of course $\tilde c$ will def\/ine a~dif\/ferent CdC.
But the two CdC's will give the same energy form, and for the following considerations it is the energy form that we
study.
So we assume now that
\begin{gather*}
c_{xy}=c_{yx}.
\end{gather*}
With this condition, the graph whose nodes are the elements of~$X$ and whose edge-weights are given by~$c$ can be
interpreted exactly as a~resistance network, with~$c$ specifying the conductances between the various nodes~\cite{JrP2,Kgm}.
It is in this way that our Riemannian metrics together with trace correspond, when~${\mathcal A}$ is commutative,
exactly to resistance networks.

We now sketch the usual development for resistance networks~\cite{JrP2, Kgm}, since we need it for our discussion of the
metric on the state space.
We f\/irst determine the corresponding Laplacian.
For this purpose it is convenient to def\/ine (as, for example, in Def\/inition~1.9 of~\cite{JrP2})
two operators on~${\mathcal A}$, which when they are viewed as operators on $L^2({\mathcal A}, \tau)$ are self-adjoint operators.
In def\/ining these operators, we will assume that $c_{xx}=0$ for all $x \in X$.
The f\/irst of these operators, the ``transfer operator''~$T$, is an ``integral operator'' def\/ined by
\begin{gather*}
(Tf)(x)=\sum\limits_y c_{xy}f(y).
\end{gather*}
For the second of these operators,~$C$, def\/ine f\/irst a~function, $\hat c$, on~$X$ by
\begin{gather*}
\hat c(x)=\sum\limits_y c_{xy}=\sum\limits_y c_{yx}.
\end{gather*}
We let~$C$ be the operator of pointwise multiplication by $\hat c$.
Then
\begin{gather*}
2\Gamma(f,g)(y)=\sum\limits_{p} (\bar f(p)-\bar f(y))(g(p)-g(y))c_{py}
\\
\phantom{2\Gamma(f,g)(y)}
=(T(\bar f g))(y)-\bar f(y)(Tg)(y)-(T\bar f)(y)g(y)+(C(\bar f g))(y)
\\
\phantom{2\Gamma(f,g)(y)}
=-((C-T)(\bar f g))(y)+((C-T)\bar f)(y) g(y)+\bar f(y)((C-T)g)(y)
\\
\phantom{2\Gamma(f,g)(y)}
=2\Gamma_{C-T}(f, g)(y).
\end{gather*}
Consequently $\Gamma=\Gamma_{C-T}$, and from the above calculation we also see that
\begin{gather*}
{\mathcal E}(f, g)=\tau(\Gamma(f,g))=\langle f, (C-T)g\rangle_\tau.
\end{gather*}
Thus the Laplace operator for~$\Gamma$ is given by
\begin{gather*}
\Delta=C-T,
\end{gather*}
and
\begin{gather*}
\Gamma=\Gamma_\Delta,
\end{gather*}
consistent with Theorem~\ref{thmbal}.
The specif\/ic formula for the Laplace operator~$\Delta$ can be written as
\begin{gather*}
(\Delta f)(x)=\sum\limits_y (f(x)-f(y))c_{xy}
\end{gather*}
Note that in the literature the Laplace operator is often taken to be the negative of the above expression, so that it
is a~non-positive operator.

We see from the above formulas why it is common to introduce a~factor of 1/2 in the def\/inition of ${\mathcal E}$.
There is another compelling reason for introducing a~factor of 1/2, namely that when the system is interpreted as
a~resistance network, and~$f$ is interpreted as giving voltages that are applied to the various nodes, the rate of
dissipation of energy caused by the resulting current is given by the earlier ${\mathcal E}$ divided by 2, basically
because the earlier formula for ${\mathcal E}$ double-counts the edges.

Let the $\delta_x$'s now be viewed as elements of $L^2({\mathcal A}, \tau)$, so that they form an orthonormal basis for
$L^2({\mathcal A}, \tau)$.
For any $z \in X$ we have
\begin{gather*}
(\Delta(\delta_x))(z)=((C-T)(\delta_x))(z)=\hat c(z)\delta_x(z)-c_{zx}.
\end{gather*}
Consequently, for $x, y \in X$ with $x \neq y$ we have
\begin{gather*}
{\mathcal E}(\delta_x, \delta_y)=\langle\delta_x,\Delta\delta_y\rangle_\tau=-c_{xy}.
\end{gather*}
Notice that this relation fails if $x=y$, but that instead we have
\begin{gather*}
{\mathcal E}(\delta_x, \delta_x)=\langle\delta_x,\Delta\delta_x\rangle_\tau=\hat c(x).
\end{gather*}

It is easily seen that even when some of the $c_{xy}$ are negative the corresponding form ${\mathcal E}$ def\/ined as
above can still be non-negative.

\begin{Theorem}
\label{thmneg}
Assume that the form ${\mathcal E}$ is non-negative, but do not require that the $c_{xy}$ are all non-negative.
Then ${\mathcal E}$ satisfies the Markov condition if and only if $c_{xy} \geq 0$ for all $x \neq y$.
\end{Theorem}

\begin{proof}
We recall here the simple argument (found, for example, in the proof of Proposition~2.1.3 of~\cite{Kgm}).
Suppose that for some given $x$, $y$ we have $c_{xy} <0$.
Set $ f=\delta_x-r\delta_y$ for some $r \in {\mathbb R}_{>0}$.
Def\/ine~$F$ on ${\mathbb R}$ by $F(t)=t$ if $t\geq 0$ and $F(t)=0$.
Then $\Lip(F)=1$, and $F\circ f=\delta_x$.
Thus if ${\mathcal E}$ were to satisfy the Lipschitz condition we should have ${\mathcal E}(\delta_x, \delta_x) \leq
{\mathcal E}(f, f)$.
But when we expand the sum for ${\mathcal E}(f,f)$ we obtain
\begin{gather*}
{\mathcal E}(\delta_x, \delta_x)+2rc_{xy}+r^2{\mathcal E}(\delta_y, \delta_y)
\end{gather*}
Since $c_{xy}$ is strictly negative, it is clear that we can choose a~positive~$r$ small enough that ${\mathcal E}(f,f)<
{\mathcal E}(d_x, d_x)$.
The converse assertion follows from Corollary~\ref{corlip}.
\end{proof}

\section{Resistance distance}
\label{strn}

As mentioned in the introduction, I have been puzzled about the nature of the ``resistance distance'' for a~resistance
network ever since I wrote Section~12 of~\cite{R5}.
In this section we will arrive at an answer that I consider satisfactory.
Let me mention that resistance distance has seen use in chemistry (e.g.~\cite{Ivn, KPRT,Kln, ZhY} and their
references), and even in evolution~\cite{McR}.

As is usual, we say that a~graph is ``connected'' if it is not the disjoint union of two non-empty subsets,~$A$ and~$B$,
such that there is no edge between any point of~$A$ and any point of~$B$.
Maximal connected subsets of~$X$ are called its ``connected components''.
These concepts are of importance to us because it is easily seen that for a~resistance network we have $L_{\mathcal
E}(f)=0$ for some $f \in {\mathcal A}$ exactly if~$f$ is constant on the connected components of~$X$.
Thus it is only when~$X$ itself is connected that we have the property that if $L_{\mathcal E}(f)=0$ then $f \in
{\mathbb C} 1_{\mathcal A}$, so that the corresponding metric on the state space takes only f\/inite values.
It is easily seen that the properties of ${\mathcal E}$ and related objects can be obtained by treating each connected
component separately.
Consequently, we will assume for the rest of this section that~$X$ is connected.
Since this depends on the choice of $c_{xy}$'s we will tend to say ``metrically connected''.

We will now develop the standard ideas about harmonic functions, which for our context go all the way back to the
seminal paper~\cite{BgD}.
(I have not noticed a~useful way to develp a~theory of ``harmonic functions'' in the non-commutative setting.) The
material in the next paragraphs, through Theorem~\ref{thmrmt} is well-known.
See for example Section~2.1 of~\cite{Kgm}.

Let~$\Delta$ be the Laplace operator for the given choice of $c_{xy}$'s for~$X$ (metrically connected).

\begin{Definition}
For given $f \in C(X)$ and $x \in X$ we say that~$f$ is \emph{harmonic} at~$x$ if $\Delta(f)(x)=0$.
\end{Definition}

If~$f$ is interpreted as an application of voltages at the points of~$X$, then being harmonic at~$x$ means that no
current is being inserted (or extracted) at~$x$.

If~$f$ is harmonic at~$x$, then
\begin{gather*}
0=\sum\limits_y(f(x)-f(y))c_{xy},
\end{gather*}
so that
\begin{gather*}
f(x)=\sum\limits_y f(y)\left(c_{xy}\Big/\left(\sum\limits_w c_{xw}\right)\right)=\sum\limits_y f(y)
(c_{xy}/\hat c(x)).
\end{gather*}
Notice that $y \mapsto c_{xy}/\hat c(x)$ is a~probability distribution on the set of points of~$X$ that share an edge with~$x$.
Thus $f(x)$ is a~weighted average of the values of~$f$ on the neighbors of~$x$.
Here we make essential use of the fact that the $c_{xy}$'s are non-negative.

For any subset~$Y$ of~$X$ let $\overline Y=Y \cup \{z \in X: c_{yz} > 0~\text{for~some}~y \in Y\}$.
Notice that this is not a~true closure operation.

\begin{Theorem}[the maximum principle]
Let~$Y$ be a~subset of~$X$, and assume that~$Y$ is connected $($for the restriction of~$c$ to~$Y)$.
Let $f \in C(X)$, with $\bar f=f$, and suppose that~$f$ is harmonic at all points of~$Y$.
Let $m=\max\{f(w): w \in \overline Y\}$.
If there is a~$y \in Y$ such that $f(y)=m$, then~$f$ is constant on~$Y$.
\end{Theorem}

\begin{proof}
Let $W=\{y \in Y: f(y)=m\}$, and suppose that~$W$ is not empty.
Let $w \in W$.
Because $f(w)$ is the weighted average of its values on $\overline{\{w\}}$, it must take value~$m$ at all points of
$\overline{\{w\}}$.
It follows that $\overline{\{w\}} \subseteq W$.
Because~$Y$ is connected, it follows easily that $W=Y$.
\end{proof}

Notice that by considering $-f$ we obtain a~similar statement about the minimum of~$f$ on $\overline Y$.

As discussed before Proposition~\ref{metlap}, because~$X$ is metrically connected, the restriction of~$\Delta$ to
${\mathcal A}_0$ is an invertible operator on ${\mathcal A}_0$, where here ${\mathcal A}_0$ consists of the
functions~$f$ such that $\tau(f)=0$.
For any distinct $p, q \in X$ the function $\delta_p-\delta_q$ is in ${\mathcal A}_0$.
Let $h_{pq}=\Delta^{-1}(\delta_p-\delta_q)$.
Then $\Delta(h_{pq})=\delta_p-\delta_q$, and so $h_{pq}$ is harmonic on the complement of $\{p,q\}$.
On applying the maximam principle to the dif\/ferent components of $X \setminus \{p, q\}$ we see that $h_{pq}$ must take
its maximum and minimum values on $\{p, q\}$.
Now
\begin{gather*}
1=\Delta(h_{pq})(p)=\sum\limits_x (h_{pq}(p)-h_{pq}(x))c_{xp},
\end{gather*}
and from this it is clear that it must be at~$p$ that $h_{pq}$ must take its maximum value, and so its minimum value
at~$q$.
That is, for all $x \in X$
\begin{gather}
h_{pq}(q) \leq h_{pq}(x) \leq h_{pq}(p).
\label{min}
\end{gather}
Note that according to Proposition~\ref{metlap} we have
\begin{gather}
\rho_{\mathcal E}(p, q)=\big(\langle\delta_p-\delta_q, h_{pq}\rangle_\tau\big)^{1/2}
=(h_{pq}(p)-h_{pq}(q))^{1/2}
\label{cons}
\end{gather}
(where here~$\tau$ is counting measure).

\begin{Definition}
Def\/ine $\rho_r$ on $X\times X$ by $\rho_r(p, p)=0$ for all $p \in X$, and
\begin{gather*}
\rho_r(p, q)=h_{pq}(p)-h_{pq}(q)
\end{gather*}
for $p, q \in X$ with $p \neq q$.
We call $\rho_r$ the \emph{resistance metric} on~$X$ (for the given $c_{xy}$'s).
\end{Definition}

\begin{Theorem}
\label{thmrmt}
The resistance metric $\rho_r$ is indeed a~metric.
\end{Theorem}

\begin{proof}
It is clear that $\rho_r$ is symmetric (because $h_{qp}=- h_{pq}$), and that $\rho_r (x,y)=0$ exactly if $x=y$.
We must show that it satisf\/ies the triangle inequality.
So let points $n$, $p$, $q$ of~$X$ be given.
By the linearity of~$\Delta$ we have $h_{pq}=h_{pn}+h_{nq}$, while $h_{nq}(p) \leq h_{nq}(n)$ and $h_{pn}(q) \geq
h_{pn}(n)$ by the inequalities~\eqref{min} from the maximum principle.
Thus
\begin{gather*}
\rho_r(p, q)=h_{pq}(p)-h_{pq}(q)=h_{pn}(p)+h_{nq}(p)-h_{pn}(q)-h_{nq}(q)
\\
\phantom{\rho_r(p, q)}
\leq h_{pn}(p)+h_{nq}(n)-h_{pn}(n)-h_{nq}(q)=\rho_r(n, p)+\rho_r(n, q).\tag*{\qed}
\end{gather*}
\renewcommand{\qed}{}
\end{proof}

Now this is strange, because from equation~\eqref{cons} we see that
\begin{gather*}
\rho_r(p,q)=h_{pq}(p)-h_{pq}(q)=(\rho_{\mathcal E}(p, q))^2,
\end{gather*}
and usually the square of a~metric is not a~metric.
Since $\rho_{\mathcal E}$ is def\/ined on the whole state space~$S({\mathcal A})$, not just on its extreme points, it is
natural to ask whether $\rho_{\mathcal E}^2$ is a~metric on all of~$S({\mathcal A})$.
We will now see that this is not the case, so the resistance metric is of a~quite dif\/ferent nature than the energy
metric.

Recall that Theorem~\ref{thmisom} tells us that $S({\mathcal A})$, equipped with the energy metric, is isometrically
embedded in a~Hilbert space.
So let us determine what kinds of subsets of a~Hilbert space have the property that when the square of the Hilbert space
metric is restricted to them the result is again a~metric.

\begin{Proposition}
\label{prohilb}
Let~$X$ be a~subset of a~Hilbert space ${\mathcal H}$, and let~$d$ be the restriction to~$X$ of the metric on ${\mathcal
H}$ that comes from its inner product.
Then $d^2$ is a~metric on~$X$ if and only if~$X$ has the property that for all $x, y, z \in X$ we have
\begin{gather*}
\Reoperator(\langle x-y, z-y\rangle) \geq 0.
\end{gather*}
\end{Proposition}

\begin{proof}
From the def\/inition of a~metric, $d^2$ is a~metric on~$X$ if and only if for all $x, y, z \in X$ we have
\begin{gather*}
\langle x-z, x-z\rangle \leq \langle x-y,x-y\rangle+\langle y-z,y-z\rangle.
\end{gather*}
On expanding these inner products and canceling some terms, we f\/ind that $0 \leq \Reoperator(\langle x-y, z-y\rangle)$ as
desired.
\end{proof}

The above proposition is related to von Neumann's embedding theorem (for which see Appendix~A.1 of~\cite{JrP2}), but the
condition of the above proposition is necessarily much stronger than the negative semi-def\/initeness of von Neumann's theorem.

The real part of an inner product is itself an inner product when the Hilbert space is considered to be a~vector space
over ${\mathbb R}$.
In view of the above proposition, we now assume that we have a~Hilbert space ${\mathcal H}$ over ${\mathbb R}$, and
a~subset~$X$ of it having the above property that for all $x, y, z \in X$ we have
\begin{gather*}
0 \leq \langle x-y, z-y\rangle.
\end{gather*}
Let~$K$ denote the closed convex hull of~$X$ in ${\mathcal H}$.
Given a~$z \in X$, choose any $y \in X$ with $y \neq z$, and let $\phi_{z;y}$ be the linear functional on ${\mathcal H}$
def\/ined by
\begin{gather*}
\phi_{z;y}(w)=\langle w, z-y\rangle.
\end{gather*}
for all $w \in {\mathcal H}$.
Then the above inequality implies that every $x \in X$ lies in the half-space of ${\mathcal H}$ def\/ined by
$\phi_{z;y}(w) \geq \langle y, z-y\rangle$.
Thus~$K$ must lie in this half-space.
Notice in particular that~$z$ itself lies strictly in the interior of this half-space.
From all of this it is clear that each $x \in X$ is an extreme point of~$K$.
Furthermore, our assumed inequality says that the angles between the lines from one point of~$X$ to any two other points
of~$X$ are acute or right angles.

In our original situation in which~$X$ is a~f\/inite subset of $S({\mathcal A})$, it follows that~$X$ consists of exactly
all the extreme points of $S({\mathcal A})$.
Thus we obtain:

\begin{Proposition}
Let~$\Delta$ be the Laplace operator for a~connected resistance network on a~finite set~$X$, let ${\mathcal A}=C(X)$
as earlier, and let $\rho_r$ be defined as above on $S({\mathcal A})$.
Then no subset of $S({\mathcal A})$ that properly contains~$X$ has the property that the restriction of $\rho_r$ to it
is a~metric.
\end{Proposition}

In view of this, I consider the resistance metric to be of a~quite dif\/ferent nature than the metrics on state spaces
(such as the energy metric) that I have been studying.
To me this is a~satisfactory resolution to my puzzlement about the resistance metric.
In particular, I believe that the ``free resistance'' def\/ined in equation~(4.73) of Section~4.9 of~\cite{JrP2} does not
satisfy the triangle inequality.

It is natural to ask whether the property given in Proposition~\ref{prohilb} characterizes the operators that arise as
Laplace operators for connected resistance networks.
A~counter-example is given in Exercise~2.5 of~\cite{Kgm}.
The basic idea is quite simple.
In Section~\ref{metric} we discussed forms more general than the energy forms coming from resistance networks, and saw
that they too embed $S({\mathcal A})$ into Hilbert spaces.
Now if~$X$ has at least~4 points, one can f\/ind conductances on~$X$ such that~$X$ is connected but one of the
conductances is~0, while $\Reoperator(\langle x-y, z-y\rangle) > 0$ for all $x, y, z \in X$.
Then we can make the~0 conductance slightly negative (so we no longer have a~resistance network), in such a~way that the
angles are changed so little that we still have $\Reoperator(\langle x-y, z-y\rangle) > 0$ for the new inner product.

\section{From Dirichlet forms to CdC's}
\label{secd2c}

We have seen that from a~CdC and a~faithful trace~$\tau$ on~${\mathcal A}$ we obtain an energy form ${\mathcal E}$.
This energy form is completely Markov in the sense that each ${\mathcal E}_n$ is Markov.
One of the central theorems of the general theory of Dirichlet forms is that conversely each completely positive and
completely Markov form on $L^2({\mathcal A}, \tau)$ comes from a~CdC (and so in our setting comes from a~Riemannian
metric and trace, and has a~corresponding quantum dynamical semigroup, etc).
In the inf\/inite-dimensional case substantial technical assumptions are needed in order to prove this.
Here we will just treat the f\/inite-dimensional case.
We will use this case later in Section~\ref{secquot}.
We assume that ${\mathcal E}$ is real (Def\/inition~\ref{defereal}), as is usually required for Dirichlet forms.
The main theorem of this section is thus:

\begin{Theorem}
\label{thmdir}
Let ${\mathcal E}$ be a~sesquilinear form on~${\mathcal A}$ which is real, completely positive and completely Markov,
and has $1_{\mathcal A}$ in its null-space.
Let~$\tau$ be a~faithful trace on~${\mathcal A}$, and let~$N$ be the operator on $L^2({\mathcal A}, \tau)$ determined by
\begin{gather*}
{\mathcal E}(a,b)=\langle a, Nb\rangle_\tau.
\end{gather*}
Then $\Gamma_N$ is a~CdC, and
\begin{gather*}
{\mathcal E}(a, b)=\tau(\Gamma_N(a, b))
\end{gather*}
for all $a, b \in {\mathcal A}$.
\end{Theorem}

\begin{proof}

It is the Markov property that is the key to the proof, but it seems to be a~bit tricky to extract useful information
from it.
We will follow the usual method, as given for example following Theorem~2.7 of~\cite{AHK}.

Since ${\mathcal E}$ is real, the argument given in the proof of Proposition~\ref{proreal}
shows that $(N(a))^*=N(a^*)$ for all $a \in {\mathcal A}$, that is, $N^\sharp=N$.
Thus according to Proposition~\ref{proncdc} we only need to show that~$\Gamma_N$ is completely positive.
A~simple calculation using Corollary~\ref{corcdn} shows that because~${\mathcal E}$ is real~${\mathcal E}_n$ also is
real.
But at f\/irst we will not use the ``completely'' aspects of ${\mathcal E}$.

Notice that~$N$ is a~positive operator on $L^2({\mathcal A},\tau)$ such that $N(1_{\mathcal A})=0$.
Since~$N$ is positive, $I+N$ is invertible, where~$I$ denotes the identity operator on the Hilbert space
$L^2({\mathcal A},\tau)$.

\begin{Lemma}[key lemma]
\label{lemkey}
Let $R=(I+N)^{-1}$.
Then, for any $a \in {\mathcal A}$ for which $a \geq 0$ we have $R(a) \geq 0$ as an element of~${\mathcal A}$, and $R(a)
\leq \|a\|$.
\end{Lemma}

\begin{proof}
Def\/ine a~positive sesquilinear form, ${\mathcal F}$, on $L^2({\mathcal A}, \tau)$ by
\begin{gather*}
{\mathcal F}(b,a)=\langle b, (I+N)a\rangle_\tau,
\end{gather*}
and note that ${\mathcal F}$ is def\/inite and that ${\mathcal F}(b, Ra)=\langle b, a\rangle_\tau$.
It is easily calculated that
\begin{gather*}
{\mathcal E}(b,b)+\|b-a\|^2_\tau={\mathcal F}(b-Ra, b-Ra)+\langle(I-R)a,a\rangle_\tau
\end{gather*}
From this it is clear that for f\/ixed~$a$ the left hand side has a~unique minimum when $b=Ra$.
Thus we f\/ind, for f\/ixed~$a$, that
\begin{gather*}
{\mathcal E}(Ra, Ra)+\|Ra-a\|^2_\tau < {\mathcal E}(b,b)+\|b-a\|^2_\tau
\end{gather*}
for all $b \in {\mathcal A}$ such that $b \neq Ra$.

Suppose now that $a=a^*$.
Because ${\mathcal E}$ is real (Def\/inition~\ref{defereal}),~$N$ preserves the involution, and thus~$R$ will also.
Consequently $Ra$ is self-adjoint.
Let~$F$ be an ${\mathbb R}$-valued Lipschitz function on ${\mathbb R}$ with $\Lip(F) \leq 1$.
Then on setting $b=F(Ra)$, which is well-def\/ined because $Ra$ is self-adjoint, we obtain
\begin{gather}
{\mathcal E}(Ra, Ra)+\|Ra-a\|^2_\tau < {\mathcal E}(F(Ra), F(Ra))+\|F(Ra)-a\|^2_\tau
\label{ineq}
\end{gather}
if $F(Ra) \neq Ra$.

We apply this last result in the following way.
Let~$F$ be def\/ined by $F(t)=\max(t, 0)$.
Notice that $\Lip(F)=1$.
Observe that for any~$a$ such that $a^*=a$ we have $F(a)=a^+$, the positive part of~$a$.
Suppose that $b \in {\mathcal A}$ with $b^*=b$.
We denote the negative part of~$a$ by $a^-$ and similarly for~$b$.
Because $a^+$ and $a^-$ are orthogonal to each other in $L^2({\mathcal A}, \tau)$ since $a^+a^-=0$, and similarly
for~$b$, we have
\begin{gather*}
\|b-a\|_\tau^2-\|b^+-a^+\|_\tau^2=\|b^--a^-\|_\tau^2+2\Reoperator(\langle b^+, a^-\rangle_\tau+\langle b^-, a^+\rangle_\tau).
\end{gather*}
Because $b^+$ and $a^-$ are positive and~$\tau$ is tracial, we have $\langle b^+, a^-\rangle_\tau \geq 0$.
Similarly $\langle b^-, a^+\rangle_\tau \geq 0$.
It follows that
\begin{gather}
\|b^+-a^+\|_\tau^2\leq\|b-a\|_\tau^2.
\label{posdif}
\end{gather}
(This is a~special case of the fact that for any real Lipschitz function~$F$ we would have
\begin{gather*}
\|F(b)-F(a)\|_\tau \leq \Lip(F)\|b-a\|_\tau,
\end{gather*}
but this general case seems not to have an easy proof.
See Lemma~2.2 of~\cite{AHK} or Proposition~2.5 of~\cite{DvL}.)

Now assume that~$a$ is positive, so that $F(a)=a$.
From inequality~\eqref{posdif} we then obtain
\begin{gather*}
\|F(Ra)-a\|_\tau=\|(Ra)^+-a\|_\tau \leq \|Ra-a\|_\tau.
\end{gather*}
Using this in the right side of inequality~\eqref{ineq} and cancelling, we obtain for our~$F$
\begin{gather*}
{\mathcal E}(Ra,Ra) < {\mathcal E}(F(Ra), F(Ra))
\end{gather*}
if $F(Ra) \neq Ra$.
But $\Lip(F) \leq 1$ and ${\mathcal E}$ is assumed to be Markov, which means that
\begin{gather*}
{\mathcal E}(Ra,Ra) \geq {\mathcal E}(F(Ra),F(Ra)).
\end{gather*}
This contradiction implies that $F(Ra)=Ra$, so that $Ra \geq 0$.
This proves the f\/irst assertion of the proposition.

To prove the second assertion, for any $r \in {\mathbb R}$ def\/ine $F_r$ by $F_r(t)=\min(t, r)$, so again
$\Lip(F_r)=1$.
Then $F_r(t)=-(t-r)^+ +r$, so that, using the fact that~${\mathcal A}$ is unital and writing~$r$ for $r1_{\mathcal
A}$, we have
\begin{gather*}
\|F_r(b)-F_r(a)\|_\tau=\|-(a-r)^+\! +r+(b-r)^+\! -r\|_\tau=\|(a-r)^+\! -(b-r)^+\|_t \leq \|b-a\|_\tau
\end{gather*}
by inequality~\eqref{posdif}.
Now let $r=\|a\|$.
Then for any~$b$ such that $b^*=b$ we have $F_r(b)=b$ exactly if $b \leq \|a||$.
But $F_r(a)=a$, and so
\begin{gather*}
\|F_r(b)-a\|_\tau \leq \|b-a\|_\tau.
\end{gather*}
When we use this in inequality~\eqref{ineq} in the way done above, we see that $F_r(Ra)=Ra$, so that $Ra \leq \|a\|$,
as desired.
\end{proof}

Notice that we did not need the full force of the Markov property.
We only needed it for the two functions $F(t)=\max(t, 0)$ and $F(t)=\min(t, \|a\|)$.
But in the context of Theorem~\ref{thmdir} the full Markov property will then be a~consequence.

Suppose now that ${\mathcal E}$ is completely Markov, where ${\mathcal E}_n$ is def\/ined by the same formula as given
just before Proposition~\ref{prolpn}.
Then as in that proposition, the ``Laplacian'' for ${\mathcal E}_n$ is just $I_n \otimes N$.
If we multiply ${\mathcal E}$ be any $t \in {\mathbb R}_{>0}$ the various ${\mathcal E}_n$'s will also be multiplied
by~$t$, and the resulting forms will still be Markov.
The corresponding ``Laplacians'' will also be multiplied by~$t$.
We can apply the key Lemma~\ref{lemkey} to all of them.
For this purpose we set $R_t=(I+tN)^{-1}$, and $R_t^{(n)}=I_n \otimes R_t$.
By the key Lemma~\ref{lemkey} each $R_t^{(n)}$ is a~positive operator on $M_n({\mathcal A})$ that is contractive on positive
elements.
This says that each $R_t$ is a~completely positive operator on~${\mathcal A}$.
Also each of these operators will carry the identity element to itself because $N(1_{\mathcal A})=0$.
The basic inequality for completely positive operators that we already used in Example~\ref{exsg} implies that for any
$A \in M_n({\mathcal A})$ we have
\begin{gather}
R_t^{(n)}(A^*A)-R_t^{(n)}(A^*)R_t^{(n)}(A)\geq0.
\label{matin}
\end{gather}
From this it easily follows that $\|R_t^{(n)}\|=1$.
Note also that $R_0$ is well-def\/ined and $R_0=I$, and that $R_t$ is actually well-def\/ined also for negative~$t$'s in
a~neighborhood of 0.

In our f\/inite-dimensional situation the function $t \mapsto R_t^{(n)}$ is clearly dif\/ferentiable for each~$n$.
Notice that the left-hand side of inequality~\eqref{matin} has value 0 at $t=0$.
It follows, much as in Example~\ref{exsg}, that the derivative at $t=0$ of the left-hand side is non-negative.
But the derivative $R'_t$ is $-N(I+tN)^{-2}$, so that $R'_0=-N$, and similarly for $R^{(n)}$.
Thus the derivative at $t=0$ of inequality~\eqref{matin} gives
\begin{gather*}
-N(a^*a)-(-N(a^*)a-a^*N(a))\geq0,
\end{gather*}
that is, for $\Gamma_N$ def\/ined as in Example~\ref{exconst} we have $\Gamma_N(a, a) \geq 0$ for all $a \in {\mathcal
A}$.
In the same way we f\/ind that $\Gamma_N^{(n)}(A,A) \geq 0$ for all $A \in M_n({\mathcal A})$, so that $\Gamma_N$ is
completely positive.
From Proposition~\ref{proncdc} it follows that $\Gamma_N$ is a~CdC.
This completes the proof of Theorem~\ref{thmdir}.
\end{proof}

On combining the above result with Theorem~\ref{thmsemi} and Corollary~\ref{matlei} we f\/ind that the completely Markov
property of ${\mathcal E}$ implies the completely Leibniz property of $L_{\mathcal E}$.

\section{Dirac operators}

In this section we show how to construct a~Hodge--Dirac operator for a~Riemannian metric, once a~trace has been chosen.
We assume throughout that~${\mathcal A}$ is a~f\/inite-dimensional $C^*$-algebra, that $(\Omega, {\partial},
\langle\cdot,\cdot\rangle_{\mathcal A})$ is a~Riemannian metric for~${\mathcal A}$, and that~$\tau$ is a~faithful trace
on~${\mathcal A}$.
We def\/ine the corresponding Hodge--Dirac operator in analogy with Def\/inition~9.24 of~\cite{GVF}.
For the case in which~${\mathcal A}$ is commutative our Dirac operator is essentially the operator used by Davies in
Theorem~4.6 of~\cite{Dav}.
We begin by def\/ining an ordinary inner product, $\langle\cdot,\cdot\rangle_\tau$, on $\Omega$ by
\begin{gather*}
\langle\omega, \omega'\rangle_\tau=\tau(\langle\omega, \omega'\rangle_{\mathcal A}).
\end{gather*}
Because $\Omega$ is f\/inite-dimensional, it is a~Hilbert space for this inner product.
We denote this Hilbert space by $L^2(\Omega, \tau)$.
Then ${\partial}$ can be viewed as an operator from $L^2({\mathcal A}, \tau)$ to $L^2(\Omega, \tau)$.
We denote the adjoint of this operator, going from $L^2(\Omega, \tau)$ to $L^2({\mathcal A}, \tau)$, by ${\partial}^*$.
Let
\begin{gather*}
{\mathcal H}=L^2({\mathcal A}, \tau) \oplus L^2(\Omega, \tau).
\end{gather*}
We def\/ine the operator~$D$ on ${\mathcal H}$ by
\begin{gather*}
D=\begin{pmatrix}
0 & {\partial}^*
\\
{\partial} & 0
\end{pmatrix}
.
\end{gather*}
We view~${\mathcal A}$ as acting on ${\mathcal H}$ by means of its left actions on~${\mathcal A}$ and $\Omega$.
We now calculate much as in the proof of Theorem~4.6 of~\cite{Dav}.
For $a \in {\mathcal A}$ and $(
\begin{smallmatrix}
b
\\
\omega
\end{smallmatrix}
) \in {\mathcal H}$ we have
\begin{gather*}
[D, a]
\begin{pmatrix}
b
\\
\omega
\end{pmatrix}
=\begin{pmatrix}
{\partial}^*(a\omega)-a{\partial}^*\omega
\\
{\partial}(ab)-a{\partial} b
\end{pmatrix}
.
\end{gather*}
By the Leibniz rule
\begin{gather}
\|{\partial}(ab)-a{\partial} b\|_\tau^2=\|({\partial} a)b\|_\tau^2=\tau(\langle({\partial} a)b, ({\partial}a)b\rangle_{\mathcal A})
=\tau(b^*\langle{\partial} a, {\partial} a\rangle_{\mathcal A} b)
\nonumber
\\
\phantom{\|{\partial}(ab)-a{\partial} b\|_\tau^2}
 \leq \|\langle{\partial} a, {\partial} a\rangle_{\mathcal A}\|_\infty \tau(b^*b),
\label{leib}
\end{gather}
where for emphasis we here denote the $C^*$-norm of~${\mathcal A}$ by $\|\cdot \|_\infty$.
Furthermore, for any $c \in {\mathcal A}$ we have
\begin{gather*}
\langle c, {\partial}^*(a\omega)-a{\partial}^* \omega\rangle_\tau
=\langle a^*{\partial} c, \omega\rangle_\tau-\langle{\partial} (a^*c), \omega\rangle_\tau=- \langle({\partial} a^*)c, \omega\rangle_\tau,
\end{gather*}
and by the Cauchy--Schwarz inequality
\begin{gather*}
|\tau(\langle({\partial} a^*)c, \omega\rangle_{\mathcal A})| \leq (\tau(\langle({\partial} a^*)c, ({\partial}
a^*)c\rangle_{\mathcal A}))^{1/2} (\tau(\langle\omega, \omega\rangle_{\mathcal A}))^{1/2},
\end{gather*}
while
\begin{gather*}
\tau(\langle({\partial} a^*)c, ({\partial} a^*)c\rangle_{\mathcal A})=\tau(c^*\langle{\partial} a^*, {\partial}
a^*\rangle_{\mathcal A} c) \leq \|\langle{\partial} a^*, {\partial} a^*\rangle_{\mathcal A}\|_\infty \tau(c^*c).
\end{gather*}
It follows that
\begin{gather*}
\|{\partial}^*(a\omega)-a{\partial}^* \omega\|_\tau \leq \|\langle{\partial} a^*, {\partial} a^*\rangle_{\mathcal
A}\|_\infty^{1/2} \|\omega\|_\tau.
\end{gather*}
From this and the calculation~\eqref{leib} above we f\/ind that
\begin{gather*}
\|[D, a]\| \leq \max \big\{\|\langle{\partial} a, {\partial} a\rangle_{\mathcal A}\|_\infty^{1/2},\|\langle{\partial}
a^*, {\partial} a^*\rangle_{\mathcal A}\|_\infty^{1/2} \big\}.
\end{gather*}
Notice that $a \mapsto \|\langle{\partial} a, {\partial} a\rangle_{\mathcal A}\|_\infty^{1/2}$ is not in general stable
under the involution, whereas $a \mapsto \|[D, a]\|$ is stable because~$D$ is self-adjoint.
Thus the form of the right-hand side of the above inequality is reasonable.
Notice further that because the representation of~${\mathcal A}$ on $L^2({\mathcal A}, \tau)$ is faithful, we have
\begin{gather*}
\sup \{\|({\partial} a)b\|_\tau^2; \|b\|_\tau \leq 1\}=\sup \{\tau(b^*\langle{\partial} a, {\partial}
a\rangle_{\mathcal A} b): \|b\|_\tau \leq 1\}=\|\langle{\partial} a, {\partial} a\rangle_{\mathcal A}\|_\infty.
\end{gather*}
Consequently $\|[D, a]\| \geq \|\langle{\partial} a, {\partial} a\rangle_{\mathcal A}\|_\infty^{1/2}$.
But because $\|[D, a]\|$ is stable under the involution, we also have $\|[D, a]\| \geq \|\langle{\partial} a^*,
{\partial} a^*\rangle_{\mathcal A}\|_\infty^{1/2}$.
We have thus obtained:

\begin{Theorem}
With notation as above we have
\begin{gather*}
\|[D, a]\|=\max \big\{\|\langle{\partial} a, {\partial} a\rangle_{\mathcal A}\|_\infty^{1/2},\|\langle{\partial} a^*,
{\partial} a^*\rangle_{\mathcal A}\|_\infty^{1/2} \big\}.
\end{gather*}
\end{Theorem}

To see that this can not be improved we have:
\begin{Example}
As in Example~\ref{exnc}, choose $v \in {\mathcal A}$ such that $[v, v^*] \neq 0$, and set ${\partial}_v a~=[v, a]$ for
all $a \in {\mathcal A}$.
Then for $a=v$ we have $\langle{\partial} a, {\partial} a\rangle_{\mathcal A}=0$ whereas $\langle{\partial} a^*,
{\partial} a^*\rangle_{\mathcal A}=[v, v^*]^*[v, v^*] \neq 0$.
\end{Example}

\begin{Example}
\label{exchd}
When ${\mathcal A}=C(X)$ and its Riemannian metric is determined by the function~$c$ on~$Z$, then Theorem~\ref{thmch}
shows that for its Hodge--Dirac operator~$D$ the corresponding seminorm,~$L$, is given by
\begin{gather*}
L(f)=\|[D, f]\|=\|\langle{\partial} f, {\partial} f\rangle_{\mathcal A}\|_\infty^{1/2}
=\sup_y\left\{\left(\sum\limits_x|f(x)-f(y)|^2c_{xy}\right)^{1/2}\right\}.
\end{gather*}
Aside from a~traditional factor of 1/2 this is the same seminorm as the seminorms $d_3$ and $d_4$ def\/ined after Lemma~4.1 of of~\cite{Dav}.
It is easily seen to be Markov (and Leibniz) in slight generalization of the seminorms def\/ined after the proof of Theorem~\ref{thmsemi}.
\end{Example}

We now show that the above seminorm seldom is the energy seminorm for a~Riemannian metric on~${\mathcal A}$.

\begin{Theorem}
\label{thmlink}
Let ${\mathcal A}=C(X)$ and let its Riemannian metric be determined by the function~$c$ on~$Z$ as in Theorem~{\rm \ref{thmch}}.
Assume that $c_{xy}=c_{yx}$ for all $x, y \in X$ and that~$X$ is metrically connected.
Then the seminorm $L(f)=\|[D, f]\|$ for its Hodge--Dirac operator can be obtained as the energy seminorm from
a~Riemannian metric on~${\mathcal A}$ if and only if there is a~point $t \in X$ such that every other point in~$X$ is
linked only to~$t$, that is, if $c_{xy} \neq 0$ exactly when $x=t$ or $y=t$.
\end{Theorem}

\begin{proof}
Suppose that the above~$L$ for~$D$ can be obtained from an energy form.
Then~$L$ must satisfy the parallelogram law.
Let~$p$ and~$q$ be any two distinct points of~$X$, and set $f=\delta_p$ and $g=\delta_q$, so that $f+g=\delta_p+\delta_q$ and $f-g=\delta_p-\delta_q$.
Let us denote $(L(f))^2$ simply by $L^2(f)$, etc.
Then it is easily calculated that $L^2(f)=\hat c(p)$ and similarly for $L^2(g)$, where $\hat c$ was def\/ined before
Theorem~\ref{thmneg} by $\hat c(p)=\sum\limits_x c_{xp}$.
For any distinct $u, v \in X$ def\/ine $m(u, v)$ by
\begin{gather*}
m(u, v)=\sup_{x \neq u, v} (c_{xu}+c_{xv}).
\end{gather*}
In terms of~$m$ one can calculate that
\begin{gather*}
L^2(f+g)=[\hat c(p)\vee \hat c(q)-c_{pq}]\vee m(p,q),
\end{gather*}
and that
\begin{gather*}
L^2(f-g)=[\hat c(p)\vee \hat c(q)+3c_{pq}]\vee m(p,q),
\end{gather*}
where the $\vee$ means ``maximum''.
Thus if~$L$ satisf\/ies the parallelogram law, then in particular we must have
\begin{gather}
[\hat c(p)\vee \hat c(q)-c_{pq}]\vee m(p,q)+[\hat c(p)\vee \hat c(q)+3c_{pq}]\vee m(p,q)=2(\hat c(p)+\hat c(q)).
\label{paral}
\end{gather}

Now, choose $t \in X$ such that $\hat c(t) \geq \hat c(x)$ for all $x \in X$.
We will show that every $x \in X$ is linked exactly to~$t$.

As a~f\/irst step, we show that every element of~$X$ is linked to~$t$, that is, $c_{tq} \neq 0$ for every $q \in X$ with
$q \neq t$.
So given~$q$, suppose that $c_{tq}=0$.
Set $p=t$ in formula~\eqref{paral}.
Then that formula becomes
\begin{gather*}
\hat c(t)\vee m(t,q)=\hat c(t)+\hat c(q).
\end{gather*}
If $\hat c(t) \geq m(t,q)$ then we obtain $\hat c(q)=0$ which contradicts connectedness.
Otherwise, by the def\/inition of $m(t,q)$ there is an $r \in X$ distinct from~$t$ and~$q$ such that
\begin{gather*}
c_{rt}+c_{rq}=\hat c(t)+\hat c(q).
\end{gather*}
It follows that $\hat c(t)=c_{rt}$ and $\hat c(q)=c_{rq}$, from which we see that, because $\hat c(q) \neq 0$ by
connectedness, we have
\begin{gather*}
\hat c(r) \geq c_{rt}+c_{rq} > \hat c(t),
\end{gather*}
which contradicts the maximality of $\hat c(t)$.
Consequently we must have $c_{tq} \neq 0$.

As the f\/inal step, we show that the elements of~$X$ are only linked to~$t$, that is, for every $p \in X$ with $p \neq t$
we have $c_{xp}=0$ for all $x \neq t$.
So, let~$p$ be given, with $p \neq t$.
Choose $q \in X$ distinct from~$p$ such that $c_{pq} \geq c_{px}$ for all $x \in X$ (so possibly $q=t$).
Then for any~$x$ we have $\hat c(q) \geq c_{xq}$ and so $[\hat c(p)\vee \hat c(q)+3c_{pq}] \geq c_{xp}+c_{xq}$.
It follows that
\begin{gather*}
[\hat c(p)\vee \hat c(q)+3c_{pq}]\geq m(p,q).
\end{gather*}
Thus formula~\eqref{paral} becomes
\begin{gather}
[\hat c(p)\vee \hat c(q)-c_{pq}]\vee m(p,q)+\hat c(p)\vee \hat c(q)+3c_{pq}=2(\hat c(p)+\hat
c(q)).
\label{wedge}
\end{gather}

Case 1.
Suppose that $[\hat c(p)\vee \hat c(q)-c_{pq}]\geq m(p,q)$, so that formula~\eqref{wedge} reduces to
\begin{gather*}
\hat c(p)\vee \hat c(q)+c_{pq}=\hat c(p)+\hat c(q).
\end{gather*}
If $\hat c(p) \geq \hat c(q)$ then we f\/ind that $\hat c(q)=c_{pq}$, so that~$q$ is linked only to~$p$.
But by the f\/irst step above,~$q$ must be linked to~$t$.
Thus we must have $q=t$, and so~$t$ is the only element to which~$p$ is linked, as desired.
If instead we have $\hat c(p) \leq \hat c(q)$ then we f\/ind that $\hat c(p)=c_{pq}$, and in a~similar way this implies
that~$p$ is only linked to~$t$.

Case 2.
Suppose instead that $[\hat c(p)\vee \hat c(q)-c_{pq}] \leq m(p,q)$, so that formula~\eqref{wedge} reduces to
\begin{gather*}
m(p,q)+\hat c(p)\vee \hat c(q)+3c_{pq}=2(\hat c(p)+\hat c(q)).
\end{gather*}
Then by the def\/inition of $m(p,q)$ there is an $r \in X$ distinct from~$p$ and~$q$ such that
\begin{gather*}
c_{pr}+c_{qr}+\hat c(p)\vee \hat c(q)+3c_{pq}=2(\hat c(p)+\hat c(q)).
\end{gather*}
If $\hat c(p) \geq \hat c(q)$ then we f\/ind that
\begin{gather*}
c_{pr}+c_{qr}=(\hat c(p)-c_{pq})+2(\hat c(q)-c_{pq}).
\end{gather*}
Because~$r$ is distinct from~$p$ and~$q$, we must have $\hat c(q)=c_{pq}$.
So~$q$ is linked only to~$p$.
But by step 1 above~$q$ is linked to~$t$.
This contradicts the assumption that $p \neq t$.
Thus we must have $\hat c(p) \leq \hat c(q)$, in which case we f\/ind that
\begin{gather*}
c_{pr}+c_{qr}=2(\hat c(p)-c_{pq})+(\hat c(q)-c_{pq}).
\end{gather*}
Thus, much as above, we must have $\hat c(p)=c_{pq}$, so that~$p$ is linked only to~$q$.
But by step 1 above~$p$ is linked to~$t$.
Thus we must have $q=t$, and so~$p$ is linked only to~$t$, as desired.

The converse assertion of the theorem is easily verif\/ied, and is closely related to Proposition~3.8 of~\cite{R28} and
the class of examples discussed in connection with standard deviation in Section~2 of~\cite{R28}, as we will discuss
again in Section~\ref{secsdev}.
\end{proof}

Returning to the general situation, we caution that if we start with a~spectral triple $({\mathcal A}, {\mathcal H}, D)$
for~${\mathcal A}$, and then form its corresponding Riemannian metric $(\Omega, {\partial},
\langle\cdot,\cdot\rangle_{\mathcal A})$, and then form the Hodge--Dirac operator $D_H$ for this spectral triple as
above, then usually we will have $D \neq D_H$.
The following commutative example is instructive.

\begin{Example}
Any f\/inite metric space $(X,\rho)$ has an essentially canonical spectral triple.
Let~$Z$ be def\/ined as in Example~\ref{excomm}, and set $c_{xy}=1/\rho(x,y)$ for all $(x,y) \in Z$.
Form ${\mathcal H}=L^2(Z)$ for counting measure, and let $M_c$ denote the operator on ${\mathcal H}$ of pointwise
multiplication by~$c$.
Let~$U$ be the self-adjoint unitary operator on ${\mathcal H}$ def\/ined by $(U\xi)(x,y)=\xi(y,x)$, and set $D=M_c U$.
Let~${\mathcal A}$ act on ${\mathcal H}$ by its left action on $C(Z)$ used in Example~\ref{excomm}, that is, by
$(f\xi)(x,y)=f(x) \xi(x,y)$.
Then $({\mathcal A}, {\mathcal H}, D)$ is a~spectral triple for~${\mathcal A}$.
Furthermore, it is easily seen that for every $f \in {\mathcal A}$ we have
\begin{gather*}
\|[D, f]|=L(f)=\sup\{|f(x)-f(y)|/\rho(x,y): (x,y) \in Z\},
\end{gather*}
which is the usual Lipschitz constant for~$f$.
The corresponding metric on the state space $S({\mathcal A})$ of~${\mathcal A}$, when restricted to~$X$ identif\/ied with
the extreme points of $S({\mathcal A})$, is exactly the original metric~$\rho$.

Let us now calculate the CdC for the Riemannian metric for the above spectral triple.
Notice that $M_c$ commutes with~$U$ and with the action of~${\mathcal A}$.
Notice further that $[U, f]=-M_{{\partial} f}$, so that
\begin{gather*}
[D, f]=-M_cM_{{\partial} f}
\end{gather*}
and
\begin{gather*}
[D, f]^*[D,g]=M_c^2M_{{\partial} \bar f}M_{{\partial} g}
\end{gather*}
for $f, g \in {\mathcal A}$.
According to Example~\ref{exspec} we must apply to this latter the conditional expectation~$E$ from ${\mathcal
B}({\mathcal H})$ onto~${\mathcal A}$ corresponding to the trace on ${\mathcal B}({\mathcal H})$.
We see that it suf\/f\/ices to determine the restriction of this conditional expectation to $C(Z)$, where $C(Z)$ is viewed
as an algebra of pointwise multiplication operators on ${\mathcal B}({\mathcal H})$.
Since~${\mathcal A}$ can be viewed as consisting of functions in $C(Z)$ that depend only on the f\/irst coordinate,
averaging over the second coordinate gives a~conditional expectation from $C(Z)$ onto ${\mathcal A} \subseteq C(Z)$.
By examining the proof of Proposition~2.36 of~\cite{Tks}
it is not dif\/f\/icult to see that this is the restriction of the
conditional expectation~$E$ from ${\mathcal B}({\mathcal H})$ onto~${\mathcal A}$.
Thus for any $F \in C(Z)$ we have
\begin{gather*}
E(F)(x)=(n-1)^{-1}\sum\limits_y \{F(x,y): (x,y) \in Z\},
\end{gather*}
where~$n$ is the number of elements of~$X$.
From this we f\/ind that the CdC is
\begin{gather*}
\Gamma(f,g)(x)=\langle{\partial} f, {\partial} g\rangle_{\mathcal A}(x)=(n-1)^{-1}\sum\limits_y (\bar f(x)-\bar
f(y))(g(x)-g(y))c^2_{xy}.
\end{gather*}
Notice that up to the constant in front, this is the expression used to determine the CdC in~\ref{thmch} except using
$c^2$ instead of~$c$.
From Example~\ref{exchd} we see that the seminorm~$L$ corresponding to the Dirac operator for this CdC is given by
\begin{gather*}
L(f)=\|\langle{\partial} f, {\partial} f\rangle_{\mathcal A}\|_\infty^{1/2}=(n-1)^{-1/2}
\sup_x\left\{\left(\sum\limits_y|f(x)-f(y)|^2 c^2_{xy}\right)^{1/2}\right\}.
\end{gather*}
This is quite dif\/ferent from the seminorm~$L$ that we started with near the beginning of this example, and this shows
that the Dirac operators themselves are quite dif\/ferent.
\end{Example}

\section{Quotients of energy metrics}
\label{secquot}

Quotients of energy forms are discussed in the literature for the case in which~${\mathcal A}$ is commutative, e.g.\
on p.~44 of~\cite{Kgm}.
But I have not seen any discussion of quotients for non-commutative~${\mathcal A}$.
We give such a~discussion here, since quotients are important for the theory of quantum metric spaces.
See Section~5 of~\cite{R21}.

Let ${\mathcal E}$ be the energy form for a~Riemannian metric and trace on~${\mathcal A}$, and let ${\mathcal B}$ be
a~quotient $C^*$-algebra of~${\mathcal A}$, with~$\pi$ the quotient homomorphism from~${\mathcal A}$ onto ${\mathcal
B}$.
It is not so clear how we should def\/ine the quotient of ${\mathcal E}$ on ${\mathcal B}$.
But we can consider the corresponding energy norm, $L_{\mathcal E}$, on~${\mathcal A}$, and we do know how to take
quotients of (semi)-norms.
The quotient, $L_{\mathcal E}^{\mathcal B}$, of $L_{\mathcal E}$ is, of course, def\/ined by
\begin{gather*}
L_{\mathcal E}^{\mathcal B}(b)=\inf \{L_{\mathcal E}(a):a \in {\mathcal A}
\;
\text{and}
\;
\pi(a)=b\}
\end{gather*}
for all $b \in {\mathcal B}$.
The following general observation is an important step in our discussion.

\begin{Proposition}
\label{proqmk}
Let~$L$ be a~seminorm on a~$C^*$-algebra~${\mathcal A}$, and let ${\mathcal B}$ be a~quotient $C^*$-algebra of
${\mathcal A}$.
Let $L^{\mathcal B}$ denote the quotient seminorm on ${\mathcal B}$.
If~$L$ is Markov then so is~$L$.
\end{Proposition}

\begin{proof}
Let $b \in {\mathcal B}$ with $b^*=b$, and let~$F$ be a~Lipschitz function from $\sigma(b)$ to ${\mathbb R}$.
Let $\hat F$ be an extension of~$F$ to all of ${\mathbb R}$ such that Lip($\hat F$)=Lip($F$)
(see Theorem~1.5.6 of~\cite{Wvr2}). Let $\varepsilon > 0$ be given.
Then there exists an $a \in {\mathcal A}$ with $a^*=a$ such that $\pi(a)=b$ and $L(a) \leq \tilde L(b)+\varepsilon$.
Note that $\pi(\hat F(a))=\hat F(b)=F(b)$.
Then
\begin{gather*}
L^{\mathcal B}(F(b)) \leq L(\hat F(a)) \leq \Lip(\hat F)L(a) \leq \Lip(F)\big(L^{\mathcal B}(b)+\varepsilon\big).
\end{gather*}
Since~$\varepsilon$ is arbitrary, we obtain $L^{\mathcal B}(F(b)) \leq \Lip(F)(L^{\mathcal B}(b))$ as desired.
\end{proof}

Now in our f\/inite-dimensional setting every 2-sided ideal of~${\mathcal A}$ is generated by a~central projection, and
the quotient by that ideal can be identif\/ied with the sub-$C^*$-algebra generated by the complementary central
projection.
Thus there is a~(proper) central projection,~$p$, in~${\mathcal A}$ such that ${\mathcal B}$ can be identif\/ied with
$p{\mathcal A}$.
We write ${\mathcal B}=p{\mathcal A}$.
The quotient map from~${\mathcal A}$ onto ${\mathcal B}$ is then simply given by $\pi(a)=pa$.

To understand $L_{\mathcal E}^{\mathcal B}$ more clearly we now express it in terms of the Laplace operator~$\Delta$ for~${\mathcal E}$.
We are primarily interested in the corresponding metric on the state space, and so to avoid unimportant complications we
will treat here only the case in which~${\mathcal A}$ is metrically connected (Def\/inition~\ref{defcon}).
Thus we assume that if $L_{\mathcal E} (a)=0$ then $a \in {\mathbb C} 1_{\mathcal A}$, so that the kernel of~$\Delta$
is exactly ${\mathbb C} 1_{\mathcal A}$.
We also require that ${\mathcal E}$ is real, as def\/ined in Def\/inition~\ref{defereal}.

We follow the argument that is given for the commutative case around Lemma~2.1.5 of~\cite{Kgm}.
Many of the calculations below will work for any positive operator on $L^2({\mathcal A}, \tau)$ whose kernel is exactly
${\mathbb C} 1_{\mathcal A}$.
Let ${\mathcal C}=(1-p){\mathcal A}$, so that ${\mathcal A}={\mathcal B} \oplus {\mathcal C}$ as $C^*$-algebras.
Let~$\tau$ also denote its restrictions to ${\mathcal B}$ and ${\mathcal C}$,
so that $L^2({\mathcal A}, \tau)=L^2({\mathcal B}, \tau) \oplus L^2({\mathcal C}, \tau)$, an orthogonal decomposition for the~$\tau$-inner-product.
(But note that if~$\tau$ happens to be normalized, the traces obtained by restricting~$\tau$ to ${\mathcal B}$ or
${\mathcal C}$ will not be normalized, but this is not a~dif\/f\/iculty.) With respect to this decomposition~$\Delta$ can be
expressed as a~matrix:
\begin{gather*}
\Delta=\begin{pmatrix}
R & J^*
\\
J & S
\end{pmatrix}
,
\end{gather*}
in which $R \geq 0$ and $S\geq 0$.
Because~${\mathcal A}$ is metrically connected,~$S$ is invertible as an operator on $L^2({\mathcal C}, \tau)$.
(To see this, note that if~$c$ is in the kernel of~$S$, and if we set $a=0_{\mathcal B} \oplus c$, then ${\mathcal
E}(a,a)=\langle a, \Delta a\rangle_\tau=0$, so that $a \in {\mathbb C} 1_{\mathcal A}$, so that $c=0$.) Then we can
use~$S$ to do ``row and column operations'' on the matrix for~$\Delta$ to obtain its Schur complement.
Specif\/ically, we have
\begin{gather*}
\begin{pmatrix}
R & J^*
\\
J & S
\end{pmatrix}
=\begin{pmatrix}
I & J^*S^{-1}
\\
0 & I
\end{pmatrix}
\begin{pmatrix}
R- J^*S^{-1}J & 0
\\
0 & S
\end{pmatrix}
\begin{pmatrix}
I & 0
\\
S^{-1}J & I
\end{pmatrix}
.
\end{gather*}
Then for every $a=b \oplus c$ we have
\begin{gather*}
{\mathcal E}(a,a)=\Big\langle
\begin{pmatrix}
b
\\
c
\end{pmatrix}
,
\begin{pmatrix}
R & J^*
\\
J & S
\end{pmatrix}
\begin{pmatrix}
b
\\
c
\end{pmatrix}
\Big\rangle_\tau
\\
\phantom{{\mathcal E}(a,a)}
=\Big\langle
\begin{pmatrix}
I & 0
\\
S^{-1}J & I
\end{pmatrix}
\begin{pmatrix}
b
\\
c
\end{pmatrix}
,
\begin{pmatrix}
R-J^*S^{-1}J & 0
\\
0 & S
\end{pmatrix}
\begin{pmatrix}
I & 0
\\
S^{-1}J & I
\end{pmatrix}
\begin{pmatrix}
b
\\
c
\end{pmatrix}
\Big\rangle_\tau
\\
\phantom{{\mathcal E}(a,a)}
=\langle b, (R -J^*S^{-1}J)b\rangle_\tau+\langle S^{-1}Jb+c,S(S^{-1}Jb+c)\rangle_\tau.
\end{gather*}
Then for f\/ixed~$b$ it is clear that ${\mathcal E}(a,a)$ is minimized by setting $c=-S^{-1}Jb$, and that the minimum is
$\langle b, (R -J^*S^{-1}J)b\rangle_\tau$.
Accordingly, let us set
\begin{gather*}
\Delta^{\mathcal B}=R -J^*S^{-1}J,
\end{gather*}
which is the Schur complement for~$S$.
From the above computations we see that $\Delta^{\mathcal B} \geq 0$ as a~Hilbert-space operator.
It is thus natural to def\/ine the quotient, ${\mathcal E}^{\mathcal B}$, of ${\mathcal E}$ by
\begin{gather*}
{\mathcal E}^{\mathcal B}(b, b')=\langle b, \Delta^{\mathcal B}(b')\rangle_\tau
\end{gather*}
for all $b, b' \in {\mathcal B}$,

Notice that $1_{\mathcal B}=p$.
Because $1_{\mathcal A}=1_{\mathcal B} \oplus 1_{\mathcal C}$, it follows from the above computations that
$\Delta^{\mathcal B}(1_{\mathcal B})=0$.
Furthermore, if $\Delta^{\mathcal B}(b)=0$, then if we set $a=b \oplus (-S^{-1}Jb)$ we see that ${\mathcal E}(a,a)=0$,
so that $a \in {\mathbb C} 1_{\mathcal A}$ and so $b \in {\mathbb C} 1_{\mathcal B}$.
Thus ${\mathcal B}$ is metrically connected.

\begin{Theorem}
\label{thmquot}
Let ${\mathcal E}$ be the energy form for a~metrically connected Riemannian metric on~${\mathcal A}$ and a~faithful
trace on~${\mathcal A}$.
Assume further that ${\mathcal E}$ is real~$($Definition~{\rm \ref{defereal})}.
Let ${\mathcal B}$ be a~quotient $C^*$-algebra of~${\mathcal A}$, and define ${\mathcal E}^{\mathcal B}$ as above.
Then ${\mathcal E}^{\mathcal B}$ is the energy form for a~metrically connected Riemannian metric on ${\mathcal B}$.
\end{Theorem}

\begin{proof}
According to Theorem~\ref{thmdir} it suf\/f\/ices to show that ${\mathcal E}^{\mathcal B}$ is a~real completely positive and
completely Markov form.

Because ${\mathcal E}$ is assumed to be real, its energy norm $L_{\mathcal E}$ is a~$*$-seminorm.
It is easily seen that then its quotient seminorm is a~$*$-seminorm, so that
\begin{gather*}
{\mathcal E}^{\mathcal B}(b^*,b^*)={\mathcal E}^{\mathcal B}(b,b)
\end{gather*}
for all $b \in {\mathcal B}$.
By the usual polarization identity it follows that ${\mathcal E}_{\mathcal B}$ is real.

For each natural number~$n$ we must show that the form $({\mathcal E}^{\mathcal B})_n$ is positive and Markov, where
$({\mathcal E}^{\mathcal B})_n$ is def\/ined on $M_n({\mathcal B})$ by
\begin{gather*}
\big({\mathcal E}^{\mathcal B}\big)_n(B,B')=\sum\limits_{jk} {\mathcal E}^{\mathcal B}(b_{jk}, b'_{jk})
\end{gather*}
for all $B,B' \in M_n({\mathcal B})$.
But, much as in the proof of Proposition~\ref{prolpn}, the right-hand side is equal to
\begin{gather*}
\sum\limits_{jk}\langle b_{jk}, \Delta^{\mathcal B} b'_{jk}\rangle_{\tau_{\mathcal B}}=\langle B, \big(I_n \otimes \Delta^{\mathcal
B}\big)B'\rangle_{\text{tr}_n \otimes \tau_{\mathcal B}}.
\end{gather*}
Now by Proposition~\ref{prolpn} the Laplacian for ${\mathcal E}_n$ is $I_n \otimes \Delta$.
Then it is easily seen that the matricial expression for this operator for
the decomposition $M_n({\mathcal A})=M_n({\mathcal B}) \oplus M_n({\mathcal C})$ is given by
\begin{gather*}
\Delta=\begin{pmatrix}
I_n\otimes R & I_n\otimes J^*
\\
I_n\otimes J & I_n\otimes S
\end{pmatrix}
.
\end{gather*}
Notice that $ I_n\otimes S$ is invertible.
Then we see that the Schur complement for $ I_n\otimes S$, which we denote by $(I_n\otimes \Delta)^{\mathcal B}$, is
\begin{gather*}
I_n\otimes R-(I_n\otimes J^*)(I_n\otimes S)^{-1}(I_n\otimes J)=I_n\otimes \Delta^{\mathcal B}.
\end{gather*}
Arguing much as above, and comparing with the expression obtained above for $({\mathcal E}^{\mathcal B})_n$, we thus
f\/ind that the quotient, $({\mathcal E}_n)^{\mathcal B}$, of ${\mathcal E}_n$ on $M_n({\mathcal B})$ is given by
$({\mathcal E}^{\mathcal B})_n$.
It follows that $({\mathcal E}^{\mathcal B})_n$ is positive, and from Proposition~\ref{proqmk} it follows that
$({\mathcal E}^{\mathcal B})_n$ is Markov.
Thus ${\mathcal E}^{\mathcal B}$ is completely positive and completely Markov, and so comes from a~CdC.
\end{proof}

For me this theorem is striking because it means that the class of Leibniz seminorms coming from Riemannian metrics that
are~$\tau$-real for a~trace has the property that the quotient of any seminorm in this class is again Leibniz, as will
be any quotient of the quotient, etc.
This is quite a~contrast with the dif\/f\/iculties I had with quotients of Leibniz seminorms not necessarily being again
Leibniz, as discussed, for example, in Section~5 of~\cite{R21}.

\section{The relationship with standard deviation}
\label{secsdev}

In the f\/irst two paragraphs of Section~2 of~\cite{R28}, for a~f\/inite set~$X$, a~Dirac operator is def\/ined on ${\mathcal
A}=C(X)$ whose corresponding seminorm is easily seen to be of the form
\begin{gather*}
L(f)=\left(\sum\limits_{x,~x \neq x_*} |f(x)-f(x_*)|^2 \beta_x\right)^{1/2},
\end{gather*}
where $x_*$ is a~special point in~$X$ and the $\beta_x$'s are strictly positive real numbers.
(For the $\alpha_x$'s of those two paragraphs in~\cite{R28} we have $b_x=|\alpha_x|^2$.) If we include a~factor of
1/2, this seminorm clearly corresponds to a~resistance network in which every point of~$x$ is connected only to $x_*$,
with the conductances given by $c_{xx_*}=\beta_x$ while $c_{xy}=0$ if $x \neq x_* \neq y$.
Notice that this is exactly the situation that was obtained in Theorem~\ref{thmlink}.
The normalization $\sum |\alpha_x|^2=1$ used in~\cite{R28} corresponds to $\hat c(x_*)=1$ while $\hat c(x)=c_{xx_*}$ for $x \neq x_*$.

In~\cite{R28} the quotient of the above seminorm~$L$ by the minimal ideal of $C(X)$ corresponding to $x_*$ is shown to
be given by standard deviations.
This result is also extended there to non-commutative $C^*$-algebras.
We now show that standard deviations, and the versions for non-commutative $C^*$-algebras, come from our Riemannian
metrics.

As the set-up we have a~f\/inite-dimensional $C^*$-algebra~${\mathcal A}$ corresponding to $C(X \setminus \{x_*\})$, and
we have a~faithful trace~$\tau$ on~${\mathcal A}$.
We also have a~positive element,~$p$, of~${\mathcal A}$ such that $\tau(p)=1$, corresponding to the function $\hat c$
restricted to the complement in~$X$ of $\{x_*\}$.
We assume that~$p$ is strictly positive, corresponding to the connectedness of the resistance network.
Then~$p$ determines a~faithful state~$\mu$ on~${\mathcal A}$ by $\mu(a)=\tau(pa)$.
We make the further quite strong requirement on~$p$ that it is in the center of~${\mathcal A}$.
Thus~$\mu$ is a~faithful tracial state.

Let ${\mathcal B}={\mathcal A} \oplus {\mathbb C}$, and let $\hat \tau$ be the extension of~$\tau$ to ${\mathcal B}$
that has value 1 on $1 \in {\mathbb C}$, so $\tilde \tau$ corresponds to counting-measure on~$X$ when~${\mathcal A}$ is
commutative.
We def\/ine an operator,~$\Delta$, on $L^2({\mathcal B}, \hat \tau)$ by
\begin{gather*}
\Delta=\begin{pmatrix}
M & J^*
\\
J & 1
\end{pmatrix}
\end{gather*}
for the decomposition ${\mathcal B}={\mathcal A} \oplus {\mathbb C}$, where~$M$ is def\/ined by $M(a)=pa$ and~$J$ is
def\/ined by $J(a)=-\mu(a)$, so that $J^*(\alpha)=-\alpha p$ for $\alpha \in {\mathbb C}$.
It is easily calculated that
\begin{gather*}
\langle(a, \alpha), \Delta(b, \beta)\rangle_{\hat \tau}=\mu((a-\alpha)^*(b-\beta)),
\end{gather*}
which shows that~$\Delta$ is a~positive operator on $L^2({\mathcal B}, \hat \tau)$, and also that $\Delta((1_{\mathcal
A}, 1))=0$.
A~straight-forward calculation shows that the qCdC for~$\Delta$ should be given by
\begin{gather*}
2\Gamma_\Delta((a, \alpha), (b, \beta))=(p(a- \alpha)^*(b- \beta), \mu((a- \alpha)^*(b- \beta))).
\end{gather*}
From these formulas it is not dif\/f\/icult to guess the form of a~Riemannian metric that leads to~$\Delta$.
To obtain it we proceed as follows.

Def\/ine a~${\mathcal B}$-bimodule by $\tilde \Omega={\mathcal A} \oplus {\mathcal A}$ with left and right actions of
${\mathcal B}$ given by
\begin{gather*}
(a, \alpha)(b, c)=(ab, \alpha c)
\qquad
\text{and}
\qquad
(b, c)(a, \alpha)=(\alpha b, ca)
\end{gather*}
for all $(a, \alpha) \in {\mathcal B}$ and all $(b,c) \in \tilde \Omega$.
Def\/ine a~${\mathcal B}$-valued inner product on $\tilde \Omega$ by
\begin{gather*}
\langle(a,b), (c,d)\rangle_{\mathcal B}=(1/2)(b^*dp, \tau(a^*cp))=(1/2)(b^*dp, \mu(a^*c)).
\end{gather*}
It is easily checked that with this inner product $\tilde \Omega$ is a~right Hilbert ${\mathcal B}$-module, and in fact
a~correspondence over ${\mathcal B}$.
We def\/ine a~derivation from ${\mathcal B}$ into $\tilde \Omega$ by
\begin{gather*}
{\partial}((a, \alpha))=(a-\alpha, -a+\alpha),
\end{gather*}
where of course here~$\alpha$ means $\alpha 1_{\mathcal A}$.
We note that if ${\partial}((a, \alpha))=0$, then $(a, \alpha)=\alpha(1_{\mathcal A}, 1)$.
Let~$\Omega$ be the sub-bimodule of $\tilde \Omega$ generated by the range of ${\partial}$.
Then we see that $(\Omega, \langle\cdot, \cdot\rangle_{\mathcal B}, {\partial})$ is a~Riemannian metric for ${\mathcal
B}$, for which ${\mathcal B}$ is metrically connected.
Let~$\Gamma$ be its CdC.
A~simple calculation shows that $\Gamma=\Gamma_\Delta$ for the $\Gamma_\Delta$ def\/ined in the previous paragraph, as
desired.
It is easy to check that~$\Gamma$ is $\tilde \tau$-real, and that~$\Delta$ is the Laplace operator for~$\Gamma$ and
$\tilde \tau$.
The energy form for~$\Gamma$ and $\tilde \tau$ is clearly given by
\begin{gather*}
{\mathcal E}_\Gamma((a, \alpha), (b, \beta))=\mu((a-\alpha)^*(b-\beta)).
\end{gather*}

We can now consider the quotient of the energy form when we factor ${\mathcal B}$ by its ideal ${\mathbb C}$.
The quotient can be identif\/ied in the evident way with~${\mathcal A}$.
We def\/ined~$\Delta$ by its matrix for the decomposition ${\mathcal B}={\mathcal A} \oplus {\mathbb C}$ of ${\mathcal
B}$.
Thus we are already in position to apply the discussion leading up to Theorem~\ref{thmquot} to obtain the Laplace
operator $\Delta^{\mathcal A}$ for the quotient~${\mathcal A}$.
(For the case in which~${\mathcal A}$ is commutative, this is closely related to the second half of Remark~4.40
of~\cite{JrP2}.) For our notation above in which~$M$ plays the role of~$R$ in Section~\ref{secquot}, we f\/ind that
\begin{gather*}
\Delta^{\mathcal A}(a)=M(a)-J^*J(a)=p(a-\mu(a)).
\end{gather*}
The corresponding CdC is given by
\begin{gather*}
\Gamma_{\Delta^{\mathcal A}}(a, b)=(1/2)p(a^*b-\mu(a^*)b-a^*\mu(b)+\mu(a^*b)),
\end{gather*}
and the corresponding energy form is
\begin{gather*}
{\mathcal E}^{\mathcal A}(a, b)=\langle a-\mu(a), b-\mu(b)\rangle_\mu=\mu(a^*b)-\mu(a^*)\mu(b).
\end{gather*}
The corresponding energy seminorm is
\begin{gather*}
L^{\mathcal A}(a)=\|a-\mu(a)\|_\mu,
\end{gather*}
which, when $a^*=a$, is exactly the standard deviation of~$a$ for the state~$\mu$ as discussed in~\cite{R28} and in
the quantum physics literature.

It is not hard to guess a~specif\/ic description of the Riemannian metric for this situation.
Def\/ine a~``slice map''~$E$ from ${\mathcal A} \otimes {\mathcal A}$ onto $1_{\mathcal A} \otimes {\mathcal A}={\mathcal A}$
by $E(a \otimes b)=\tau(a)b$.
It is easy to see that, up to normalization of~$\tau$, this is a~faithful conditional expectation.
Much as in Examples~\ref{exexp} and~\ref{exmulti} we def\/ine a~corresponding~${\mathcal A}$-valued inner product on
${\mathcal A} \otimes {\mathcal A}$, determined on elementary tensors by
\begin{gather*}
\langle a \otimes b, c \otimes d\rangle_{\mathcal A}=(1/2)E(a^*c \otimes b^*d)=(1/2)b^*\tau(a^*c)d.
\end{gather*}
With this inner product ${\mathcal A} \otimes {\mathcal A}$ becomes a~correspondence over~${\mathcal A}$.
We def\/ine ${\partial}$ somewhat as in Theorem~\ref{thmcst} by
\begin{gather*}
{\partial}(a)=p(a \otimes 1_{\mathcal A}-1_{\mathcal A} \otimes a)p.
\end{gather*}
Because~$p$ in invertible, the sub-bimodule generated by the range of ${\partial}$ is, as in the proof of
Proposition~\ref{progcal}, seen to be the kernel of the bimodule homomorphism from ${\mathcal A} \otimes {\mathcal A}$
onto~${\mathcal A}$ sending $a \otimes b$ to $ab$.
We denote this bimodule by $\Omega$.
Then $(\Omega, \langle\cdot, \cdot \rangle_{\mathcal A}, {\partial})$ is a~Riemannian metric for~${\mathcal A}$.
Its CdC is given by
\begin{gather*}
\Gamma(a,b)=(1/2)\langle p(a\otimes 1_{\mathcal A}-1_{\mathcal A}\otimes a)p,p(b \otimes 1_{\mathcal A}-1_{\mathcal A}\otimes b)p\rangle_{\mathcal A}
\\
\phantom{\Gamma(a,b)}
=(1/2)(\mu(a^*b)-\mu(a^*)b-a^*\mu(b)+a^*b)p.
\end{gather*}
This agrees with the CdC obtained just above.
All of this is related to the ``independent copies trick'' discussed before Proposition~3.6 of~\cite{R28}.

We remark that it does not seem easy to guess this Riemannian metric directly from that preceding it for ${\mathcal B}$,
of which~${\mathcal A}$ is the quotient, without using the Laplace operators.
This seems to be related to the fact that so far there does not seem to be known a~useful general way to obtain from
a~spectral triple on a~$C^*$-algebra a~spectral triple on a~quotient $C^*$-algebra of that $C^*$-algebra.

We also remark that in~\cite{R28} the case in which~$\mu$ is not tracial is treated, but we do not pursue that aspect here.

\subsection*{Acknowledgements}

The research reported here was supported in part by National Science Foundation grant DMS-1066368.

\LastPageEnding

\end{document}